\documentclass[11pt]{article}
\usepackage[T1]{fontenc}
\usepackage[utf8]{inputenc}
\usepackage{lmodern,amsmath,amssymb,amsthm,mathtools}
\usepackage[a4paper,margin=2.5cm]{geometry}
\usepackage{longtable,booktabs,array,calc,graphicx,xcolor,newunicodechar,microtype}
\usepackage{fancyvrb}
\usepackage{xurl}
\usepackage[hidelinks]{hyperref}
\providecommand{\tightlist}{\setlength{\itemsep}{0pt}\setlength{\parskip}{0pt}}

\def\textsc#1{#1}
\newunicodechar{π}{\ensuremath{\pi}}\newunicodechar{γ}{\ensuremath{\gamma}}\newunicodechar{δ}{\ensuremath{\delta}}
\newunicodechar{ε}{\ensuremath{\varepsilon}}\newunicodechar{α}{\ensuremath{\alpha}}\newunicodechar{β}{\ensuremath{\beta}}
\newunicodechar{θ}{\ensuremath{\theta}}\newunicodechar{κ}{\ensuremath{\kappa}}\newunicodechar{λ}{\ensuremath{\lambda}}
\newunicodechar{μ}{\ensuremath{\mu}}\newunicodechar{ν}{\ensuremath{\nu}}\newunicodechar{ρ}{\ensuremath{\rho}}
\newunicodechar{τ}{\ensuremath{\tau}}\newunicodechar{Δ}{\ensuremath{\Delta}}\newunicodechar{Σ}{\ensuremath{\Sigma}}
\newunicodechar{ℓ}{\ensuremath{\ell}}\newunicodechar{ℚ}{\ensuremath{\mathbb{Q}}}\newunicodechar{ℝ}{\ensuremath{\mathbb{R}}}
\newunicodechar{−}{\ensuremath{-}}\newunicodechar{∓}{\ensuremath{\mp}}\newunicodechar{±}{\ensuremath{\pm}}
\newunicodechar{·}{\ensuremath{\cdot}}\newunicodechar{×}{\ensuremath{\times}}\newunicodechar{√}{\ensuremath{\surd}}
\newunicodechar{→}{\ensuremath{\to}}\newunicodechar{↦}{\ensuremath{\mapsto}}\newunicodechar{∅}{\ensuremath{\emptyset}}
\newunicodechar{∈}{\ensuremath{\in}}\newunicodechar{∉}{\ensuremath{\notin}}\newunicodechar{∎}{\ensuremath{\blacksquare}}
\newunicodechar{∏}{\ensuremath{\prod}}\newunicodechar{∖}{\ensuremath{\setminus}}\newunicodechar{∞}{\ensuremath{\infty}}
\newunicodechar{∩}{\ensuremath{\cap}}\newunicodechar{∪}{\ensuremath{\cup}}\newunicodechar{≈}{\ensuremath{\approx}}
\newunicodechar{≠}{\ensuremath{\neq}}\newunicodechar{≤}{\ensuremath{\le}}\newunicodechar{≥}{\ensuremath{\ge}}
\newunicodechar{⊂}{\ensuremath{\subset}}\newunicodechar{⊆}{\ensuremath{\subseteq}}\newunicodechar{⊇}{\ensuremath{\supseteq}}
\newunicodechar{⊔}{\ensuremath{\sqcup}}\newunicodechar{⌊}{\ensuremath{\lfloor}}\newunicodechar{⌋}{\ensuremath{\rfloor}}
\newunicodechar{′}{\ensuremath{'}}\newunicodechar{½}{\ensuremath{\tfrac12}}\newunicodechar{̂}{\^{}}
\newunicodechar{⁰}{\ensuremath{^{0}}}\newunicodechar{¹}{\ensuremath{^{1}}}\newunicodechar{²}{\ensuremath{^{2}}}\newunicodechar{³}{\ensuremath{^{3}}}\newunicodechar{⁴}{\ensuremath{^{4}}}\newunicodechar{⁵}{\ensuremath{^{5}}}\newunicodechar{⁶}{\ensuremath{^{6}}}\newunicodechar{⁷}{\ensuremath{^{7}}}\newunicodechar{⁸}{\ensuremath{^{8}}}\newunicodechar{⁹}{\ensuremath{^{9}}}\newunicodechar{⁺}{\ensuremath{^{+}}}\newunicodechar{⁻}{\ensuremath{^{-}}}\newunicodechar{ᵏ}{\ensuremath{^{k}}}
\newunicodechar{∨}{\ensuremath{\vee}}\newunicodechar{∧}{\ensuremath{\wedge}}\newunicodechar{⊗}{\ensuremath{\otimes}}\newunicodechar{Λ}{\ensuremath{\Lambda}}\newunicodechar{⇒}{\ensuremath{\Rightarrow}}
\newunicodechar{₀}{\ensuremath{_{0}}}\newunicodechar{₁}{\ensuremath{_{1}}}\newunicodechar{₂}{\ensuremath{_{2}}}\newunicodechar{₃}{\ensuremath{_{3}}}\newunicodechar{₊}{\ensuremath{_{+}}}\newunicodechar{ₙ}{\ensuremath{_{n}}}\newunicodechar{ᵢ}{\ensuremath{_{i}}}
\newunicodechar{ }{\,\allowbreak}
\title{Certified pandigital approximations of $\pi$, $\gamma$, $\delta$ and $e$}
\author{Daniel Bendjema Garcia}
\date{20 November 2026}
\begin{document}
\maketitle

\begin{abstract}
In 2004 Erich Friedman asked for the best approximations of π, e, the Euler--Mascheroni constant γ and Feigenbaum's constant δ that use each of the digits 1, \ldots, n exactly once, with + − × / \^{}, parentheses and the decimal point. We determine the optimal expressions for π, γ and δ for every n ≤ 9, with the digits in any order and in increasing order, and for e when n ≤ 6. The optima are found by exhaustive search and certified by an independent algorithm that encloses every expression rigorously and decides each candidate with Arb and exact arithmetic; a second certified computation covers intermediate values of any size. Twenty-one entries improve the records published before this work; one of them, the optimal formula for π with 1--9 (error 2.81·\(10^{-14}\)), was found independently by M. Verheyen, and we prove it optimal. For e, whose records come from identities between power towers, we prove that Bamberger's formula is optimal within the natural family of such identities. All code and certificates are public.
\end{abstract}

\section*{1. Introduction}\addcontentsline{toc}{section}{1. Introduction}

\subsection*{1.1 The problem}\addcontentsline{toc}{subsection}{1.1 The problem}

Erich Friedman's Math Magic problem for August 2004 {[}1{]} reads:

\begin{quote}
This month's problem is to approximate famous mathematical constants using only the first n digits (each used exactly once) and the mathematical symbols + -- × / ( ) . and \^{}. For each n, what is the best approximation you can get to the following constants?
\end{quote}

The constants are γ, e, π and Feigenbaum's constant δ. The page records, for each constant and each n = 1, \ldots, 9, the best answer sent by readers, mostly by Bill Clagett, Gerrit de Blaauw, Richard Sabey and Friedman himself. A second set of tables answers Sabey's variant, ``What is the best approximation to these constants using the digits 1 through n in that order?''. The page also has variants that we do not treat (the digits of the constant itself, the digits 0--n, no decimal point, n copies of one digit). Two decades of contributions followed, some on Friedman's page and others elsewhere:

\begin{itemize}
\tightlist
\item
  \textbf{e.} Richard Sabey (2004) found \({\left(1+{9}^{-{4}^{6\cdot 7}}\right)}^{{3}^{{2}^{85}}}\), correct to 18 457 734 525 360 901 453 873 569 ≈ 1.8·\(10^{25}\) decimals (with the definition of Appendix B.1; N itself has \ldots570 digits) {[}1, 2, 3{]}. In March 2024 Daniel Bamberger raised this to 8 368 428 989 068 425 943 817 590 916 445 001 887 164 ≈ 8.4·\(10^{39}\) decimals with \({\left(1+{.2}^{{9}^{6\cdot 7}}\right)}^{{5}^{{3}^{84}}}\), in an answer on Mathematics Stack Exchange {[}16{]} that also explains how it was found {[}1, 2{]}. With the digits 0--9, the Reddit user ``Fastfaxr'' reached about 5.4·\(10^{45}\) decimals with \({({2}^{0}+{9}^{-{4}^{76}})}^{{3}^{{8}^{51}}}\) {[}2, 4{]}; in the same thread Bamberger points out that this solves a different problem, since Friedman's digits are 1 to n.~A. McKernon (2026) reports about 3.1·\(10^{30}\) decimals with 0--9 {[}11{]}; his formula uses 0! = 1, so it is outside Friedman's rules.
\item
  \textbf{π.} On Friedman's page the 1--9 entry is \({2}^{{5}^{.4}}-.6-{({.3}^{9}/7)}^{{.8}^{.1}}\), error 6.60·\(10^{-13}\), credited to Sabey; MathWorld credits the same formula to B. Ziv (2004) {[}1, 2{]}, and the Contest Center π page to Richard Hess {[}14{]}. Independently, Richard Hess posed the 1--9 problem for π in the Puzzle Corner of \emph{Technology Review} (July/August 2006), starting from \(3+(16-{8}^{-5})/(97+{2}^{4})\), error 3.3·\(10^{-9}\) {[}5{]}. The best answer printed in November/December 2006 was Joel Karnofsky's \(3.14+{\left({7}^{-{.9}^{-6}}+2/8\right)}^{5}\) (9.3·\(10^{-11}\)), and in July/August 2007 the column printed Sergey Ioffe's \(3+{5}^{-{7}^{.1}}+{.49}^{8}/{2}^{6}\) (3.8·\(10^{-13}\), found with a genetic algorithm) {[}5, 3{]}. Before November 2026, the best published 1--9 entry we found was Oleg Vlasii's (4.11·\(10^{-14}\), \((2/.98-.3)\cdot \left(.4+{5}^{{7}^{-.6}-.1}\right)\), from the Contest Center π page {[}3, 14{]}). On 8 November 2026 M. Verheyen independently communicated the optimal formula (2.81·\(10^{-14}\)) to MathWorld {[}2{]}. We recomputed every error quoted here with mpmath at 50 digits.
\end{itemize}

\subsection*{1.2 What is new here}\addcontentsline{toc}{subsection}{1.2 What is new here}

\begin{enumerate}
\def\labelenumi{\arabic{enumi}.}
\tightlist
\item
  \textbf{Optimality, not just records.} For π, γ and δ, every table entry for n ≤ 9, in both variants, is shown to be optimal, and the optimal value is unique. The optima are found by a fast exhaustive search in double precision and then proved optimal by an independent certifier with rigorous enclosures and exact verification (Theorem 1, Appendix A). A second certified computation extends the result to expressions with intermediate values out of range (Theorem 2, Appendix C).
\item
  \textbf{New records.} Twenty-one entries improve on the records published before this work (Table 1), including π with 1--9 (2.81·\(10^{-14}\) against 4.11·\(10^{-14}\); this formula was found independently by M. Verheyen {[}2{]} and is proved optimal here), γ with 1--9 (1.14·\(10^{-15}\) against 4.25·\(10^{-12}\)) and δ with 1--9 (1.38·\(10^{-14}\) against 3.04·\(10^{-10}\)).
\item
  \textbf{e.} For e the problem changes nature beyond n = 6: the best formulas use the limit \({(1+1/N)}^{N}\) with N a tower of powers, so they leave the range of Theorem 1. The certificate of item 1 covers e for n ≤ 6 (and n = 7 in order); beyond that we determine the best formulas inside the family F of identities (U + A\^{}(±\(B^{X}\)))\^{}(C\^{}(±\(D^{Y}\))) and (U + \(A^{±X}\))\^{}(\(C^{±Y}\)). Inside F, Bamberger's formula is optimal for 1--9, and new formulas are optimal for 1--7, 1--8, the in-order variants and 0--9.
\item
  \textbf{Variants.} We give the optima without unary minus (Section 8), with the same two certificates. We also explain why π and e behave so differently (Section 9).
\end{enumerate}

\subsection*{1.3 Structure}\addcontentsline{toc}{subsection}{1.3 Structure}

Section 2 fixes the rules. Section 3 gives the results. Sections 4--6 describe the search, the certifier and the treatment of out-of-range expressions. Section 7 treats e and Section 8 the variant without unary minus. Section 9 compares π and e, and Section 10 discusses limitations and open problems. The appendices contain the proofs.

\section*{2. Rules}\addcontentsline{toc}{section}{2. Rules}

Friedman's statement fixes the digits and the symbols. We make explicit the conventions that the statement takes for granted, so that the optimality claims have a precise meaning. Each convention follows from the statement or from the examples on Friedman's page (Section 2.3). Where Friedman fixes nothing, we checked that no reasonable choice changes any result (Section 2.4).

\subsection*{2.1 Expressions}\addcontentsline{toc}{subsection}{2.1 Expressions}

An \emph{expression} over a set S of digits is one of the following:

\begin{itemize}
\tightlist
\item
  \textbf{Literal:} the digits of S written consecutively in some order (in the in-order mode, increasing order), with at most one decimal point placed before any of them: 12, 1.2, .12. Its value is the decimal number read.
\item
  \textbf{−E}, where E is an expression over S (unary minus).
\item
  \textbf{(\(E_{1}\) op \(E_{2}\))}, where op ∈ \{+, −, ×, /, \^{}\}, \(E_{1}\) is over \(S_{1}\) and \(E_{2}\) over \(S_{2}\), and \(S_{1}\), \(S_{2}\) are disjoint, nonempty and \(S_{1}\) ∪ \(S_{2}\) = S. In the in-order mode \(S_{1}\) contains the smaller digits.
\end{itemize}

Nothing else is allowed: no roots, factorials, logarithms, repeating decimals or constants, and no concatenation of results.

\subsection*{2.2 Values and the problem}\addcontentsline{toc}{subsection}{2.2 Values and the problem}

Every expression has an exact real value; there is no rounding. \({x}^{y}\) is defined only if x \textgreater{} 0 (y real), if x = 0 and y \textgreater{} 0, or if x \textless{} 0 and y is an integer. x/y is defined only if y ≠ 0. In particular \({0}^{0}\) is undefined and complex numbers never occur. An expression with an undefined subexpression is invalid.

Given a constant c ∈ \{π, e, γ, δ\} and n, we seek the minimum of \textbar val(E) − c\textbar{} over the valid expressions E that use each digit 1, \ldots, n exactly once. A formula is \emph{optimal} if it attains that minimum. Two different values with the same exact error are a \emph{tie}, and both are given.

For γ we use the Euler--Mascheroni constant and for δ Feigenbaum's first constant, 4.669 201 609 102 990 671 853 203 820 466\ldots. The certifier takes γ from Arb (\texttt{arb.const\_\allowbreak{}\allowbreak{}euler}, a rigorous ball at any precision {[}7{]}); its first 60 decimals agree with OEIS A001620 {[}13{]}. Feigenbaum's δ is not available in Arb, so the certifier uses the ball 4.669201609102990671853203820466 ± \(10^{-30}\) (30 decimals). The 105 digits of OEIS A006890 {[}12{]}, computed by Briggs {[}10{]}, begin with exactly these 31 digits and continue \ldots4662016\ldots, so the true value is the center plus 2.02·\(10^{-31}\) and lies inside the ball. These digits are also proved: Burbanks, Osbaldestin and Thurlby give a computer-assisted rigorous enclosure of δ that confirms 403 significant digits {[}19{]}. The radius \(10^{-30}\) is at least \(10^{15}\) times smaller than every optimal error for δ (Tables 2--3), so it changes no decision. The independent re-evaluation of the optimal formulas with mpmath also uses δ with these 30 decimals.

\subsection*{2.3 Correspondence with Friedman's statement}\addcontentsline{toc}{subsection}{2.3 Correspondence with Friedman's statement}

{\small
\begin{longtable}[]{@{}>{\raggedright\arraybackslash}p{(\linewidth - 2\tabcolsep) * \real{0.5000}}>{\raggedright\arraybackslash}p{(\linewidth - 2\tabcolsep) * \real{0.5000}}@{}}
\toprule\noalign{}
Rule & Origin \\
\midrule\noalign{}
\endhead
\bottomrule\noalign{}
\endlastfoot
Each digit exactly once, all used & the statement \\
Only + − × / \^{} and parentheses & its list of symbols \\
Digits written together (85) and one decimal point per number (.49, 1.2) & writing digits side by side and the ``.'' of the list; used in his examples \\
Unary minus (−x, \({9}^{-4}\)) & the ``−'' of the list; used in the tables of his page \\
No roots, factorial, logarithm, constants or repeating decimals & not in the list \\
``In order'': digits written in increasing order & his variant ``in order'' \\
Absolute error \textbar x − c\textbar{} & his measure \\
a\^{}\(b^{c}\) = a\^{}(\(b^{c}\)) (only for the notation) & his notation \(({3}^{{2}^{85}})\) \\
\end{longtable}
}

\subsection*{2.4 Conventions that do not affect the results}\addcontentsline{toc}{subsection}{2.4 Conventions that do not affect the results}

\begin{itemize}
\tightlist
\item
  \textbf{Negative base with a non-integer exponent.} We allow only integer exponents. If odd roots were accepted, (−a)\^{}(p/q) = ±\(a^{p/q}\) with q odd, and that value is already reachable with the same digits by unary minus.
\item
  \textbf{\({0}^{0}\).} We leave it undefined. If it were 1, nothing new would appear with 1--n: if some set of digits \(S_{2}\) gives 0, the literal of \(S_{1}\) (never 0) raised to that expression already gives 1 with the same digits, also in order.
\item
  \textbf{Implicit multiplication}, as in 2(3+4), is the same as × and adds no values.
\item
  \textbf{Zeros (variant 0--9).} We use the most permissive rule: any order of the digits is a literal, leading zeros included (05, 0.5, .05). Whatever is optimal under it is optimal under any stricter convention in which the formula is still valid, and our formula for e (with 70 and .6) is valid under any of them.
\end{itemize}

\section*{3. Results}\addcontentsline{toc}{section}{3. Results}

\subsection*{3.1 New records}\addcontentsline{toc}{subsection}{3.1 New records}

Table 1 lists the 21 entries that improve on the records published before this work. The optimal formula for π with 1--9 in any order was also found independently by M. Verheyen (MathWorld {[}2{]}, communicated on 8 November 2026); this work proves that it is optimal. Every error was recomputed from the formula with mpmath at 60--120 digits and with Arb.

\textbf{Table 1.} Improved records (absolute error). The column ``Previous'' gives the best record published before this work.

{\footnotesize
\begin{longtable}[]{@{}>{\raggedright\arraybackslash}p{(\linewidth - 10\tabcolsep) * \real{0.0961}}>{\raggedright\arraybackslash}p{(\linewidth - 10\tabcolsep) * \real{0.0961}}>{\raggedright\arraybackslash}p{(\linewidth - 10\tabcolsep) * \real{0.0961}}>{\raggedright\arraybackslash}p{(\linewidth - 10\tabcolsep) * \real{0.2849}}>{\raggedright\arraybackslash}p{(\linewidth - 10\tabcolsep) * \real{0.1748}}>{\raggedright\arraybackslash}p{(\linewidth - 10\tabcolsep) * \real{0.2519}}@{}}
\toprule\noalign{}
Constant & Variant & n & Formula & Error & Previous \\
\midrule\noalign{}
\endhead
\bottomrule\noalign{}
\endlastfoot
π & any order & 9 & \(5-{\left({.6}^{{(2+9)}^{-.1-.3}}/{4}^{-.7}\right)}^{.8}\) & 2.81·\(10^{-14}\) & 4.11·\(10^{-14}\) (Vlasii) \\
π & any order & 8 & \(2+.6\cdot {.7}^{-.8-{.1}^{-({.3}^{5.4})}}\) & 1.60·\(10^{-12}\) & 5.15·\(10^{-10}\) \\
π & any order & 7 & \({\left(.4\cdot ({(.5-.2/3)}^{-.1}-.6)\right)}^{-.7}\) & 3.02·\(10^{-10}\) & 4.33·\(10^{-8}\) \\
π & in order & 9 & \({(-.1+2)}^{3/\left(-.4-(5-{.6}^{-7+.8})/9\right)}\) & 2.37·\(10^{-10}\) & 2.43·\(10^{-8}\) \\
π & in order & 4 & \(-({.1}^{.2/3})+4\) & 7.11·\(10^{-4}\) & 1.8·\(10^{-3}\) \\
γ & any order & 9 & \({4}^{{.7}^{{6}^{.2}}/.8-1/(.5+{3}^{-.9})}\) & 1.14·\(10^{-15}\) & 4.25·\(10^{-12}\) \\
γ & any order & 8 & \({.4}^{({3}^{.5}-1)/7}/{2}^{-.6}-.8\) & 6.08·\(10^{-13}\) & 8.10·\(10^{-12}\) \\
γ & any order & 7 & \({\left(.4+{\left(-.3/(.1-{.5}^{-.2})\right)}^{7}\right)}^{.6}\) & 2.58·\(10^{-11}\) & 1.21·\(10^{-9}\) \\
γ & in order & 9 & \({\left((.1\cdot {.2}^{-.3}-{4}^{-5})\cdot 6+7/8\right)}^{-.9}\) & 5.82·\(10^{-11}\) & 5.01·\(10^{-10}\) \\
γ & in order & 8 & \({({(-.1)}^{-2}+3)}^{-4-.5+{6}^{.7}/.8}\) & 1.48·\(10^{-9}\) & 1.34·\(10^{-8}\) \\
δ & any order & 9 & \({5}^{{.7}^{-.2}-{({6}^{-3})}^{.4+{.1}^{9.8}}}\) & 1.38·\(10^{-14}\) & 3.04·\(10^{-10}\) \\
δ & any order & 8 & \(5+.67-1-{82}^{-({.3}^{-.4})}\) & 2.23·\(10^{-12}\) & 1.01·\(10^{-8}\) \\
δ & any order & 7 & \((.1+{7}^{-.2}/(.3+.5/6))/.4\) & 3.06·\(10^{-11}\) & 8.59·\(10^{-9}\) \\
δ & in order & 9 & \({\left({.1}^{-({(-2)}^{-3\cdot 4})}+5\right)}^{-.6\cdot .7}\cdot 8+.9\) & 1.16·\(10^{-9}\) & 4.02·\(10^{-9}\) \\
δ & in order & 8 & \((1-2)/.3+{.4}^{5-.6}/7+8\) & 1.27·\(10^{-8}\) & 2.04·\(10^{-8}\) \\
δ & in order & 7 & \(-\left({.1}^{{(-.2/(.3-4))}^{.5}-.6}\right)+7\) & 6.93·\(10^{-8}\) & 3.07·\(10^{-7}\) \\
e & any order & 8 & \({(1+{2}^{-{8}^{7}})}^{{4}^{{.5}^{-6/.3}}}\) & ≈ \(10^{-631\,306}\) & 3.96·\(10^{-47}\) \\
e & any order & 7 & \({(1+{4}^{-36})}^{{.5}^{-72}}\) & 2.88·\(10^{-22}\) & 1.57·\(10^{-8}\) \\
e & in order & 9 & \({(1+{.2}^{{3}^{4}})}^{{5}^{6\cdot (7+8)\cdot .9}}\) & 3.3·\(10^{-57}\) & 1.62·\(10^{-9}\) \\
e & in order & 8 & \({(1+{.2}^{34})}^{{5}^{6\cdot 7-8}}\) & 2.33·\(10^{-24}\) & 2.77·\(10^{-9}\) \\
e & 0--9 & 10 & \({(1+{.2}^{{9}^{70}})}^{{5}^{{3}^{84/.6}}}\) & ≈ \(10^{-4.4\cdot 10^{66}}\) & ≈ \(10^{-5.4\cdot 10^{45}}\) (Fastfaxr) \\
\end{longtable}
}

The ``Previous'' column is the entry of Friedman's page {[}1{]} for that cell, checked cell by cell on 2 November 2026, except for π with 1--9 in any order (Vlasii {[}3, 14{]}; Friedman's page has 6.60·\(10^{-13}\)) and e with 0--9 (Fastfaxr {[}2, 4{]}). On 8 November 2026, after these checks, M. Verheyen communicated to MathWorld {[}2{]} the same formula for π with 1--9 in any order; it is the certified optimum, found independently.

\subsection*{3.2 Complete tables}\addcontentsline{toc}{subsection}{3.2 Complete tables}

Tables 2 and 3 give, for every cell, the minimum error and an optimal formula. The certification (Section 5) also proves that each optimal value is unique: there are no ties in any cell.

\textbf{Table 2.} Digits 1--n in increasing order (certified optima; error to 15 significant digits).

{\footnotesize
\begin{longtable}[]{@{}>{\raggedright\arraybackslash}p{(\linewidth - 8\tabcolsep) * \real{0.0705}}>{\raggedright\arraybackslash}p{(\linewidth - 8\tabcolsep) * \real{0.2324}}>{\raggedright\arraybackslash}p{(\linewidth - 8\tabcolsep) * \real{0.2324}}>{\raggedright\arraybackslash}p{(\linewidth - 8\tabcolsep) * \real{0.2324}}>{\raggedright\arraybackslash}p{(\linewidth - 8\tabcolsep) * \real{0.2324}}@{}}
\toprule\noalign{}
n & π & γ & δ & e \\
\midrule\noalign{}
\endhead
\bottomrule\noalign{}
\endlastfoot
2 & 0.141592653589793 \texttt{(\allowbreak{}1+\allowbreak{}2)} & 0.0537416795786604 \texttt{(\allowbreak{}.1\^{}.2)} & 0.330798390897009 \texttt{(\allowbreak{}1/\allowbreak{}.2)} & 0.281718171540955 \texttt{(\allowbreak{}1+\allowbreak{}2)} \\
3 & 0.0206850065785861 \texttt{(\allowbreak{}.1\^{}-(\allowbreak{}.2+\allowbreak{}.3))} & 0.00148803880217084 \texttt{(\allowbreak{}1.2\^{}-3)} & 0.0276127754902118 \texttt{(\allowbreak{}.1\^{}-(\allowbreak{}.2/\allowbreak{}.3))} & 0.0182818284590452 \texttt{-(\allowbreak{}.1-\/-\allowbreak{}(\allowbreak{}.2-3))} \\
4 & 7.11447819312645·\(10^{-4}\) \texttt{-(\allowbreak{}(\allowbreak{}.1\^{}(\allowbreak{}.2/\allowbreak{}3))-4)} & 6.45368460757851·\(10^{-4}\) \texttt{(\allowbreak{}-(\allowbreak{}1-(\allowbreak{}2\^{}.3))/\allowbreak{}.4)} & 0.00253494243632401 \texttt{(\allowbreak{}1*\allowbreak{}(\allowbreak{}(\allowbreak{}.2/\allowbreak{}.3)+\allowbreak{}4))} & 4.80977891633317·\(10^{-4}\) \texttt{(\allowbreak{}.1+\allowbreak{}(\allowbreak{}2+\allowbreak{}(\allowbreak{}.3\^{}.4)))} \\
5 & 1.39606375849592·\(10^{-4}\) \texttt{(\allowbreak{}(\allowbreak{}.1*\allowbreak{}(\allowbreak{}.2\^{}-(\allowbreak{}.3\^{}4)))\^{}-.5)} & 8.6432061531656·\(10^{-5}\) \texttt{(\allowbreak{}(\allowbreak{}1+\allowbreak{}(\allowbreak{}(\allowbreak{}.2\^{}.3)/\allowbreak{}4))*\allowbreak{}.5)} & 1.17259949560626·\(10^{-4}\) \texttt{(\allowbreak{}(\allowbreak{}.1\^{}.2)*\allowbreak{}(\allowbreak{}-(\allowbreak{}.3-4)/\allowbreak{}.5))} & 1.01102805408584·\(10^{-4}\) \texttt{-(\allowbreak{}.1-(\allowbreak{}(\allowbreak{}-(\allowbreak{}.2-.3))\^{}-.45))} \\
6 & 3.81997701434605·\(10^{-6}\) \texttt{(\allowbreak{}(\allowbreak{}.1\^{}-(\allowbreak{}.2/\allowbreak{}.3))-(\allowbreak{}.4+\allowbreak{}(\allowbreak{}.5+\allowbreak{}.6)))} & 2.72350515411717·\(10^{-7}\) \texttt{(\allowbreak{}1/\allowbreak{}(\allowbreak{}.2+\allowbreak{}(\allowbreak{}.3+\allowbreak{}(\allowbreak{}(\allowbreak{}.4\^{}.5)+\allowbreak{}.6))))} & 8.19606838516996·\(10^{-6}\) \texttt{(\allowbreak{}(\allowbreak{}1/\allowbreak{}(\allowbreak{}.2-(\allowbreak{}3\^{}-4)))-(\allowbreak{}.5\^{}.6))} & 2.30955450623124·\(10^{-6}\) \texttt{(\allowbreak{}(\allowbreak{}1+\allowbreak{}2)\^{}(\allowbreak{}.3+\allowbreak{}(\allowbreak{}(\allowbreak{}.4\^{}5)+\allowbreak{}.6)))} \\
7 & 1.88921690576907·\(10^{-7}\) \texttt{(\allowbreak{}(\allowbreak{}1+\allowbreak{}(\allowbreak{}2\^{}-.3))\^{}(\allowbreak{}.4\^{}-(\allowbreak{}.5/\allowbreak{}(\allowbreak{}.6\^{}.7))))} & 7.13032466917039·\(10^{-8}\) \texttt{(\allowbreak{}(\allowbreak{}-(\allowbreak{}-(\allowbreak{}.1-(\allowbreak{}23\^{}-.45))-.6))\^{}.7)} & 6.93189698729377·\(10^{-8}\) \texttt{-(\allowbreak{}(\allowbreak{}.1\^{}-\/-\allowbreak{}(\allowbreak{}(\allowbreak{}(\allowbreak{}.2/\allowbreak{}-(\allowbreak{}.3-4))\^{}.5)-.6))-7)} & 2.16560552103072·\(10^{-8}\) \texttt{(\allowbreak{}(\allowbreak{}(\allowbreak{}-(\allowbreak{}.1-2)*\allowbreak{}(\allowbreak{}.3\^{}-.4))+\allowbreak{}.5)-(\allowbreak{}.6/\allowbreak{}.7))} \\
8 & 5.27974872105257·\(10^{-9}\) \texttt{(\allowbreak{}(\allowbreak{}.1\^{}-(\allowbreak{}.2/\allowbreak{}.3))-\/-\allowbreak{}(\allowbreak{}(\allowbreak{}(\allowbreak{}4/\allowbreak{}.5)\^{}-6)-(\allowbreak{}.7+\allowbreak{}.8)))} & 1.47786850637223·\(10^{-9}\) \texttt{(\allowbreak{}(\allowbreak{}(\allowbreak{}.1\^{}-2)+\allowbreak{}3)\^{}-(\allowbreak{}4-\/-\allowbreak{}(\allowbreak{}.5-(\allowbreak{}(\allowbreak{}6\^{}.7)/\allowbreak{}.8))))} & 1.27242137555211·\(10^{-8}\) \texttt{-(\allowbreak{}(\allowbreak{}(\allowbreak{}1\^{}.2)/\allowbreak{}.3)-(\allowbreak{}(\allowbreak{}(\allowbreak{}.4\^{}(\allowbreak{}5-.6))/\allowbreak{}7)+\allowbreak{}8))} & 2.33·\(10^{-24}\) (Section 7) \\
9 & 2.3672128172767·\(10^{-10}\) \texttt{(\allowbreak{}(\allowbreak{}-(\allowbreak{}.1-2))\^{}(\allowbreak{}3/\allowbreak{}-(\allowbreak{}.4-(\allowbreak{}-(\allowbreak{}5-(\allowbreak{}.6\^{}-(\allowbreak{}7-.8)))/\allowbreak{}9))))} & 5.82083094786006·\(10^{-11}\) \texttt{(\allowbreak{}(\allowbreak{}(\allowbreak{}(\allowbreak{}(\allowbreak{}.1/\allowbreak{}(\allowbreak{}.2\^{}.3))-(\allowbreak{}4\^{}-5))*\allowbreak{}6)+\allowbreak{}(\allowbreak{}.7/\allowbreak{}.8))\^{}-.9)} & 1.15911801428249·\(10^{-9}\) \texttt{(\allowbreak{}(\allowbreak{}(\allowbreak{}(\allowbreak{}(\allowbreak{}(\allowbreak{}.1\^{}-(\allowbreak{}2\^{}-(\allowbreak{}3*\allowbreak{}4)))+\allowbreak{}5)\^{}-.6)\^{}.7)*\allowbreak{}8)+\allowbreak{}.9)} & 3.3·\(10^{-57}\) (Section 7) \\
\end{longtable}
}

\textbf{Table 3.} Digits 1--n in any order.

{\footnotesize
\begin{longtable}[]{@{}>{\raggedright\arraybackslash}p{(\linewidth - 8\tabcolsep) * \real{0.0705}}>{\raggedright\arraybackslash}p{(\linewidth - 8\tabcolsep) * \real{0.2324}}>{\raggedright\arraybackslash}p{(\linewidth - 8\tabcolsep) * \real{0.2324}}>{\raggedright\arraybackslash}p{(\linewidth - 8\tabcolsep) * \real{0.2324}}>{\raggedright\arraybackslash}p{(\linewidth - 8\tabcolsep) * \real{0.2324}}@{}}
\toprule\noalign{}
n & π & γ & δ & e \\
\midrule\noalign{}
\endhead
\bottomrule\noalign{}
\endlastfoot
2 & 0.141592653589793 \texttt{(\allowbreak{}1+\allowbreak{}2)} & 0.0537416795786604 \texttt{(\allowbreak{}.1\^{}.2)} & 0.330798390897009 \texttt{(\allowbreak{}1/\allowbreak{}.2)} & 0.281718171540955 \texttt{(\allowbreak{}1+\allowbreak{}2)} \\
3 & 0.00825932025645991 \texttt{(\allowbreak{}(\allowbreak{}1/\allowbreak{}.3)-.2)} & 1.34604288092904·\(10^{-4}\) \texttt{(\allowbreak{}3\^{}-(\allowbreak{}.1/\allowbreak{}.2))} & 0.0276127754902118 \texttt{(\allowbreak{}.1\^{}-(\allowbreak{}.2/\allowbreak{}.3))} & 0.0182818284590452 \texttt{(\allowbreak{}(\allowbreak{}3-.2)-.1)} \\
4 & 3.26178699200199·\(10^{-5}\) \texttt{(\allowbreak{}.4\^{}-(\allowbreak{}2.1\^{}.3))} & 1.04078238186349·\(10^{-5}\) \texttt{(\allowbreak{}(\allowbreak{}.1+\allowbreak{}(\allowbreak{}4\^{}-.3))\^{}2)} & 0.00253494243632401 \texttt{(\allowbreak{}(\allowbreak{}.4+\allowbreak{}(\allowbreak{}.2/\allowbreak{}3))/\allowbreak{}.1)} & 9.16653390173138·\(10^{-6}\) \texttt{(\allowbreak{}2\^{}(\allowbreak{}.4\^{}-(\allowbreak{}.1+\allowbreak{}.3)))} \\
5 & 5.92704284501969·\(10^{-6}\) \texttt{(\allowbreak{}3+\allowbreak{}(\allowbreak{}.2/\allowbreak{}(\allowbreak{}.5+\allowbreak{}(\allowbreak{}.4\^{}.1))))} & 2.72350515411717·\(10^{-7}\) \texttt{(\allowbreak{}(\allowbreak{}1+\allowbreak{}(\allowbreak{}(\allowbreak{}.4-.3)\^{}.5))\^{}-2)} & 2.94607845925291·\(10^{-5}\) \texttt{(\allowbreak{}(\allowbreak{}(\allowbreak{}.4+\allowbreak{}(\allowbreak{}.3\^{}-.1))\^{}2)/\allowbreak{}.5)} & 1.68360319123107·\(10^{-6}\) \texttt{(\allowbreak{}(\allowbreak{}3\^{}.5)+\allowbreak{}(\allowbreak{}4\^{}-(\allowbreak{}.1\^{}2)))} \\
6 & 1.38010145419644·\(10^{-8}\) \texttt{(\allowbreak{}(\allowbreak{}4+\allowbreak{}(\allowbreak{}6*\allowbreak{}(\allowbreak{}(\allowbreak{}3\^{}.2)\^{}-.1)))\^{}.5)} & 8.38119699777165·\(10^{-9}\) \texttt{(\allowbreak{}(\allowbreak{}(\allowbreak{}.6\^{}(\allowbreak{}.5\^{}-.3))-.2)\^{}(\allowbreak{}.1+\allowbreak{}.4))} & 1.27021747009822·\(10^{-7}\) \texttt{(\allowbreak{}(\allowbreak{}(\allowbreak{}.3\^{}-(\allowbreak{}.61\^{}.5))-.4)\^{}2)} & 2.66169527336182·\(10^{-8}\) \texttt{(\allowbreak{}1+\allowbreak{}(\allowbreak{}5/\allowbreak{}(\allowbreak{}3/\allowbreak{}(\allowbreak{}(\allowbreak{}6\^{}.2)-.4))))} \\
7 & 3.01840488663709·\(10^{-10}\) \texttt{(\allowbreak{}(\allowbreak{}.4*\allowbreak{}(\allowbreak{}(\allowbreak{}(\allowbreak{}.5-(\allowbreak{}.2/\allowbreak{}3))\^{}-.1)-.6))\^{}-.7)} & 2.57795418531953·\(10^{-11}\) \texttt{(\allowbreak{}(\allowbreak{}.4+\allowbreak{}(\allowbreak{}(\allowbreak{}.3/\allowbreak{}(\allowbreak{}(\allowbreak{}.5\^{}-.2)-.1))\^{}7))\^{}.6)} & 3.0580735760576·\(10^{-11}\) \texttt{(\allowbreak{}(\allowbreak{}.1+\allowbreak{}(\allowbreak{}(\allowbreak{}7\^{}-.2)/\allowbreak{}(\allowbreak{}.3+\allowbreak{}(\allowbreak{}.5/\allowbreak{}6))))/\allowbreak{}.4)} & 2.88·\(10^{-22}\) (Section 7) \\
8 & 1.59872892562979·\(10^{-12}\) \texttt{(\allowbreak{}2+\allowbreak{}(\allowbreak{}.6/\allowbreak{}(\allowbreak{}.7\^{}(\allowbreak{}.8+\allowbreak{}(\allowbreak{}.1\^{}-(\allowbreak{}.3\^{}5.4))))))} & 6.08386083711289·\(10^{-13}\) \texttt{(\allowbreak{}(\allowbreak{}(\allowbreak{}2\^{}.6)*\allowbreak{}(\allowbreak{}.4\^{}(\allowbreak{}(\allowbreak{}(\allowbreak{}3\^{}.5)-1)/\allowbreak{}7)))-.8)} & 2.22546658664114·\(10^{-12}\) \texttt{(\allowbreak{}(\allowbreak{}5.67-1)-(\allowbreak{}82\^{}-(\allowbreak{}.3\^{}-.4)))} & ≈ \(10^{-631\,306}\) (Section 7) \\
9 & 2.81216334342069·\(10^{-14}\) \texttt{(\allowbreak{}5-(\allowbreak{}(\allowbreak{}(\allowbreak{}(\allowbreak{}1+\allowbreak{}3)\^{}.7)*\allowbreak{}(\allowbreak{}.6\^{}(\allowbreak{}(\allowbreak{}2+\allowbreak{}9)\^{}-.4)))\^{}.8))} & 1.13861460673895·\(10^{-15}\) \texttt{(\allowbreak{}4\^{}-(\allowbreak{}(\allowbreak{}1/\allowbreak{}(\allowbreak{}.5+\allowbreak{}(\allowbreak{}3\^{}-.9)))-(\allowbreak{}(\allowbreak{}.7\^{}(\allowbreak{}6\^{}.2))/\allowbreak{}.8)))} & 1.37757967725404·\(10^{-14}\) \texttt{(\allowbreak{}5\^{}(\allowbreak{}(\allowbreak{}.7\^{}-.2)-(\allowbreak{}(\allowbreak{}6\^{}-3)\^{}(\allowbreak{}.4+\allowbreak{}(\allowbreak{}.1\^{}9.8)))))} & ≈ 10\^{}(−8.4·\(10^{39}\)) (Bamberger, Section 7) \\
\end{longtable}
}

Notes. (i) The formulas are written as the certifier prints them, with every parenthesis and with ``−−'' for minus a negative. Table 1 gives shorter forms of the same values. (ii) For e with n ≥ 7 in any order and n ≥ 8 in order, the entries are optimal inside the family F of Section 7, not proved optimal over all expressions. (iii) With n = 1 the only values are ±1 and ±.1, so the optimum is the literal 1 for every constant (errors 2.14, 0.423, 3.67 and 1.72 for π, γ, δ and e); the tables start at n = 2.

\subsection*{3.3 Agreement with Friedman's tables}\addcontentsline{toc}{subsection}{3.3 Agreement with Friedman's tables}

We compared every cell of Tables 2--3 with Friedman's page {[}1{]}. No published entry beats a certified optimum, and every published error matches the error we recompute from its formula. The published entries are optimal in all cells with n ≤ 6, in both variants and for all four constants, except π in order with n = 4; they are also optimal in the following cells with n ≥ 7: π in order n = 7 and 8 (Sabey; de Blaauw), γ in order n = 7 (de Blaauw), e in order n = 7 (de Blaauw) and e in any order n = 9 within the family F (Bamberger). The other 20 cells with n ≤ 9 are improved (Table 1; the formula for π with 1--9 in any order was also found independently by M. Verheyen {[}2{]}); the 21st entry of Table 1 is the 0--9 variant. Table 4 lists the authors of the optimal published entries.

\textbf{Table 4.} Published entries of Friedman's page that are optimal, with n ≤ 6. Initials as on the page: BC Bill Clagett, EF Erich Friedman, GB Gerrit de Blaauw, RS Richard Sabey.

{\small
\begin{longtable}[]{@{}>{\raggedright\arraybackslash}p{(\linewidth - 4\tabcolsep) * \real{0.2123}}>{\raggedright\arraybackslash}p{(\linewidth - 4\tabcolsep) * \real{0.3696}}>{\raggedright\arraybackslash}p{(\linewidth - 4\tabcolsep) * \real{0.4181}}@{}}
\toprule\noalign{}
Constant & Any order, n = 2--6 & In order, n = 2--6 \\
\midrule\noalign{}
\endhead
\bottomrule\noalign{}
\endlastfoot
π & EF, BC, GB, BC, BC & RS (all except n = 4) \\
γ & BC, EF, RS, RS, GB & RS (all) \\
δ & BC, RS, RS, RS, GB & RS (all) \\
e & EF, EF, RS, GB, BC & RS (all) \\
\end{longtable}
}

\section*{4. The search}\addcontentsline{toc}{section}{4. The search}

The search only produces candidates and the window radius E of Section 5. It is not part of any proof, so we describe it briefly.

The search is a dynamic program over subsets of digits combined with a meet-in-the-middle step.

\begin{enumerate}
\def\labelenumi{\arabic{enumi}.}
\tightlist
\item
  \textbf{Values per subset.} For every subset S of the digits, V{[}S{]} is the set of values obtained using each digit of S once. It is built from V{[}A{]} and V{[}S ∖ A{]} with + − × / \^{}, plus the literals and sign.
\item
  \textbf{Pruning.} Duplicates are merged after rounding to about 12 significant digits, and values outside {[}\(10^{-6}\), \(10^{6}\){]} are dropped, as are powers that overflow. Without pruning, V grows too fast beyond 5 digits. The pruning is the main reason why the search proves nothing.
\item
  \textbf{Final crossing.} For every split (A, B) of the n digits and every operation, the value needed in B is solved for: π − x, π/x, x − π, x/π, log\_\allowbreak{}\allowbreak{}x π or \(π^{1/y}\). The nearest value in the sorted table V{[}B{]} is then found. This costs about \textbar V{[}A{]}\textbar·log\textbar V{[}B{]}\textbar{} instead of \textbar V{[}A{]}\textbar·\textbar V{[}B{]}\textbar.
\item
  \textbf{Verification.} Candidates with error below \(10^{-13}\) are rebuilt as expressions and re-evaluated with mpmath at 50 digits.
\end{enumerate}

The version used for the tables (v4) runs in double precision. The corrected version (v5) uses the correctly rounded pow of CORE-MATH {[}8{]}, no fused multiply--add, and deterministic tie-breaking, so its output is the same on Linux and Windows and for any number of threads. With n = 8 in any order the search takes about 11 minutes with 2 threads; with n = 9 in any order it took 27.4 hours for the three constants on an ordinary desktop machine.

The search alone proves nothing: double precision can create values that do not exist (``ghosts'', such as \textbar.3 − (.1 + .2)\textbar{} = 5.6·\(10^{-17}\) instead of 0), lose values that do exist, and amplify rounding errors (\({1.0001}^{9876}\) has a relative error near \(10^{-13}\)). The certifier of Section 5 removes all these gaps.

\section*{5. The certifier}\addcontentsline{toc}{section}{5. The certifier}

\subsection*{5.1 Statement}\addcontentsline{toc}{subsection}{5.1 Statement}

\textbf{Theorem 1.} Let c ∈ \{π, γ, δ\} and n ∈ \{2, \ldots, 9\}, in either mode, or c = e and n ≤ 6 (and n = 7 in order). Let \(R_{n}\) be the set of valid expressions over 1, \ldots, n all of whose intermediate values are 0 or have absolute value in {[}\(10^{-300}\), \(10^{300}\){]}. Then the minimum of \textbar val(E) − c\textbar{} over \(R_{n}\) is the value given in Tables 2--3, and it is attained exactly by the listed value; when several writings are listed, their values are proved equal.

Theorem 1 is the in-range statement; Theorem 2 (Section 6) removes the range restriction. Together they give the main result: the listed value is optimal among all valid expressions.

The cut for e is not an omission. From n = 7 in any order and n = 8 in order, the best approximations of e use the limit \({(1+1/N)}^{N}\) with N a tower of powers, so their intermediate values are far out of range, and the in-range optimum is no longer the true one. For these cells we prove optimality within that family of towers (Section 7 and Appendix B); optimality over all expressions remains open.

\subsection*{5.2 Idea}\addcontentsline{toc}{subsection}{5.2 Idea}

We do not look for the best value; we prove that no expression comes closer to c than the record. It is enough to find every expression whose exact value \emph{could} lie in the window W = {[}c − E\emph{, c + E}{]} and to decide each of them exactly. Two simplifications, both proved, shrink the problem without changing the set of values: the sign is handled apart (the value set is symmetric thanks to unary minus), and powers with a negative base can be reduced to positive bases.

The same \emph{window query} --- which expressions over the digits S can have a value in W? --- is answered by two engines with different arithmetic:

\begin{itemize}
\tightlist
\item
  \textbf{Motor A} (C, fast). For sets of up to 5 digits it keeps tables. Each entry is a \emph{class}: a double that identifies it and an interval that contains the exact value of all its expressions. Values near 1 are stored by an enclosure of their logarithm. For larger sets it splits S = \(S_{1}\) ∪ \(S_{2}\) with \(S_{1}\) in a table. For each value s of \(S_{1}\) and each operation it computes the window W′ where the value of \(S_{2}\) would have to lie, and queries \(S_{2}\) recursively. It uses only correctly rounded double arithmetic (IEEE 754 and CORE-MATH's pow), so every bound is the predecessor or successor of a correctly rounded result.
\item
  \textbf{Motor B} (Python with Arb, exact). It answers the same query with Arb balls of arbitrary precision. It decides whether a value is exactly 0 or 1 with exact arithmetic on rationals and normalized radicals; distinct radicals are linearly independent over ℚ (Besicovitch), which also decides sums. Motor B resolves the \emph{singular} classes of motor A's tables (classes whose interval contains 0, or 1 when they act as the base of a power), verifies one by one the candidate paths that motor A outputs, and gives the concrete formulas.
\item
  \textbf{Motor C} (C with Arb) is a pre-filter for the large cells. It drops a candidate path only when Arb proves that every value the path can produce is outside W, out of range or undefined.
\end{itemize}

A cell is certified when motor A's candidate list, checked exactly by motor B, contains no value with error below the record except the record itself, and no case is left unresolved. The optimal formulas are finally re-evaluated by an independent evaluator (mpmath, parsing the string), which also checks the digits, the validity and the range of every subexpression.

\subsection*{5.3 Correctness}\addcontentsline{toc}{subsection}{5.3 Correctness}

The proof is in Appendix A. Its steps are:

\begin{itemize}
\tightlist
\item
  \textbf{Reductions:} the sign is symmetric (Lemma 1), bases can be taken positive (Lemma 2), and nine forms suffice (Lemma 3).
\item
  \textbf{Arithmetic:} correct rounding gives outward enclosures (Lemma 4); extrema over a box are at the corners (Lemma 5); exact rationals are handled with 128-bit integers (Lemma 6); exp and ln are bounded through CORE-MATH's pow (Lemma 7), with sharper bounds for small arguments (Lemma 8).
\item
  \textbf{Tables:} every derivation rule is sound (Lemma 9), and singular derivations are resolved exactly by motor B (Lemma 10). By induction on \textbar S\textbar, every in-range value over S lies in a class of the table (Theorem A).
\item
  \textbf{Queries:} the preimage windows contain every admissible partner (Lemma 11), the leaf search misses no class (Lemma 12), and by induction the query covers every in-range expression with value in the window (Theorem B).
\item
  \textbf{Verification:} motor B's arithmetic decides signs and equalities (Lemma 13) and its queries are complete (Lemma 14). Each candidate path is then decided exactly (Theorem C), and motor C only drops paths that motor B would find empty (Lemma 15).
\end{itemize}

Appendix A.8 lists the independent checks and two bugs they found, and A.10 maps each lemma to the code.

\subsection*{5.4 What is trusted}\addcontentsline{toc}{subsection}{5.4 What is trusted}

IEEE 754 binary64 hardware (SSE2), compilation without contraction (−ffp-contract=off), the correct rounding of CORE-MATH's pow (proved by its authors {[}8{]}; we also checked 28 737 cases against Arb), Arb (python-flint 0.9.0; FLINT 3.0.1 for motor C), Python's exact integers and fractions, and the linear independence of real radicals (Besicovitch 1940; Mordell 1953; Siegel 1972) {[}9, 17, 18{]}. Appendix A.9 states it precisely.

\subsection*{5.5 Runs}\addcontentsline{toc}{subsection}{5.5 Runs}

Fifty-six cells were certified on an ordinary PC: π, γ and δ for n = 2--9 in order and n = 2--8 in any order, and e for n ≤ 6 in both modes and n = 7 in order. All have 0 unresolved cases and a unique optimal value. The three cells with n = 9 in any order ran on a desktop PC (Intel Core i5-7600, 4 cores, 32 GB of RAM, Windows 10) with the optimized motors A and C described after Lemma 15 (Appendix A.6). Motor A (first version) wrote one candidate file for the three windows in a single pass of 47 hours with 4 threads (168 678 s); the optimized version of Lemma 15 needs about 7 hours for one window. Motor C reduced its 379 019 117 paths to 499 690 at 1024 bits (32.5 hours, 4 processes) and to 137 445 at 2048 bits (0.3 hours). Motor B then decided the remaining 18 371 / 30 904 / 22 911 distinct paths for π / γ / δ in 34 / 60 / 45 minutes with 2 processes. Linux and Windows builds give byte-identical outputs.

\textbf{Table 5.} Certification effort (selected cells).

{\small
\begin{longtable}[]{@{}>{\raggedright\arraybackslash}p{(\linewidth - 6\tabcolsep) * \real{0.2249}}>{\raggedright\arraybackslash}p{(\linewidth - 6\tabcolsep) * \real{0.3385}}>{\raggedright\arraybackslash}p{(\linewidth - 6\tabcolsep) * \real{0.1117}}>{\raggedright\arraybackslash}p{(\linewidth - 6\tabcolsep) * \real{0.3249}}@{}}
\toprule\noalign{}
Cell & Candidates (motor A) & Values in window & Time \\
\midrule\noalign{}
\endhead
\bottomrule\noalign{}
\endlastfoot
π, γ, δ; n = 9 in order & 3797 / 3224 / 4090 & 1 / 1 / 3 & ≈ 2.5 min each \\
π, γ, δ; n = 7 any order & 299 / 309 / 315 & 1 / 2 / 1 & ≈ 2.3 min each \\
π, γ, δ; n = 8 any order & 125 847 / 112 880 / 137 963 & 2 / 3 / 13 & ≈ 45 min each (motor B) \\
π, γ, δ; n = 9 any order & 18 371 / 30 904 / 22 911 (after motor C; 379 019 117 before it, the three together) & 9 / 3 / 5 & motor C 33 h; motor B 34 / 60 / 45 min \\
\end{longtable}
}

\section*{6. Expressions with huge or tiny intermediate values}\addcontentsline{toc}{section}{6. Expressions with huge or tiny intermediate values}

Theorem 1 restricts the intermediate values to 0 or {[}\(10^{-300}\), \(10^{300}\){]}. The restriction is not cosmetic. For example, \({2}^{{3}^{7}}/{9}^{{5}^{4}}\) ≈ 8.9·\(10^{61}\) is a perfectly good value whose two intermediates exceed \(10^{300}\), and the best formulas for e are built from such towers. This section removes the restriction.

\textbf{Theorem 2.} For every cell of Theorem 1, no valid expression, with no restriction on its intermediate values, comes closer to c than the optimum of Theorem 1. So the optima of Tables 2--3 are optimal over all expressions of Section 2. The same computation, with the window of the optimum without unary minus, gives the corresponding statement for the variant of Section 8 (Appendix C.6).

\textbf{How rare such expressions are.} An approximate enumeration in floating point of every expression, with no range limit, shows how little they add:

\begin{itemize}
\tightlist
\item
  \textbf{Up to 4 digits:} no value within the range {[}\(10^{-300}\), \(10^{300}\){]} is reached only through an extreme intermediate (extreme values themselves, such as \({9}^{{8}^{7}}\), do occur).
\item
  \textbf{Five digits:} 20 to 300 new values per set, all with \textbar lo\(g_{10}\) u\textbar{} ≥ 30.
\item
  \textbf{Digits 1--6:} 41 247 new values. The closest are at distance 3.8·\(10^{-3}\) from π, 5.9·\(10^{-3}\) from γ, 3.3·\(10^{-2}\) from δ and 1.1·\(10^{-2}\) from e, against certified optima of order \(10^{-8}\).
\end{itemize}

The mechanisms that bring an extreme value back near c are four:

\begin{enumerate}
\def\labelenumi{\arabic{enumi}.}
\tightlist
\item
  \textbf{Damping:} \({x}^{y}\) with x extreme and y small.
\item
  \textbf{Pairs of extremes:} quotients, or near-cancelling differences.
\item
  \textbf{Thresholds:} values just beyond \(10^{300}\) times moderate factors.
\item
  \textbf{Shadows:} y + ε with ε tiny, and 1 + ε raised to a huge power, which is the mechanism of e.
\end{enumerate}

\textbf{The proof} is a second certified computation, motor X (Appendix C):

\begin{itemize}
\tightlist
\item
  \textbf{Generation lemma.} In an expression with an out-of-range intermediate, that intermediate is in the left operand, in the right operand, or is the result of combining two in-range operands. So only the expressions that contain an extreme value need to be explored. The lemma is about the structure of expressions, not their values, so it also covers expressions whose value happens to coincide with that of an in-range one.
\item
  \textbf{Representation.} Motor X represents every value by an enclosure of its logarithm, rounded outward with CORE-MATH, plus three channels for values a double cannot hold:

  \begin{itemize}
  \tightlist
  \item
    ln\textbar ln u\textbar{} for values within \(10^{-300}\) of 1 or beyond \(e^{10^{308}}\);
  \item
    exact rational anchors for ``rational + ε'';
  \item
    a third level for (1 + t)\^{}s with t ≈ e\^{}(−e\^{}1125).
  \end{itemize}
\item
  \textbf{Exact decisions.} Coincidences between in-range values, such as \({4}^{-.3}\) = \({.5}^{.6}\), and in-range bases that might equal 1 are decided exactly by motor B.
\item
  \textbf{Queries.} The window queries of motor A, with motor A's tables of in-range classes, extend to these values. A query returns every value that needs an out-of-range intermediate and may lie in the window of the record.
\end{itemize}

The query returns a finite list of enclosures (\emph{hits}). Every expression of every hit is then evaluated by \texttt{lvl.py}, an independent interval arithmetic with unbounded exponents. For every cell, every such value is proved to lie outside the window.

\textbf{Table 6.} Motor X on the cells of Theorem 1: number of top-level splits, hits per constant (π / γ / δ / e; values that may lie in the window, before the check with \texttt{lvl}), and CPU time for all constants together. The window has the radius of the larger of the two optima, with and without unary minus, enlarged by \(10^{-9}\) in relative terms, or of an independently checked upper bound of it (see below), and each value is checked against both windows. Every cell listed has 0 values left undecided for either variant, 0 values certainly inside either window, and 0 undecided operations of motor B. n = 2 needs no computation: two digits cannot produce an out-of-range value.

{\footnotesize
\begin{longtable}[]{@{}>{\raggedright\arraybackslash}p{(\linewidth - 8\tabcolsep) * \real{0.1364}}>{\raggedright\arraybackslash}p{(\linewidth - 8\tabcolsep) * \real{0.1364}}>{\raggedright\arraybackslash}p{(\linewidth - 8\tabcolsep) * \real{0.1364}}>{\raggedright\arraybackslash}p{(\linewidth - 8\tabcolsep) * \real{0.3856}}>{\raggedright\arraybackslash}p{(\linewidth - 8\tabcolsep) * \real{0.2052}}@{}}
\toprule\noalign{}
Mode & n & Splits & Hits & CPU (s) \\
\midrule\noalign{}
\endhead
\bottomrule\noalign{}
\endlastfoot
in order & 3--8 & 2--7 & 0 in every cell & ≤ 5 \\
in order & 9 & 8 & 40 / 40 / 40 & 279 \\
any order & 3 & 3 & 0 / 0 / 0 / 0 & 0.0 \\
any order & 4 & 7 & 0 / 0 / 0 / 0 & 0.1 \\
any order & 5 & 15 & 0 / 0 / 0 / 0 & 1.9 \\
any order & 6 & 31 & 0 / 0 / 0 / 0 & 221 \\
any order & 7 & 63 & 0 / 0 / 0 & 674 \\
any order & 8 & 127 & 490 / 482 / 492 & 1030 \\
any order & 9 & 255 & 1 845 905 / 1 757 123 / 1 848 564 & 1.21·\(10^{6}\) ‡ \\
\end{longtable}
}

‡ 336 CPU hours for the three constants together: 89 h for 54 splits on the PC of Section 5.5 and 247 h for 201 splits on Linux machines with 4 cores. The other rows ran on one Linux machine; the cell n = 8 in any order took 1 919 s on the PC of Section 5.5.

For n = 9 in order, the 40 hits of each constant were rebuilt with all their derivations and evaluated with \texttt{lvl}: all are provably outside both windows. The last cases needed exact factorized powers, such as \textbar{}\({.1}^{-{2}^{3\cdot (4+5)}}-{(.6-.7)}^{-{8}^{9}}\)\textbar{} = 0 (both powers equal \({10}^{134217728}\)). The first version of motor X gave 22 hits per constant here; Appendix C.9 explains the difference, which changes no result. For e, Table 6 covers the cells of Theorem 1 (n ≤ 6, and n = 7 in order).

For the variant without unary minus the radius may be an upper bound of its optimum \(E_{11}\) rather than \(E_{11}\) itself: the error of a minus-free formula checked independently (δ with n = 8, and π and δ with n = 9). This is enough, because a window of radius r ≥ \(E_{11}\) contains W(\(E_{11}\)). With n = 9 in any order the radius of the window with unary minus is also an upper bound of E*: 3.23·\(10^{-14}\), 2.12·\(10^{-15}\) and 2.2·\(10^{-14}\) for π, γ and δ (the optima of Table 3 are 2.81·\(10^{-14}\), 1.14·\(10^{-15}\) and 1.38·\(10^{-14}\)). The cell n = 8 in any order was also run on the Windows machine of Section 5.5, and gave the same hits as the Linux run for each of the 381 pairs (constant, split). The independent review of this work repeated δ with n = 8 and the minus-free window on its own machine: 127 splits, 492 hits, 0 undecided values.

\textbf{n = 9 in any order.} The 255 splits were computed on the PC of Section 5.5 (54 splits) and on Linux machines (201 splits), with the same code and radii. Each Linux run first computed a control split with known hits, and every result was checked again on the PC of Section 5.5 before the closure. All 255 splits are closed for π, γ and δ and for both variants: 0 undecided values, 0 values certainly inside either window and 0 undecided operations of motor B. During this run a sign error was found in one version of the evaluator \texttt{lvl} (Appendix C.9, error 8). Every split evaluated with that version was recomputed with the corrected one; the 219 pairs (split, constant) that the faulty run had completed have the same hits in both runs. The records of n = 9 were checked by five versions of \texttt{lvl}, listed in Appendix C.9. A sample of whole splits made with each earlier version (splits 9, 17, 70 and 191) was computed again with the final version: the same hits, and no undecided value.

\section*{7. The constant e}\addcontentsline{toc}{section}{7. The constant e}

\subsection*{7.1 Why e is different}\addcontentsline{toc}{subsection}{7.1 Why e is different}

Unlike π, e has a limit that the allowed operations can build: \({(1+1/N)}^{N}\) → e, with error about e/(2N). The precision therefore grows with the size of N, and power towers make N enormous with a few digits. The art is to write the same N twice with disjoint digits. The known records rest on exact identities between towers:

{\small
\begin{longtable}[]{@{}>{\raggedright\arraybackslash}p{(\linewidth - 4\tabcolsep) * \real{0.2711}}>{\raggedright\arraybackslash}p{(\linewidth - 4\tabcolsep) * \real{0.5375}}>{\raggedright\arraybackslash}p{(\linewidth - 4\tabcolsep) * \real{0.1914}}@{}}
\toprule\noalign{}
Formula & Key identity & N \\
\midrule\noalign{}
\endhead
\bottomrule\noalign{}
\endlastfoot
Sabey 2004 & \({9}^{{4}^{42}}\) = \({3}^{2\cdot {2}^{84}}\) = \({3}^{{2}^{85}}\) & \({3}^{{2}^{85}}\) \\
Bamberger 2024 & \({.2}^{{9}^{42}}\) = \({5}^{-{3}^{84}}\) & \({5}^{{3}^{84}}\) \\
Fastfaxr (Reddit) & \({2}^{0}\) = 1, and \({9}^{{4}^{76}}\) = \({3}^{{2}^{153}}\) = \({3}^{{8}^{51}}\) & \({3}^{{2}^{153}}\) \\
\end{longtable}
}

Three devices appear. Bamberger uses the decimal .2 = 1/5 as an inverse base, which saves the minus sign and gives base 5. Fastfaxr writes 1 as \({2}^{0}\) so that the 0 is used up and the 1 stays free. And base 3 in the exponent pays off: \({3}^{84}\) grows much faster than \({2}^{85}\).

\subsection*{7.2 Error bounds}\addcontentsline{toc}{subsection}{7.2 Error bounds}

For every real N ≥ 1, (e/2N)(1 − 11/(12N)) \textless{} e \(-{(1+1/N)}^{N}\) \textless{} e/(2N) (Lemma B1). A formula is correct to k decimals if its error is below ½·\(10^{-k}\); with u = lo\(g_{10}\)(N/e), k = ⌊u⌋ in every case (Appendix B.1). Table B1 gives the exact count for each record: 21, 631 305 and 8 368 428 989 068 425 943 817 590 916 445 001 887 164 decimals for 1--7, 1--8 and 1--9 (Bamberger), 23 and 56 for 1--8 and 1--9 in order, and a 67-digit count for 0--9.

\subsection*{7.3 The family F and its optima}\addcontentsline{toc}{subsection}{7.3 The family F and its optima}

F consists of the expressions (U + A\^{}(±\(B^{X}\)))\^{}(C\^{}(±\(D^{Y}\))) and (U + \(A^{±X}\))\^{}(\(C^{±Y}\)), where U, A, B, C, D, X and Y are arbitrary expressions over disjoint groups of digits that together use every digit, U = 1, and the two large numbers are exactly reciprocal, so that the expression is \({(1+1/N)}^{N}\). The proof of optimality inside F (Appendix B.3) has no artificial bounds. It uses a size lemma, an exact enumeration of all blocks of up to 3 digits, an Arb sieve and exact verification of every overlap with prime exponents; no case is left unresolved and no transcendence theorem is needed. Its results are:

{\small
\begin{longtable}[]{@{}>{\raggedright\arraybackslash}p{(\linewidth - 4\tabcolsep) * \real{0.2258}}>{\raggedright\arraybackslash}p{(\linewidth - 4\tabcolsep) * \real{0.5230}}>{\raggedright\arraybackslash}p{(\linewidth - 4\tabcolsep) * \real{0.2513}}@{}}
\toprule\noalign{}
Digits & Maximum of F & Correct decimals \\
\midrule\noalign{}
\endhead
\bottomrule\noalign{}
\endlastfoot
1--7 & N = \({2}^{72}\): \({(1+{4}^{-36})}^{{.5}^{-72}}\) & 21 \\
1--8 & N = \({2}^{{2}^{21}}\): \({(1+{2}^{-{8}^{7}})}^{{4}^{{.5}^{-6/.3}}}\) & 631 305 \\
1--9 & N = \({5}^{{3}^{84}}\) (Bamberger) & ≈ 8.37·\(10^{39}\) \\
1--8, in order & N = \({5}^{34}\): \({(1+{.2}^{34})}^{{5}^{6\cdot 7-8}}\) & 23 \\
1--9, in order & N = \({5}^{81}\): \({(1+{.2}^{{3}^{4}})}^{{5}^{6\cdot (7+8)\cdot .9}}\) & 56 \\
0--9 & N = \({5}^{{3}^{140}}\): \({(1+{.2}^{{9}^{70}})}^{{5}^{{3}^{84/.6}}}\) & ≈ 4.38·\(10^{66}\) \\
\end{longtable}
}

Towers with three or more stories exist with 1--9. One example is \({\left(1+{.5}^{{3}^{{7}^{2}}}\right)}^{{8}^{{({9}^{6})}^{4}}}\), with N = \({2}^{{3}^{49}}\) and about 7.2·\(10^{22}\) decimals, but none beats Bamberger's formula. With four or more stories, ln N ≤ 2.1·\(10^{9}\) (Appendix B.4, Table B2).

\subsection*{7.4 The variant 0--9}\addcontentsline{toc}{subsection}{7.4 The variant 0--9}

With the digit 0 and the permissive rule of Section 2.4, the maximum of F is

\({(1+{.2}^{{9}^{70}})}^{{5}^{{3}^{84/.6}}}\),

using \({9}^{70}\) = \({3}^{140}\) and 84/.6 = 140, so that N = \({5}^{{3}^{140}}\). N has 4 379 597 503 586 508 305 506 762 455 287 865 626 188 758 462 213 921 632 491 859 893 687 digits, and the formula is correct to 4 379 597 503 586 508 305 506 762 455 287 865 626 188 758 462 213 921 632 491 859 893 686 decimals (about 4.4·\(10^{66}\)) in the sense of Appendix B.1, which removes the ``one fewer'' ambiguity. The best previously published value we found is Fastfaxr's 5.4·\(10^{45}\).

\subsection*{7.5 What remains open for e}\addcontentsline{toc}{subsection}{7.5 What remains open for e}

Optimality outside F --- for instance with other limit forms, such as \({(1+u/N)}^{N/u}\), or sums and products of towers --- is open. So is a proof that no exact identity with 1--9 can beat N = \({5}^{{3}^{84}}\).

\section*{8. The variant without unary minus}\addcontentsline{toc}{section}{8. The variant without unary minus}

Unary minus is the only operation of Section 2 that is not a binary operator, and some of the published entries look as if they were found without it. We therefore also solve the variant in which the rule −E is removed. Negative values still arise from subtraction; everything else is unchanged, including the range of Theorem 1.

\textbf{Method.} Minus-free expressions are a subset of those of Section 2, so the certified optimum with unary minus is a lower bound. Without unary minus the value set is no longer symmetric. For each magnitude u that motor B finds in the window we therefore track two flags: whether +u is the value of a minus-free expression, and whether −u is. The flags are propagated through exactly the same enumeration as the certificate (tables, window queries, path verification) and are OR-ed over all derivations. The search starts with the window of radius E* and, if it finds no positively reachable value, multiplies the radius by 4 and repeats. For n = 8 and 9 in any order it makes a single pass instead, with a window whose radius is the error of a minus-free formula found by a meet-in-the-middle search and checked independently at 80 digits; that radius is an upper bound of the optimum, so the minimum of the positively reachable values in the window is the optimum. (For δ with n = 8 the formula found, \((.1/(.5\cdot .8))+3/({7}^{.2}\cdot .46)\), has exactly the value of the optimum with n = 7 in any order, so the radius coincides with the window of that cell.) With n = 9 the pipeline is that of Theorem 1: motor A with this window (about 7 hours per constant on the desktop PC of Section 5.5), motor C (about 4.5 hours) and the flagged verification of the 18 451 / 31 016 / 23 230 surviving paths (4--5 hours, 0 unresolved). Every path that survived motor C with the window of E* also survives with the wider window, except two paths of γ, which the exact verification shows to have no value in either window. A power with a negative base is accepted only when its exponent is an exact rational integer; an exponent that is an integer only through an identity between irrational powers, such as (7\^{}√2)\^{}√8 = 2401, would be rejected, and a minus-free writing of that form would be missed. An instrumented rerun of the verification of all 59 cells, which counts every negative base whose exponent is not an exact rational but whose ball may contain an integer, found no such case and reproduced every result, so no writing was missed this way. The flagged engine agrees exactly with a brute-force enumeration of all minus-free expressions for every set of up to 4 digits (any order) and the in-order sets of 3--4 digits.

\textbf{Results.} Table 7 gives the optima. In 22 of the 59 cells the optimum does not change (with n ≥ 7 only in two: δ in order and γ in any order, both with n = 8); in the others the record uses an essential unary minus, such as \(4-{.1}^{.2/3}\) written in order as \(-({.1}^{.2/3})+4\). Three published entries that our search improved coincide exactly with the minus-free optimum:

\begin{itemize}
\tightlist
\item
  π in order with n = 4: 1.80·\(10^{-3}\), \({((.1+{2}^{.3}))}^{4}\);
\item
  γ in order with n = 8: 1.34·\(10^{-8}\);
\item
  δ in order with n = 9: 4.02·\(10^{-9}\), \(1+{2}^{{\left((3-.4\cdot (.5-{.6}^{7}))-.8\right)}^{.9}}\).
\end{itemize}

This is not a systematic convention of the contributors. Their in-order entries use unary minus elsewhere, even in front of the first digit: \(-.1+{(-.2+.3)}^{-.45}\) for e with n = 5 (Sabey) and \(-((.1-2)/{.3}^{.4}-.5+6/7)\) for e with n = 7 (de Blaauw), both optimal. The other improved in-order entries are not minus-free optima either: for π with n = 9 and γ with n = 9 the published error (2.43·\(10^{-8}\) and 5.01·\(10^{-10}\)) is worse than the minus-free optimum (1.05·\(10^{-9}\) and 9.34·\(10^{-11}\)), and for δ with n = 7 and 8 the published formulas use a unary minus. The three coincidences are better read as a sign that hand searches found the best formulas whose sign structure is simple, and missed those that need a minus in front of a whole group, such as \(-({.1}^{.2/3})+4\).

\textbf{Table 7.} Optimal errors without unary minus (optimal over all minus-free expressions: Section 8 and Theorem 2, Appendix C.6). In brackets, the optimum with unary minus when it is different.

{\footnotesize
\begin{longtable}[]{@{}>{\raggedright\arraybackslash}p{(\linewidth - 16\tabcolsep) * \real{0.0432}}>{\raggedright\arraybackslash}p{(\linewidth - 16\tabcolsep) * \real{0.1192}}>{\raggedright\arraybackslash}p{(\linewidth - 16\tabcolsep) * \real{0.1222}}>{\raggedright\arraybackslash}p{(\linewidth - 16\tabcolsep) * \real{0.1162}}>{\raggedright\arraybackslash}p{(\linewidth - 16\tabcolsep) * \real{0.1162}}>{\raggedright\arraybackslash}p{(\linewidth - 16\tabcolsep) * \real{0.1222}}>{\raggedright\arraybackslash}p{(\linewidth - 16\tabcolsep) * \real{0.1222}}>{\raggedright\arraybackslash}p{(\linewidth - 16\tabcolsep) * \real{0.1222}}>{\raggedright\arraybackslash}p{(\linewidth - 16\tabcolsep) * \real{0.1162}}@{}}
\toprule\noalign{}
n & π in order & γ in order & δ in order & e in order & π any order & γ any order & δ any order & e any order \\
\midrule\noalign{}
\endhead
\bottomrule\noalign{}
\endlastfoot
2 & 1.42·\(10^{-1}\) & 5.37·\(10^{-2}\) & 3.31·\(10^{-1}\) & 2.82·\(10^{-1}\) & 1.42·\(10^{-1}\) & 5.37·\(10^{-2}\) & 3.31·\(10^{-1}\) & 2.82·\(10^{-1}\) \\
3 & 2.16·\(10^{-2}\) (2.07·\(10^{-2}\)) & 1.10·\(10^{-2}\) (1.49·\(10^{-3}\)) & 3.08·\(10^{-2}\) (2.76·\(10^{-2}\)) & 1.83·\(10^{-2}\) & 8.26·\(10^{-3}\) & 2.78·\(10^{-3}\) (1.35·\(10^{-4}\)) & 3.08·\(10^{-2}\) (2.76·\(10^{-2}\)) & 1.83·\(10^{-2}\) \\
4 & 1.80·\(10^{-3}\) (7.11·\(10^{-4}\)) & 1.35·\(10^{-3}\) (6.45·\(10^{-4}\)) & 2.53·\(10^{-3}\) & 4.81·\(10^{-4}\) & 4.07·\(10^{-4}\) (3.26·\(10^{-5}\)) & 2.89·\(10^{-5}\) (1.04·\(10^{-5}\)) & 2.53·\(10^{-3}\) & 4.31·\(10^{-4}\) (9.17·\(10^{-6}\)) \\
5 & 3.74·\(10^{-4}\) (1.40·\(10^{-4}\)) & 8.64·\(10^{-5}\) & 7.98·\(10^{-4}\) (1.17·\(10^{-4}\)) & 2.82·\(10^{-4}\) (1.01·\(10^{-4}\)) & 5.93·\(10^{-6}\) & 2.07·\(10^{-6}\) (2.72·\(10^{-7}\)) & 4.56·\(10^{-5}\) (2.95·\(10^{-5}\)) & 2.67·\(10^{-6}\) (1.68·\(10^{-6}\)) \\
6 & 7.35·\(10^{-6}\) (3.82·\(10^{-6}\)) & 2.72·\(10^{-7}\) & 2.62·\(10^{-5}\) (8.20·\(10^{-6}\)) & 2.31·\(10^{-6}\) & 1.38·\(10^{-8}\) & 1.70·\(10^{-8}\) (8.38·\(10^{-9}\)) & 6.17·\(10^{-7}\) (1.27·\(10^{-7}\)) & 2.66·\(10^{-8}\) \\
7 & 6.78·\(10^{-7}\) (1.89·\(10^{-7}\)) & 1.78·\(10^{-7}\) (7.13·\(10^{-8}\)) & 4.51·\(10^{-7}\) (6.93·\(10^{-8}\)) & 7.90·\(10^{-8}\) (2.17·\(10^{-8}\)) & 2.43·\(10^{-9}\) (3.02·\(10^{-10}\)) & 7.30·\(10^{-10}\) (2.58·\(10^{-11}\)) & 3.03·\(10^{-9}\) (3.06·\(10^{-11}\)) & --- \\
8 & 1.64·\(10^{-8}\) (5.28·\(10^{-9}\)) & 1.34·\(10^{-8}\) (1.48·\(10^{-9}\)) & 1.27·\(10^{-8}\) & --- & 1.78·\(10^{-12}\) (1.60·\(10^{-12}\)) & 6.08·\(10^{-13}\) & 1.20·\(10^{-11}\) (2.23·\(10^{-12}\)) & --- \\
9 & 1.05·\(10^{-9}\) (2.37·\(10^{-10}\)) & 9.34·\(10^{-11}\) (5.82·\(10^{-11}\)) & 4.02·\(10^{-9}\) (1.16·\(10^{-9}\)) & --- & 8.12·\(10^{-14}\) (2.81·\(10^{-14}\)) & 1.98·\(10^{-15}\) (1.14·\(10^{-15}\)) & 4.96·\(10^{-14}\) (1.38·\(10^{-14}\)) & --- \\
\end{longtable}
}

A few optimal minus-free formulas:

\begin{itemize}
\tightlist
\item
  π, any order, n = 4: the literal 3.142;
\item
  π, any order, n = 7: \({5}^{{(.7/.6-.3)}^{1/.42}}\);
\item
  γ, in order, n = 9: \({(1-.2)}^{.3/(.4+.5)}\cdot .6\cdot {.7}^{.8-.9}\);
\item
  π, in order, n = 9: \({({(.1\cdot .2-3)}^{4})}^{{({.5}^{6}+.7)}^{.8}}/9\);
\item
  π, any order, n = 9: \(2/\left({\left({.4}^{5/(9-.3)}/6-.7/8\right)}^{.1}\right)\);
\item
  γ, any order, n = 9: \({.4}^{.6-{(.2/{.1}^{.3}-.5/{8}^{7})}^{9}}\);
\item
  δ, any order, n = 9: \(.4+{(({8}^{.1}+{5}^{.23})/.6)}^{.97}\).
\end{itemize}

The complete list, with every formula, is in the results file \texttt{c11.jsonl}.

\section*{9. Why π and e are different problems}\addcontentsline{toc}{section}{9. Why π and e are different problems}

For e the records improve by thousands of orders of magnitude from one n to the next; for π, γ and δ by one to three. This section explains the contrast. It has three parts: a structural one (what the operations can build), a quantitative one (a random model that predicts the optima of Tables 2--3), and a counterpoint (with complex intermediates, π behaves like e).

\subsection*{9.1 Limits that the operations can build}\addcontentsline{toc}{subsection}{9.1 Limits that the operations can build}

With + − × / and \({x}^{y}\) = exp(y ln x) for x \textgreater{} 0, every value is a real exponential--logarithmic number: it lies in the smallest subfield of ℝ that contains ℚ and is closed under exp and under ln of positive elements. Some constants are limits of very short expressions of this kind. The clearest case is e = lim \({(1+1/N)}^{N}\), with error about e/(2N) (Lemma B1). Power towers make N astronomically large with few digits, provided that N can be written twice with disjoint digits, which is what the identities of Section 7 do. Logarithms are limits too: \(N\cdot ({x}^{1/N}-1)\) → ln x.

No such limit is known for π. Chow {[}6{]} works with the complex field of EL numbers, which contains π = −i·log(−1), and asks in his final section what the right real analogue would be. Whether π lies in the real exponential--logarithmic closure of ℚ is open {[}6, 15{]}, although one would expect that it does not; we know of no proof of this, even conditional on Schanuel's conjecture. We use none of this: the point is only that, without a limit form, an approximation of π with these operations is a numerical coincidence, and its quality is governed by how many distinct values the digits produce near π.

\subsection*{9.2 A random model}\addcontentsline{toc}{subsection}{9.2 A random model}

Let \(D_{n}\)(c) be the number of distinct values over the digits 1, \ldots, n that lie within ±5 \% of c (in range, as in Theorem 1). If these values were independent and uniformly spread over that interval of half-width Δ = 0.05·c, the smallest error would be close to an exponential variable with mean Δ/(\(D_{n}\) + 1). Its 5 \% and 95 \% quantiles are 0.051 and 3.0 times the mean.

We counted \(D_{n}\) by an exhaustive enumeration in double precision, merging values closer than \(10^{-12}\) in relative terms (an approximate count, heuristic and not part of any proof), for n ≤ 6 in any order and n ≤ 7 in order. For π, \(D_{n}\) = 21, 900, 46 260 and 4 234 516 for n = 3--6 in any order. The same computation reproduces the certified optima of Tables 2--3 in all those cells. For larger n we extrapolate: the ratio \(D_{n}\)/\(N_{n}\), where \(N_{n}\) is the number of expressions of Section 2.1 counted with at most one unary minus per node (the grammar allows −−E, so without this convention the number would be infinite), is assumed to keep falling at the rate of the last measured step. The shaded band uses the rate of the step before.

\begin{figure}[t]
\centering
\includegraphics[width=\linewidth]{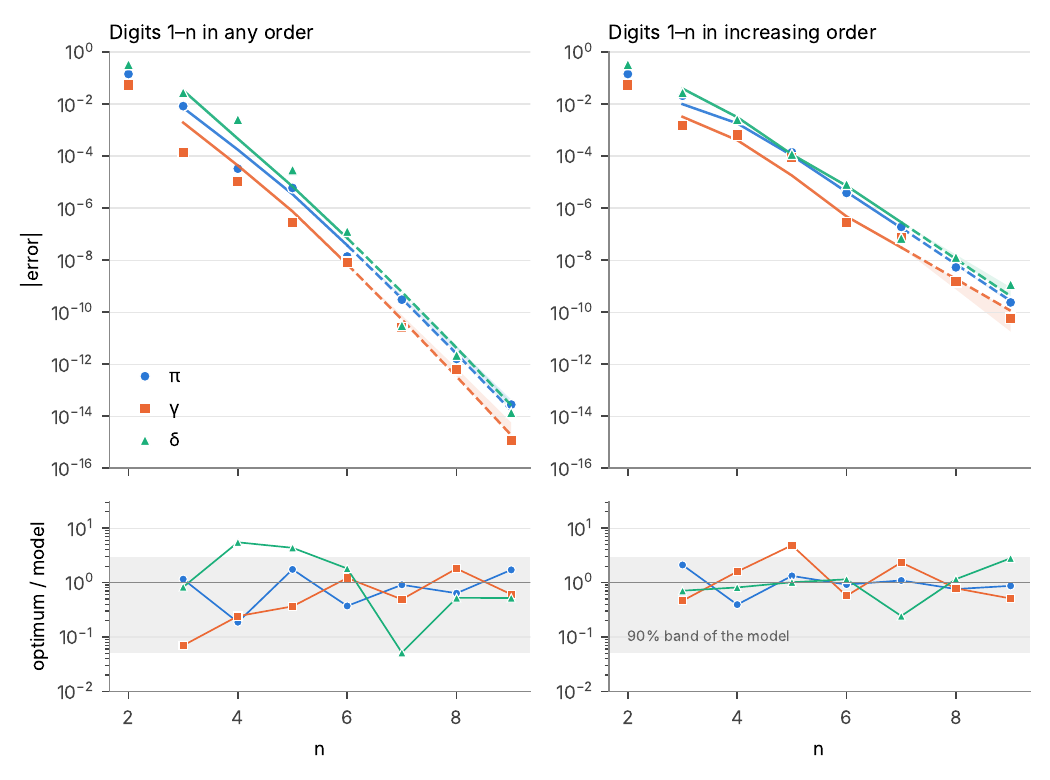}
\caption{Certified optimal error against n (top) and its ratio to the random model (bottom), for π, γ and δ, with the digits in any order (left) and in increasing order (right). Lines: model, solid where \(D_{n}\) is counted exactly and dashed where it is extrapolated, with the extrapolation band shaded. Gray band: 90 \% range of the model.}\label{fig:error}
\end{figure}

Figure~\ref{fig:error} compares the model with the certified optima. In 39 of the 42 cells with n = 3--9, the ratio between the optimum and the model lies in the 90 \% band, against 37.8 expected; the geometric mean of the ratio is 0.82. For an exponential variable it would be \(e^{-γ}\) ≈ 0.56, so in geometric mean the optima are about 1.5 times worse than the model predicts; the cells are not independent, so this is only indicative. The three cells outside the band are δ in any order with n = 4 and 5 (the optimum is 5.5 and 4.4 times worse than the model) and γ in order with n = 5 (4.8 times worse). For n = 9 in any order the model predicts 1.6·\(10^{-14}\) for π, 1.9·\(10^{-15}\) for γ and 2.6·\(10^{-14}\) for δ, against optima of 2.81·\(10^{-14}\), 1.14·\(10^{-15}\) and 1.38·\(10^{-14}\). The optima for π, γ and δ are therefore within the range that a random spread of the available values predicts, and no cell is much better than predicted. The three cells outside the band are all among the 27 with an exact count, and all on the side of worse optima; 15 of the 42 cells rely on extrapolated counts, and the model is too coarse to detect finer structure. The model also accounts for part of γ's apparent luck: for a given \(D_{n}\) the expected error is proportional to Δ = 0.05·c, the counts \(D_{n}\) turn out to be similar for the three constants (4.23, 4.18 and 3.33 million with n = 6 in any order), and γ is the smallest of the three.

\subsection*{9.3 The complex counterpoint}\addcontentsline{toc}{subsection}{9.3 The complex counterpoint}

The restriction to real values is what keeps π in the ``coincidence'' regime. In 2019 Avraham Aizenbud sent Friedman {[}1{]} the 0--9 expression

\(-{(-{8}^{96})}^{.5}\cdot \left({(-1)}^{{4}^{-72}}-{3}^{0}\right)\),

with principal values. With ε = \({4}^{-72}\) = \({2}^{-144}\) and \({(-{8}^{96})}^{.5}\) = i·\({2}^{144}\) = i/ε, its value is −i·(\(e^{iπε}\) − 1)/ε = sin(πε)/ε + i·(1 − cos(πε))/ε. That is π − π\(^{3}\)ε\(^{2}\)/6 + O(ε\(^{4}\)) plus an imaginary part π\(^{2}\)ε/2 + O(ε\(^{3}\)). We checked it with mpmath at 200 digits: the real part differs from π by 1.04·\(10^{-86}\) and the imaginary part is 2.21·\(10^{-43}\), so the value approximates π to about 43 digits, as Friedman reports. Here \({(-1)}^{\varepsilon}-1\) ≈ iπε plays the role of \({(1+1/N)}^{N}\) for e: once complex intermediates are allowed, iπ = log(−1) is a limit of short expressions, and π becomes a problem of the same kind as e. Our rules, like Friedman's tables, exclude it, since \({(-{8}^{96})}^{.5}\) and \({(-1)}^{{4}^{-72}}\) are not real.

\section*{10. Limitations and open problems}\addcontentsline{toc}{section}{10. Limitations and open problems}

\begin{enumerate}
\def\labelenumi{\arabic{enumi}.}
\tightlist
\item
  \textbf{Scope of the certificate.} Theorems 1 and 2 cover π, γ and δ for n ≤ 9 in both modes and e for n ≤ 6 (and n = 7 in order). For e with more digits, optimality is proved only inside the family F (Section 7).
\item
  \textbf{e outside F} (Section 7.5).
\item
  \textbf{0--9 for π, γ and δ.} Not computed: about ten times the cost of 1--9, i.e.~weeks of computation with the present engines.
\item
  \textbf{Other rule sets.} With √, ! or ln, much better formulas are known (MathWorld lists formulas with 17--18 digits); those problems are different and not treated here.
\end{enumerate}

\section*{Reproducibility}\addcontentsline{toc}{section}{Reproducibility}

Code, data, logs and hashes are archived in Zenodo, \url{https://doi.org/10.5281/zenodo.24356239}, and on GitHub, \url{https://github.com/Bendjex/pandigital-certified} (version 3.0). The package contains:

\begin{itemize}
\tightlist
\item
  motors A, B and C, in the optimized versions and in the first versions with which most cells were certified;
\item
  motor X and its tests, every version of its evaluator \texttt{lvl} that checked a record, and the records of n = 9 in any order;
\item
  the verifier of the variant without unary minus;
\item
  the scripts of the proofs for e;
\item
  the Windows packages used for n = 9;
\item
  the certified results of every cell (one JSON line per cell, with the error as an Arb ball: 25 significant digits for π, γ and e, and fewer for the smallest errors of δ, whose ball has radius about \(10^{-30}\) because δ is used with 30 decimals), with the logs and the paths that survived each pre-filter;
\item
  the hashes of the resolution files and of the Windows binaries.
\end{itemize}

A cell is recomputed with two commands (\texttt{build\_\allowbreak{}\allowbreak{}res.sh}, \texttt{certify.py}); the README lists the versions used. Resolution files built with several threads contain the same lines in a different order; their sorted content, and the certificate, are identical, and the published hashes correspond to one thread.

The README also documents the limitations of the code, none of which affects any run of the article. For example:

\begin{itemize}
\tightlist
\item
  motor C's block mode must not be resumed blindly after a crash of the operating system, because its checkpoints are flushed but not synced;
\item
  a test option of motor A (\texttt{MOTORA\_\allowbreak{}\allowbreak{}MAXITEMS}), which limits the enumeration, must not be set;
\item
  to reproduce the records of motor X for n = 9 in any order, the radii of the runs must be used (the package includes them).
\end{itemize}

\section*{Use of AI}\addcontentsline{toc}{section}{Use of AI}

An AI assistant (Claude, by Anthropic) was used under the author's direction to help with the implementation. The author audited and edited all code, proofs, and final text, and every result is independently established by deterministic, reproducible computation.

\section*{Funding}\addcontentsline{toc}{section}{Funding}

No funding was obtained for this work.

\section*{Disclosure statement}\addcontentsline{toc}{section}{Disclosure statement}

The author reports there are no competing interests to declare.

\section*{Data availability}\addcontentsline{toc}{section}{Data availability}

All code, data, logs and certificates supporting the results are openly available in Zenodo at \url{https://doi.org/10.5281/zenodo.24356239}, under the MIT license (code) and CC BY 4.0 (data).

\section*{References}\addcontentsline{toc}{section}{References}

\begin{enumerate}
\def\labelenumi{\arabic{enumi}.}
\tightlist
\item
  E. Friedman, \emph{Problem of the Month (August 2004)}, Math Magic. \url{https://erich-friedman.github.io/mathmagic/0804.html} (accessed 2 November 2026).
\item
  E. W. Weisstein, \emph{e Approximations} and \emph{Pi Approximations}, MathWorld --- A Wolfram Resource. \url{https://mathworld.wolfram.com/eApproximations.html}, \url{https://mathworld.wolfram.com/PiApproximations.html} (accessed 20 November 2026).
\item
  G. Ross, \emph{Pandigital Approximations}, Futility Closet, 12 May 2010, \url{https://www.futilitycloset.com/2010/05/12/pandigital-approximations/}; and \emph{Pandigital Pi}, Futility Closet, 16 March 2017, \url{https://www.futilitycloset.com/2017/03/16/pandigital-pi}.
\item
  Fastfaxr, \emph{New Pandigital Formula for e}, r/math, Reddit, 6 March 2025, \url{https://www.reddit.com/r/math/comments/1j568wk/} (accessed 2 November 2026).
\item
  A. Gottlieb (ed.), Puzzle Corner, \emph{Technology Review} (MIT News section): problem J/A 2 by R. Hess, July/August 2006, pp.~M62--M63; solutions in the November/December 2006 issue and, under ``Better late than never'', in the July/August 2007 issue. Archive: \url{https://cs.nyu.edu/~gottlieb/tr/back-issues/} (accessed 14 November 2026).
\item
  T. Y. Chow, What is a closed-form number?, \emph{Amer. Math. Monthly} 106 (1999), no. 5, 440--448.
\item
  F. Johansson, Arb: efficient arbitrary-precision midpoint-radius interval arithmetic, \emph{IEEE Trans. Comput.} 66 (2017), no. 8, 1281--1292. \href{https://doi.org/10.1109/TC.2017.2690633}{\nolinkurl{https://doi.org/}\allowbreak\mbox{\nolinkurl{10.1109/TC.2017.2690633}}}
\item
  T. Hubrecht, C.-P. Jeannerod, P. Zimmermann, Towards a correctly-rounded and fast power function in binary64 arithmetic, \emph{Proc. 30th IEEE Symposium on Computer Arithmetic (ARITH 2023)}, 111--118, \url{https://doi.org/10.1109/ARITH58626.2023.00028}; and A. Sibidanov, P. Zimmermann, S. Glondu, The CORE-MATH project, \emph{Proc. 29th IEEE Symposium on Computer Arithmetic (ARITH 2022)}, 26--34.
\item
  A. S. Besicovitch, On the linear independence of fractional powers of integers, \emph{J. London Math. Soc.} 15 (1940) 3--6.
\item
  K. Briggs, A precise calculation of the Feigenbaum constants, \emph{Math. Comp.} 57 (1991), no. 195, 435--439. \url{https://doi.org/10.1090/S0025-5718-1991-1079009-6}
\item
  A. McKernon, \emph{pandigital-e}, GitHub repository (2026), \url{https://github.com/andrewtmac/pandigital-e}; and \emph{I Asked an AI to Beat a 23-Year-Old Math Record at 11pm\ldots{}}, Medium, 28 April 2026.
\item
  OEIS Foundation, \emph{A006890: Decimal expansion of Feigenbaum bifurcation velocity}, The On-Line Encyclopedia of Integer Sequences, \url{https://oeis.org/A006890}.
\item
  OEIS Foundation, \emph{A001620: Decimal expansion of Euler's constant (or the Euler-Mascheroni constant), gamma}, The On-Line Encyclopedia of Integer Sequences, \url{https://oeis.org/A001620}.
\item
  The Contest Center, \emph{PI Competition}, \url{http://www.contestcen.com/pi.htm} (page © 2025; accessed 12 November 2026). Its table ``Pandigital: uses all digits once each'' lists the entries of Karnofsky (9.3·\(10^{-11}\)), Ioffe (2.81·\(10^{-11}\)), Vlasii (1.88·\(10^{-12}\) and 4.11·\(10^{-14}\)) and Hess (6.60·\(10^{-13}\)) for 1--9.
\item
  A. P. Goucher, \emph{Closed-form numbers}, Complex Projective 4-Space (blog), 17 October 2020, \url{https://cp4space.hatsya.com/2020/10/17/closed-form-numbers/}.
\item
  D. Bamberger, answer of 14 March 2024 to \emph{An amazing approximation of e}, Mathematics Stack Exchange, question 1945026, \url{https://math.stackexchange.com/questions/1945026} (accessed 2 November 2026).
\item
  L. J. Mordell, On the linear independence of algebraic numbers, \emph{Pacific J. Math.} 3 (1953) 625--630.
\item
  C. L. Siegel, Algebraische Abhängigkeit von Wurzeln, \emph{Acta Arith.} 21 (1972) 59--64.
\item
  A. Burbanks, A. Osbaldestin, J. Thurlby, Rigorous computer-assisted bounds on renormalisation fixed point functions, eigenfunctions, and universal constants, preprint (2021), arXiv:2103.05991, doi:10.48550/arXiv.2103.05991.
\end{enumerate}

\clearpage
\section*{Appendix A. Correctness of the certifier}\addcontentsline{toc}{section}{Appendix A. Correctness of the certifier}

This appendix proves, for each certified cell, that the record is optimal. It follows the order of the computation:

\begin{itemize}
\tightlist
\item
  reductions (§A.2) and rigorous arithmetic (§A.3);
\item
  motor A's tables (§A.4) and window queries (§A.5);
\item
  exact verification with motor B (§A.6) and the main theorem (§A.7).
\end{itemize}

§A.8 summarizes the independent checks, which are not part of the proof. §A.9 fixes the trusted base, and §A.10 links each lemma to the function that implements it.

\subsection*{A.0 Statement}\addcontentsline{toc}{subsection}{A.0 Statement}

A \emph{cell} is a triple (c, n, mode): c is π, e, γ or δ, D = \{1, \ldots, n\}, and the mode is ``in order'' or ``any order''. For δ we use δ₃₀ = 4.669201609102990671853203820466, with an uncertainty of ±\(10^{-30}\) (Section 2.2).

\textbf{Main theorem (Theorem 1 of the text).} For each cell of Tables 2--3, let m* = inf\{\textbar val E − c\textbar{} : E ∈ R(D)\}. Then:

\begin{enumerate}
\def\labelenumi{\arabic{enumi}.}
\tightlist
\item
  the infimum is attained, and m* lies in the ball given in the results file (25 significant digits for π, γ and e; for δ the radius of the ball is about \(10^{-30}\), from the 30 decimals of δ);
\item
  the values v ∈ V(D) with \textbar v − c\textbar{} = m* are exactly those listed;
\item
  the formula given with each of those values is in R(D) and attains it.
\end{enumerate}

This assumes the trusted base of §A.9 and that the computation of the cell ended with 0 unresolved cases.

\subsection*{A.1 Definitions}\addcontentsline{toc}{subsection}{A.1 Definitions}

\begin{itemize}
\item
  \textbf{Expressions.} For a nonempty S ⊆ D, Expr(S) is the set of expressions of the formal rules (Section 2) that use each digit of S exactly once. They are:

  \begin{itemize}
  \tightlist
  \item
    literals: concatenated digits, with an optional decimal point;
  \item
    −E;
  \item
    (\(E_{1}\) op \(E_{2}\)), with op one of +, −, ×, / and \^{}.
  \end{itemize}

  In the in-order mode the digits appear in increasing order: S is an interval and, in (\(E_{1}\) op \(E_{2}\)), \(E_{1}\) is over the left block.
\item
  \textbf{Value.} val(E) is the exact real number. E is \emph{valid} if no subexpression divides by 0, raises 0 to an exponent ≤ 0, or raises a negative base to a non-integer exponent.
\item
  \textbf{Range.} E is \emph{in range} if every subexpression E′ has val(E′) = 0 or \(10^{-300}\) ≤ \textbar val(E′)\textbar{} ≤ \(10^{300}\). We write:

  \begin{itemize}
  \tightlist
  \item
    R(S) for the set of valid, in-range expressions;
  \item
    V(S) = \{val E : E ∈ R(S)\};
  \item
    V\(^{+}\)(S) = \{\textbar x\textbar{} : x ∈ V(S)\}.
  \end{itemize}
\item
  \textbf{Doubles.} RN(z) is the double nearest to z (ties to even), and pred(x), succ(x) are the previous and next doubles. ulp(r) = \(2^{e-53}\) if \(2^{e-1}\) ≤ r \textless{} 2\^{}e. We write τ = \(2^{-10}\).
\item
  \textbf{Classes.} A class κ describes a set H(κ) ⊆ {[}0, ∞) containing values of V\(^{+}\)(S). There are three kinds:

  \begin{itemize}
  \tightlist
  \item
    \emph{ordinary} (r, k), with r \textgreater{} 0 a double and 0 ≤ k ≤ 254. If k = 0, H = \{r\}. If k ≥ 1, H = {[}max(0, pred(r − ρ)), succ(r + ρ){]} with ρ = \(2^{k-1}\)·ulp(r).
  \item
    \emph{rational} p/q, a reduced fraction with \textbar p\textbar, q ≤ \(2^{53}\): H = \{\textbar p\textbar/q\}.
  \item
    \emph{anchored} (q, {[}μ₁, μ₂{]}), with q \textgreater{} 0 rational as above and μ₁ ≤ μ₂ doubles such that 0 ∉ {[}μ₁, μ₂{]} ⊂ (−τ, τ): H = \{q·\(e^{μ}\) : μ ∈ {[}μ₁, μ₂{]}\}. It holds values that are closer to a simple rational than a double can resolve, such as 1 + \(10^{-80}\) or 0.2 + 0.\(1^{81}\).
  \end{itemize}

  Each class carries an enclosure {[}lo, hi{]} ⊇ H and, if lo \textgreater{} 0, an enclosure {[}ℓ₁, ℓ₂{]} ⊇ ln H, both computed outward (§A.3). For a rational class p/q with \textbar p/q − 1\textbar{} \textless{} \(2^{-9}\), the enclosure of ln H is also taken from Lemma 8 applied to u = (p − q)/q, which is much tighter near 1.
\item
  \textbf{Singular.} An intermediate result is singular in three cases:

  \begin{itemize}
  \tightlist
  \item
    Z: its enclosure contains 0 without being the exact zero;
  \item
    U: it contains 1 without being the exact one, or it is anchored with 0 ∈ {[}μ₁, μ₂{]};
  \item
    W: its relative width exceeds \(2^{-20}\).
  \end{itemize}
\item
  \textbf{Windows.} A window W has an interval {[}a, b{]} of doubles, 0 ≤ a ≤ b. Optionally it also has a logarithmic interval {[}\(L_{1}\), \(L_{2}\){]} and an anchored form A·e\footnote{ν₁, ν₂}, with A \textgreater{} 0 rational and {[}ν₁, ν₂{]} ⊂ (−τ, τ). W describes the set D(W) of the x ∈ {[}a, b{]} that satisfy the conditions present: ln x ∈ {[}\(L_{1}\), \(L_{2}\){]} if x \textgreater{} 0, and x ∈ A·e\footnote{ν₁, ν₂}. If there is a logarithmic interval and a = 0, then \(L_{1}\) = −∞.
\item
  \textbf{Paths.} A path is a leaf or a quadruple (\(S_{1}\), κ, f, π).

  \begin{itemize}
  \tightlist
  \item
    A leaf is a class of T(S), a literal of S, or ``constant k''.
  \item
    In the quadruple, \(S_{1}\) ⊆ S, κ is a class of T(\(S_{1}\)), f is one of the seven families of §A.5 and π is a path over \(S_{2}\) = S ∖ \(S_{1}\).
  \end{itemize}
\item
  \textbf{Covering.}

  \begin{itemize}
  \tightlist
  \item
    A class leaf κ covers E if \textbar val E\textbar{} ∈ H(κ). A literal leaf covers that literal.
  \item
    The quadruple (\(S_{1}\), κ, f, π) covers E = (\(E_{1}\) op \(E_{2}\)) if three things hold: the part of E over \(S_{1}\) has absolute value s ∈ H(κ); \textbar val E\textbar{} = f(s, l), with l the absolute value of the other part; and π covers that other part.
  \item
    If f(s, ·) is constant, equal to k, the leaf ``constant k'' covers any other part.
  \end{itemize}
\end{itemize}

\subsection*{A.2 Reductions}\addcontentsline{toc}{subsection}{A.2 Reductions}

\textbf{Lemma 1 (symmetry).} V(S) = −V(S).

\emph{Proof.} If E ∈ R(S), then −E is valid and uses the same digits in the same order. Its subexpressions are those of E plus the root, whose absolute value is \textbar val E\textbar. Hence −E ∈ R(S), and val(−E) = −val E. ∎

\textbf{Lemma 2 (positive bases).} Let R\(^{+}\)(S) be the set of E ∈ R(S) in which every power \({x}^{y}\) has x \textgreater{} 0, or x = 0 and y \textgreater{} 0. For every valid expression E over S (in range or not) there is a valid E′ with the same digits in the same order, val E′ = val E, every power of E′ with a positive or zero base, and every node of E′ with the absolute value of a node of E. In particular, if E ∈ R(S) then E′ ∈ R\(^{+}\)(S).

\emph{Proof.} By induction on the size of E, replacing from the inside out; by induction val \(E_{1}\)′ = val \(E_{1}\) and val \(E_{2}\)′ = val \(E_{2}\). If E = (\(E_{1}\))\^{}(\(E_{2}\)) with val \(E_{1}\) \textless{} 0, validity forces val \(E_{2}\) = k to be an integer, and then val E = (−1)\^{}k·\textbar val \(E_{1}\)\textbar\^{}k. Take E′ = ((−\(E_{1}\)′)\^{}(\(E_{2}\)′)) if k is even and E′ = −((−\(E_{1}\)′)\^{}(\(E_{2}\)′)) if k is odd; in both cases val E′ = val E. The new base −\(E_{1}\)′ has value \textbar val \(E_{1}\)\textbar{} \textgreater{} 0, and every new node has the absolute value of a node of E. In the other cases E′ = (\(E_{1}\)′ op \(E_{2}\)′), with the same value as E. ∎

\textbf{Lemma 3 (forms).} Let E = (\(E_{1}\) op \(E_{2}\)) ∈ R\(^{+}\)(S), with \(E_{1}\) over \(S_{1}\), s = \textbar val \(E_{1}\)\textbar{} and l = \textbar val \(E_{2}\)\textbar. Then \textbar val E\textbar{} is one of

\[
s + l,\quad |s - l|,\quad s\,l,\quad s/l,\quad l/s,\quad s^{\pm l},\quad l^{\pm s}
\]

In the in-order mode only s + l, \textbar s − l\textbar, s·l, s/l and \(s^{±l}\) occur, because \(E_{1}\) is the left block. Conversely, let g(s, l) be one of those values, with s = \textbar val \(E_{1}\)\textbar, l = \textbar val \(E_{2}\)\textbar{} and \(E_{i}\) ∈ R\(^{+}\)(\(S_{i}\)). If g(s, l) is defined and in range, there is E ∈ R\(^{+}\)(S) with \textbar val E\textbar{} = g(s, l).

\emph{Proof.} With val \(E_{1}\) = ±s and val \(E_{2}\) = ±l we have \textbar±s ± l\textbar{} ∈ \{s + l, \textbar s − l\textbar\}, \textbar(±s)(±l)\textbar{} = s·l and \textbar(±s)/(±l)\textbar{} = s/l. In R\(^{+}\) the base of a power is positive or zero, hence s, and the exponent is ±l: \textbar val E\textbar{} = \(s^{±l}\). For the converse, the sign of each part is chosen with unary minus: −\(E_{1}\) and −\(E_{2}\) stay in R\(^{+}\) and use the same digits in the same order (Lemma 1). With those signs every form is obtained, and the new subexpressions are ±s, ±l or ±g(s, l), all in range. ∎

\subsection*{A.3 Rigorous arithmetic}\addcontentsline{toc}{subsection}{A.3 Rigorous arithmetic}

\textbf{Lemma 4 (rounding).} The two parts of the lemma are:

\begin{enumerate}
\def\labelenumi{\arabic{enumi}.}
\tightlist
\item
  The IEEE 754 binary64 operations +, −, × and /, in round-to-nearest mode, return RN of the exact result, and CORE-MATH's cr\_\allowbreak{}\allowbreak{}pow(x, y) returns RN\(({x}^{y})\).
\item
  If r = RN(z), then pred(r) ≤ z ≤ succ(r).
\end{enumerate}

Hence pred(RN(z)) and succ(RN(z)) are lower and upper bounds of z, and applied to every operation they give \emph{outward} enclosures.

\emph{Proof.} Part 1 is the IEEE 754 standard plus the CORE-MATH theorem {[}8{]}. The program is compiled with \texttt{-ffp-contract=off} and without \texttt{-ffast-math}, so every operation written in the code is an IEEE operation. For part 2, \textbar RN(z) − z\textbar{} ≤ ulp/2 in the binade of z, and the gaps r − pred(r) and succ(r) − r are at least that. This also holds at a binade boundary, where the lower gap is half, and in the subnormal range. ∎

\textbf{Lemma 5 (corners).} On a box {[}\(x_{1}\), \(x_{2}\){]} × {[}\(y_{1}\), \(y_{2}\){]}, the functions x + y, x − y, x·y (with x, y ≥ 0) and x/y (with \(y_{1}\) \textgreater{} 0) attain their minimum and maximum at corners of the box. So does \({x}^{y}\) if \(x_{1}\) \textgreater{} 0.

\emph{Proof.} The first four are monotone in each variable, with a fixed sign on the box. For the power, \({x}^{y}\) = exp(y·ln x), and y·ln x is bilinear in (ln x, y), so its extrema on the rectangle {[}ln \(x_{1}\), ln \(x_{2}\){]} × {[}\(y_{1}\), \(y_{2}\){]} are at the corners; the exponential is increasing. ∎

\textbf{Corollary.} If x ∈ {[}\(x_{1}\), \(x_{2}\){]}, y ∈ {[}\(y_{1}\), \(y_{2}\){]}, and m and M are the minimum and maximum of the RN of x op y at the four corners, then x op y ∈ {[}pred(m), succ(M){]}. This follows from the monotonicity of RN and Lemma 4.

\textbf{Lemma 6 (exact rationals).} Let p/q be a reduced fraction with \textbar p\textbar, q ≤ \(2^{53}\). Then:

\begin{itemize}
\tightlist
\item
  p and q are exact doubles, and RN(p/q) is obtained with a single division;
\item
  p/q is a double if and only if q is a power of 2.
\end{itemize}

Operations on rationals are done with 128-bit integers, which do not overflow, and are accepted only if the reduced result again satisfies the bounds. The power (p/q)\^{}(u/v), with v ≤ 64 and \textbar u\textbar{} ≤ 4096, is accepted as rational only if p and q are exact v-th powers; this is checked by raising the candidate in integers.

\emph{Proof.} IEEE division is correctly rounded. A reduced fraction p/q with \textbar p\textbar{} ≤ \(2^{53}\) is a double if and only if q is a power of 2; with q ≤ \(2^{53}\) there is no underflow. If p and q are coprime, (p/q)\^{}(1/v) is rational if and only if p and q are v-th powers. The candidate given by the floating root is only an estimate, so the three neighbouring integers are checked. ∎

\textbf{Lemma 7 (exponential and logarithm).} Let \(λ_{1}\) = 0x1.71547652b82fep+0 and \(λ_{2}\) = succ(\(λ_{1}\)); then \(λ_{1}\) \textless{} lo\(g_{2}\) e \textless{} \(λ_{2}\) (checked with mpmath). Define:

\begin{itemize}
\tightlist
\item
  exp\_\allowbreak{}\allowbreak{}hi(z) = succ(cr\_\allowbreak{}\allowbreak{}pow(2, succ(RN(z·λ)))), with λ = \(λ_{2}\) if z ≥ 0 and λ = \(λ_{1}\) if z \textless{} 0;
\item
  exp\_\allowbreak{}\allowbreak{}lo(z) = pred(cr\_\allowbreak{}\allowbreak{}pow(2, pred(RN(z·λ′)))), with λ′ = \(λ_{1}\) if z ≥ 0 and λ′ = \(λ_{2}\) if z \textless{} 0, or 0 if cr\_\allowbreak{}\allowbreak{}pow returns 0;
\item
  ln\_\allowbreak{}\allowbreak{}lo(x) = the first g of the sequence \(g_{0}\) = estimate, \(g_{i+1}\) = \(g_{i}\) − 2\^{}i·h (h \textgreater{} 0 a fixed small step) with exp\_\allowbreak{}\allowbreak{}hi(g) ≤ x (−∞ if there is none);
\item
  ln\_\allowbreak{}\allowbreak{}hi(x) = the same with exp\_\allowbreak{}\allowbreak{}lo(g) ≥ x.
\end{itemize}

Then exp\_\allowbreak{}\allowbreak{}lo(z) ≤ \(e^{z}\) ≤ exp\_\allowbreak{}\allowbreak{}hi(z) and ln\_\allowbreak{}\allowbreak{}lo(x) ≤ ln x ≤ ln\_\allowbreak{}\allowbreak{}hi(x).

\emph{Proof.} \(e^{z}\) = \(2^{z·log_{2} e}\). If z ≥ 0, z·\(λ_{2}\) ≥ z·lo\(g_{2}\) e; if z \textless{} 0, z·\(λ_{1}\) ≥ z·lo\(g_{2}\) e. Rounding up is preserved because 2\^{}w is increasing and cr\_\allowbreak{}\allowbreak{}pow is correctly rounded (Lemma 4). The lower bound is the same. For the logarithm: if exp\_\allowbreak{}\allowbreak{}hi(g) ≤ x, then \(e^{g}\) ≤ x and g ≤ ln x. The initial estimate affects only speed, never validity. ∎

\textbf{Lemma 8 (small arguments).} If \textbar m\textbar{} ≤ 1, then m ≤ \(e^{m}\) − 1 ≤ m + m\(^{2}\). If u ≥ −1/2, then u − u\(^{2}\) ≤ ln(1 + u) ≤ u. On intervals, with outward rounding:

\begin{itemize}
\tightlist
\item
  expm1({[}\(m_{1}\), \(m_{2}\){]}) ⊆ {[}\(m_{1}\), succ(\(m_{2}\) + succ(\(m_{2}\)\(^{2}\))){]};
\item
  log1p({[}\(u_{1}\), \(u_{2}\){]}) ⊆ {[}pred(\(u_{1}\) − succ(\(u_{1}\)\(^{2}\))), \(u_{2}\){]}.
\end{itemize}

\emph{Proof.} \(e^{m}\) − 1 − m = Σ\_\allowbreak{}\allowbreak{}\{k≥2\} \(m^{k}\)/k!. If 0 ≤ m ≤ 1, the sum is between 0 and m\(^{2}\)(e − 2) ≤ m\(^{2}\). If −1 ≤ m \textless{} 0, the series alternates with decreasing terms and lies between 0 and m\(^{2}\)/2. For the logarithm, the upper bound is concavity. Let h(u) = ln(1 + u) − u + u\(^{2}\). Then h(0) = 0 and h′(u) = u(1 + 2u)/(1 + u), which is ≤ 0 on {[}−1/2, 0{]} and ≥ 0 on {[}0, ∞); hence h ≥ 0 on {[}−1/2, ∞). Both functions are increasing, so on an interval the endpoints suffice, and the roundings go outward (Lemma 4). ∎

\subsection*{A.4 Motor A's tables}\addcontentsline{toc}{subsection}{A.4 Motor A's tables}

\textbf{Construction.} The tables T(S) are built for \textbar S\textbar{} ≤ TBL, in increasing order of \textbar S\textbar; TBL is a fixed parameter (5 in the runs with n = 9), and the guarantees do not depend on its value. T(S) gathers three things:

\begin{enumerate}
\def\labelenumi{\arabic{enumi}.}
\tightlist
\item
  the literals of S, as rational classes;
\item
  the \emph{derivations} (κ₁, κ₂, f): for every allowed split S = \(S_{1}\) ∪ \(S_{2}\), every κᵢ ∈ T(\(S_{i}\)) and every form f of Lemma 3 allowed in the mode, with x the value over \(S_{1}\) and y the value over \(S_{2}\). The forms are x + y, \textbar x − y\textbar, x·y, x/y, y/x, \({x}^{y}\), \({x}^{-y}\), \(y^{x}\) and \(y^{-x}\). In order, \(S_{1}\) is the left block and only the five forms of the shape x op y occur, counting \({x}^{y}\) and \({x}^{-y}\) as two;
\item
  the items that motor B returns for every singular derivation (Lemma 10).
\end{enumerate}

When merging, ordinary classes with the same key are joined (keeping the larger code), and identical rational or anchored classes are identified. In any order every unordered split is visited once, with the nine forms; in order, once per cut point.

\textbf{Lemma 9 (derivations).} Let κ₁ and κ₂ be non-singular classes, f a form and P = \{(x, y) ∈ H(κ₁) × H(κ₂) : f(x, y) is defined\}. The derivation returns ``undefined'', ``out of range'', ``singular'' or a class κ, with these guarantees:

\begin{itemize}
\tightlist
\item
  ``undefined'' only if P = ∅;
\item
  ``out of range'' only if no f(x, y) with (x, y) ∈ P is 0 or lies in {[}\(10^{-300}\), \(10^{300}\){]};
\item
  if it returns κ, then f(P) ⊆ H(κ).
\end{itemize}

The rules, in the order in which they apply:

{\small
\begin{longtable}[]{@{}>{\raggedright\arraybackslash}p{(\linewidth - 2\tabcolsep) * \real{0.5000}}>{\raggedright\arraybackslash}p{(\linewidth - 2\tabcolsep) * \real{0.5000}}@{}}
\toprule\noalign{}
Case & Rule \\
\midrule\noalign{}
\endhead
\bottomrule\noalign{}
\endlastfoot
Exact zeros and ones & x + 0 = x, 0·y = 0, 0/y = 0, y/1 = y, a\(^{0}\) = 1 if a ≠ 0, \({0}^{y}\) = 0 if y \textgreater{} 0, \({1}^{y}\) = 1, a\(^{1}\) = a. Undefined: x/0, \(0^{0}\) and \({0}^{-y}\) \\
Two rationals & The exact rational result, if it fits (Lemma 6) \\
Sum or difference with a rational or anchored operand & Write x = \(Q_{1}\) + δ₁ and y = \(Q_{2}\) + δ₂ with Q exact rational: δ = 0 for a rational, δ ∈ q·expm1({[}μ₁, μ₂{]}) for an anchored class, and Q = 0 with δ ∈ {[}lo, hi{]} for an ordinary one. Then Q = \(Q_{1}\) ± \(Q_{2}\) exactly and δ = δ₁ ± δ₂ outward. If Q ≠ 0 and \textbar δ\textbar{} ≤ \textbar Q\textbar·\(2^{-12}\), the result is anchored at \textbar Q\textbar{} \\
Product or quotient of rationals or anchored classes & \(q_{1}\)e\^{}μ₁ · \(q_{2}\)e\^{}μ₂ = (\(q_{1}\)\(q_{2}\))·\(e^{μ_{1}+μ_{2}}\), and similarly for the quotient (a rational has μ = 0) \\
Anchored base & (q \(e^{μ}\))\^{}y = \(q^{y}\)·\(e^{yμ}\) if y and \(q^{y}\) are rational; otherwise from the logarithm y·ln(q \(e^{μ}\)) \\
Rational base, anchored exponent & r\^{}(q \(e^{μ}\)) = \(r^{q}\)·e\^{}(q(\(e^{μ}\) − 1)·ln r), if \(r^{q}\) is rational \\
General case & Corollary of Lemma 5 on the enclosures. In a power the base may not contain 0 or 1 (otherwise singular); if the result lies within τ of 1, it is redone from ln x = y·ln a \\
From a logarithm {[}ℓ₁, ℓ₂{]} & If {[}ℓ₁, ℓ₂{]} ⊂ (−τ, τ), the anchored class (1, {[}ℓ₁, ℓ₂{]}); otherwise the interval {[}exp\_\allowbreak{}\allowbreak{}lo(ℓ₁), exp\_\allowbreak{}\allowbreak{}hi(ℓ₂){]} \\
Anchored (q, {[}μ₁, μ₂{]}) & Singular U if 0 ∈ {[}μ₁, μ₂{]}. If {[}μ₁, μ₂{]} is not inside (−τ, τ), it becomes the interval {[}pred(\(q_{1}\)·exp\_\allowbreak{}\allowbreak{}lo(μ₁)), succ(\(q_{2}\)·exp\_\allowbreak{}\allowbreak{}hi(μ₂)){]}, with {[}\(q_{1}\), \(q_{2}\){]} the enclosure of q \\
Interval {[}lo, hi{]} & Out of range if hi \textless{} pred(\(10^{-300}\)) or lo \textgreater{} succ(\(10^{300}\)). Singular Z if lo ≤ 0. Otherwise key r ∈ {[}lo, hi{]} and the smallest code k with H(r, k) ⊇ {[}lo, hi{]}; singular U if H(r, k) contains 1, and W if its width exceeds r·\(2^{-20}\) \\
\end{longtable}
}

\emph{Proof.} Each rule is an identity valid for all pairs of P, followed by outward enclosures (Lemmas 4--8). The delicate points are these:

\begin{itemize}
\tightlist
\item
  \textbf{Exactness.} The zero and one flags are set only on rational classes, with H = \{0\} or \{1\}, so the rules of the first row are exact.
\item
  \textbf{Sums with an anchor.} If Q ≠ 0 and \textbar δ\textbar{} ≤ \textbar Q\textbar·\(2^{-12}\), then Q + δ = Q(1 + u) with u = δ/Q and \textbar u\textbar{} ≤ \(2^{-12}\). Hence \textbar Q + δ\textbar{} = \textbar Q\textbar·e\^{}(ln(1+u)), and Lemma 8 encloses ln(1 + u). Otherwise \textbar Q + δ\textbar{} lies in {[}lo, hi{]} if lo \textgreater{} 0, in {[}−hi, −lo{]} if hi \textless{} 0, and else it is singular Z.
\item
  \textbf{Powers.} A non-singular class does not contain 1 (U), so in the general case the base satisfies \(x_{2}\) \textless{} 1 or \(x_{1}\) \textgreater{} 1, and the exponent does not contain 0. The passage through the logarithm uses ln x = y·ln a, valid for a \textgreater{} 0.
\item
  \textbf{Code.} k = code\_\allowbreak{}\allowbreak{}for(r, lo, hi) is the smallest integer with \(2^{k-1}\)·ulp(r) ≥ max(r − lo, hi − r), rounded up. Hence H(r, k) ⊇ {[}lo, hi{]} whenever hi is finite. If k \textgreater{} 254 were needed, the relative width would be enormous and rule W marks it singular. The case hi = +∞ (an upper bound that overflows, possible only for a power whose exponent multiplies a logarithm of large relative width) is not covered by this argument: code\_\allowbreak{}\allowbreak{}for would then return a class that does not contain {[}lo, hi{]}. The independent review of this work found it by reading the code, and counted the calls with an instrumented motor A over every table used by the cells (up to 5 digits; n ≤ 9 in order, and n = 6, 7, 8 and 9 in any order: about 1.5·\(10^{9}\) calls). None had hi = +∞, so the lemma holds for every derivation of those tables. The published code is unchanged, and the README describes the issue and the guard that would remove it.
\item
  \textbf{Range.} pred(\(10^{-300}\)) \textless{} \(10^{-300}\) and succ(\(10^{300}\)) \textgreater{} \(10^{300}\), so a derivation is discarded only if all its values are out of range. ∎
\end{itemize}

\textbf{Theorem A (complete tables).} Let \textbar S\textbar{} ≤ TBL. For every E ∈ R\(^{+}\)(S) there is κ ∈ T(S) with \textbar val E\textbar{} ∈ H(κ). No class of T(S) is singular.

\emph{Proof.} By induction on \textbar S\textbar{} and, within each S, on the size of E:

\begin{itemize}
\tightlist
\item
  If E is a literal, its rational class is in T(S).
\item
  If E = −E′, it has the absolute value of the smaller E′.
\item
  If E = (\(E_{1}\) op \(E_{2}\)), its parts are in R\(^{+}\)(\(S_{1}\)) and R\(^{+}\)(\(S_{2}\)), with \textbar{}\(S_{i}\)\textbar{} \textless{} \textbar S\textbar. By induction their absolute values lie in H(κ₁) and H(κ₂) for non-singular κᵢ ∈ T(\(S_{i}\)). By Lemma 3, \textbar val E\textbar{} = f(x, y) for some form f, with (x, y) ∈ P because E is valid. By Lemma 9 the derivation (κ₁, κ₂, f) is neither ``undefined'' nor ``out of range'', the latter because E is in range. If it is a class, that class contains \textbar val E\textbar. If it is singular, one of the items of Lemma 10 contains it.
\end{itemize}

Merging only enlarges sets. Singular results are never stored, and if an item from B were singular the computation stops. ∎

\textbf{Lemma 10 (resolution of singular derivations).} Let (κ₁, κ₂, f) be a singular derivation of T(S). Motor B computes exactly \(X_{i}\) = V\(^{+}\)(\(S_{i}\)) ∩ H(κᵢ) (Lemma 14). For every pair with f defined and possibly in range, it writes one item:

\begin{itemize}
\tightlist
\item
  the exact rational, if it fits;
\item
  an anchored class (q, enclosure of μ), if the value is q·\(e^{μ}\) with \textbar μ\textbar{} \textless{} τ and μ ≠ 0, where q is the rational part of an algebraic value;
\item
  an anchored class at 1, if \textbar ln v\textbar{} \textless{} τ and ln v ≠ 0;
\item
  an outward enclosure {[}lo, hi{]} with \textbar ln v\textbar{} ≥ τ, containing neither 0 nor 1.
\end{itemize}

Every f(x, y) with x ∈ \(X_{1}\), y ∈ \(X_{2}\), defined and in range, lies in the set H of some item. If B cannot separate a value from 0 or 1 at 16 384 bits, it reports it unresolved and the computation stops. Motor A reads the items with the classification rules of Lemma 9.

\emph{Proof.} \(X_{i}\) contains all relevant values (Lemma 14). Arb balls contain the exact value, and their conversion to doubles is outward. The test \textbar ln v\textbar{} ≥ τ excludes 1, and the separation from 0 is checked before it. ∎

\subsection*{A.5 Window queries}\addcontentsline{toc}{subsection}{A.5 Window queries}

\textbf{Construction.} query(S, W) returns a list of paths:

\begin{itemize}
\tightlist
\item
  If \textbar S\textbar{} ≤ TBL, the leaves of T(S) that pass the filters of Lemma 12.
\item
  If \textbar S\textbar{} \textgreater{} TBL, the literal leaves of S that pass those filters, plus the paths obtained as follows:

  \begin{itemize}
  \tightlist
  \item
    every allowed split S = \(S_{1}\) ∪ \(S_{2}\) is visited with \(S_{1}\) the smaller part: in order, each cut point, with \(S_{1}\) the shorter block; in any order, each split with \textbar{}\(S_{1}\)\textbar{} ≤ \textbar S\textbar/2, once;
  \item
    for every class κ ∈ T(\(S_{1}\)) (ordinary, rational or anchored) and every family f allowed on the side of \(S_{1}\), the windows P\_\allowbreak{}\allowbreak{}f(κ, W) of Lemma 11 are computed, and query(\(S_{2}\), W′) is called for each W′ ∈ P\_\allowbreak{}\allowbreak{}f(κ, W). In any order all seven families are allowed. In order, +, \textbar−\textbar, ×, s/l and \(s^{±l}\) are allowed if \(S_{1}\) is the left block, and +, \textbar−\textbar, ×, l/s and \(l^{±s}\) if it is the right one;
  \item
    if f(s, ·) is constant, equal to k, and k ∈ D(W), the leaf ``constant k'' is added.
  \end{itemize}
\end{itemize}

As \textbar S\textbar{} ≤ 9, \textbar{}\(S_{1}\)\textbar{} ≤ 4 ≤ TBL: the smaller part is always in a table. The initial call is query(D, \(W_{0}\)), with \(W_{0}\) = {[}a, b{]}, a ≤ c − ε and b ≥ c + ε (§A.7).

\textbf{Lemma 11 (preimages).} Let κ ∈ T(\(S_{1}\)) be non-singular, s ∈ H(κ), l ∈ V\(^{+}\)(\(S_{2}\)) and f a family with f(s, l) defined, in range and in D(W). If f(s, ·) is not constant, then l ∈ D(W′) for some W′ ∈ P\_\allowbreak{}\allowbreak{}f(κ, W).

Notation:

\begin{itemize}
\tightlist
\item
  {[}\(s_{1}\), \(s_{2}\){]} is the enclosure of κ, and {[}ℓ₁, ℓ₂{]} ⊇ ln H(κ).
\item
  {[}\(L_{1}\), \(L_{2}\){]} is the logarithmic interval of W: \(L_{1}\) = ln\_\allowbreak{}\allowbreak{}lo(a) (or −∞ if a = 0) and \(L_{2}\) = ln\_\allowbreak{}\allowbreak{}hi(b), intersected with that of W if it has one.
\item
  If κ is rational or anchored, s = \(Q_{s}\) + \(δ_{s}\) as in Lemma 9. If W has no anchored form, x = 0 + {[}a, b{]}; if it has one, x = A + A·expm1({[}ν₁, ν₂{]}).
\end{itemize}

{\small
\begin{longtable}[]{@{}>{\raggedright\arraybackslash}p{(\linewidth - 4\tabcolsep) * \real{0.1559}}>{\raggedright\arraybackslash}p{(\linewidth - 4\tabcolsep) * \real{0.3301}}>{\raggedright\arraybackslash}p{(\linewidth - 4\tabcolsep) * \real{0.5140}}@{}}
\toprule\noalign{}
Family & l as a function of x = f(s, l) & Windows for l \\
\midrule\noalign{}
\endhead
\bottomrule\noalign{}
\endlastfoot
s + l & l = x − s & If κ is rational or anchored: l = (\(Q_{W}\) − \(Q_{s}\)) + (\(δ_{W}\) − \(δ_{s}\)), an affine window (anchored when \textbar δ\textbar{} is small against the rational part). Otherwise {[}a − \(s_{2}\), b − \(s_{1}\){]} \\
\textbar s − l\textbar{} & l = s + x or l = s − x & Two affine windows, (\(Q_{s}\) + \(Q_{W}\), \(δ_{s}\) + \(δ_{W}\)) and (\(Q_{s}\) − \(Q_{W}\), \(δ_{s}\) − \(δ_{W}\)). Otherwise {[}\(s_{1}\) + a, \(s_{2}\) + b{]} and {[}\(s_{1}\) − b, \(s_{2}\) − a{]} \\
s·l & l = x/s & With W and κ anchored or rational: (A/q)·\(e^{[ν_{1}, ν_{2}] - [μ_{1}, μ_{2}]}\). If a \textgreater{} 0: ln l ∈ {[}\(L_{1}\) − ℓ₂, \(L_{2}\) − ℓ₁{]}. If a = 0: {[}0, b/\(s_{1}\){]} \\
s/l & l = s/x & (q/A)·\(e^{[μ_{1}, μ_{2}] - [ν_{1}, ν_{2}]}\). If a \textgreater{} 0: ln l ∈ {[}ℓ₁ − \(L_{2}\), ℓ₂ − \(L_{1}\){]}. If a = 0: {[}\(s_{1}\)/b, ∞). None if b = 0 \\
l/s & l = x·s & (A·q)·\(e^{[ν_{1}, ν_{2}] + [μ_{1}, μ_{2}]}\). If a \textgreater{} 0: ln l ∈ {[}\(L_{1}\) + ℓ₁, \(L_{2}\) + ℓ₂{]}. If a = 0: {[}0, b·\(s_{2}\){]} \\
\(s^{±l}\) & l = \textbar t\textbar{} with \(s^{t}\) = x & t ∈ {[}\(L_{1}\), \(L_{2}\){]}/{[}ℓ₁, ℓ₂{]} (if \(L_{1}\) = −∞, only the bound given by \(L_{2}\)); l = \textbar t\textbar. l = 0 only if 1 ∈ {[}a, b{]}. None if b = 0 \\
\(l^{±s}\) & ln l = ±(ln x)/s & If κ is rational and W anchored: \(A^{±1/s}\)·\(e^{±[ν_{1}, ν_{2}]/s}\), when the root is rational. Otherwise ln l ∈ ±{[}\(L_{1}\), \(L_{2}\){]}/{[}\(s_{1}\), \(s_{2}\){]}. l = 0 only with +s and a = 0. If b = 0, the window \{0\} \\
\end{longtable}
}

The constant families are:

\begin{itemize}
\tightlist
\item
  s = 0 in s·l = 0, s/l = 0 (l \textgreater{} 0) and \(s^{l}\) = 0 (l \textgreater{} 0);
\item
  s = 0 in \(l^{±0}\) = 1 (l \textgreater{} 0);
\item
  s = 1 in \(s^{±l}\) = 1.
\end{itemize}

Before use, the windows are rounded outward and cut to {[}0, succ(\(10^{300}\)){]}. Their lower end is raised to pred(\(10^{-300}\)) where l = 0 is impossible.

\emph{Proof.} For fixed s, the map l ↦ f(s, l) is injective on its domain, and the table writes its inverse with each of the descriptions of D(W):

\begin{itemize}
\tightlist
\item
  \textbf{interval:} the classical formula;
\item
  \textbf{logarithm:} ln l = ln x − ln s and the like, exact identities for positive values;
\item
  \textbf{anchored form:} the rational parts are combined exactly, and the exponents μ and ν are added or subtracted as intervals.
\end{itemize}

The union over s ∈ {[}\(s_{1}\), \(s_{2}\){]} is covered by interval arithmetic at the corners (Lemma 5) with outward rounding (Lemma 4). In \(s^{t}\) = x, ln s has a fixed sign because κ does not contain 1, so t = ln x / ln s. The cut loses nothing, because l ∈ V\(^{+}\)(\(S_{2}\)) is 0 or lies in {[}\(10^{-300}\), \(10^{300}\){]}. The lower end is raised only where l = 0 would give a value outside D(W), or an undefined operation: s/0, \(0^{-s}\), or \(s^{0}\) = 1 with 1 ∉ {[}a, b{]}. ∎

\textbf{Lemma 12 (leaves).} The leaf search of query(S, W) emits every class κ of the table, or of the literal table, with H(κ) ∩ D(W) ≠ ∅.

\emph{Proof.} If an ordinary class (r, k) has an enclosure meeting {[}a, b{]}, then r ∈ {[}a(1 − \(2^{k_{m}-51}\)), b(1 + \(2^{k_{m}-51}\)){]}, with k\_\allowbreak{}\allowbreak{}m the largest code of the table. The reason is that pred and succ add at most 3·ulp(r), and ρ + 3·ulp(r) ≤ r·\(2^{k-51}\) for k ≥ 1. So the binary search by key misses none. An anchored class meeting {[}a, b{]} has q ∈ {[}0.998·a, 1.002·b{]}, because \(e^{±τ}\) ∈ {[}0.99902, 1.00098{]}.

The filters discard κ only when H(κ) ∩ D(W) = ∅:

\begin{itemize}
\tightlist
\item
  the enclosure of κ does not meet {[}a, b{]};
\item
  the exact value r of a class with code 0 equals the anchor A, and 0 ∉ {[}ν₁, ν₂{]};
\item
  κ is the exact 1 and 0 ∉ {[}\(L_{1}\), \(L_{2}\){]};
\item
  the enclosure of ln H(κ) does not meet {[}\(L_{1}\), \(L_{2}\){]};
\item
  κ is anchored at the same anchor as W and {[}μ₁, μ₂{]} ∩ {[}ν₁, ν₂{]} = ∅. ∎
\end{itemize}

\textbf{Theorem B (complete windows).} If E ∈ R\(^{+}\)(S) and \textbar val E\textbar{} ∈ D(W), some path of query(S, W) covers E.

\emph{Proof.} By induction on \textbar S\textbar{} and, within each S, on the size of E:

\begin{itemize}
\tightlist
\item
  If \textbar S\textbar{} ≤ TBL, Theorem A and Lemma 12 suffice.
\item
  If E is a literal, its literal leaf covers it (Lemma 12).
\item
  If E = −E′, apply the argument to E′, which has the same absolute value.
\item
  If E = (\(E_{1}\) op \(E_{2}\)), its split of the digits is one of those the query visits. Let \(S_{1}\) be the smaller part, \(E_{a}\) the part of E over \(S_{1}\) and \(E_{b}\) the other. As \textbar{}\(S_{1}\)\textbar{} ≤ TBL, Theorem A gives κ ∈ T(\(S_{1}\)) with s = \textbar val \(E_{a}\)\textbar{} ∈ H(κ). By Lemma 3, \textbar val E\textbar{} = f(s, l), with l = \textbar val \(E_{b}\)\textbar{} and f a family allowed on the side of \(S_{1}\).

  \begin{itemize}
  \tightlist
  \item
    If f(s, ·) is constant, equal to k, then k = \textbar val E\textbar{} ∈ D(W), and the leaf ``constant k'' covers \(E_{b}\).
  \item
    Otherwise Lemma 11 gives W′ with l ∈ D(W′). By induction query(\(S_{2}\), W′) has a path π covering \(E_{b}\), and then (\(S_{1}\), κ, f, π) covers E. ∎
  \end{itemize}
\end{itemize}

\subsection*{A.6 Motor B: exact verification}\addcontentsline{toc}{subsection}{A.6 Motor B: exact verification}

Motor B works with two kinds of values:

\begin{itemize}
\tightlist
\item
  \textbf{S}, exact algebraic values: finite sums Σ c\_\allowbreak{}\allowbreak{}r·r of normalized radicals r = ∏ \(p^{e_{p}}\), with p prime, e\_\allowbreak{}\allowbreak{}p ∈ (0, 1) rational and nonzero rational coefficients;
\item
  \textbf{N}, expression trees over S values, with +, −, ×, inverse and powers with positive base, evaluated with Arb balls at any precision.
\end{itemize}

\textbf{Lemma 13 (exact arithmetic and sign).}

\begin{enumerate}
\def\labelenumi{\arabic{enumi}.}
\tightlist
\item
  Operations on S values are exact, and their normal form is unique.
\item
  An S value with at least one term is not 0. Its sign is decided with Arb balls of increasing precision, and the process terminates.
\item
  The sign of an N value is decided with Arb up to 8192 bits. If it is not separated from 0, its canonical form is computed: if that is the zero form, the value is exactly 0; otherwise the case is \emph{unresolved}.
\item
  Two values with the same canonical form are equal.
\item
  Arb balls contain the exact value. Comparisons with a window raise the precision (256, 1024, 4096 and 16 384 bits) while the ball meets an endpoint, and a value is discarded only if its ball lies entirely outside.
\end{enumerate}

\emph{Proof.}

\begin{itemize}
\item
  \begin{enumerate}
  \def\labelenumi{(\arabic{enumi})}
  \tightlist
  \item
    These are operations on fractions and rational exponents; two distinct normalized radicals have an irrational quotient, by unique factorization.
  \end{enumerate}
\item
  \begin{enumerate}
  \def\labelenumi{(\arabic{enumi})}
  \setcounter{enumi}{1}
  \tightlist
  \item
    This is the theorem of Besicovitch, Mordell and Siegel {[}9, 17, 18{]}: positive real radicals whose pairwise quotients are irrational are linearly independent over ℚ.
  \end{enumerate}
\item
  (3)--(4) The canonical form writes the value as a sum of coefficient times monomial. Each monomial is a product of powers \(p^{E}\) (p prime), \(q^{E}\) (q a rational that is not factored) or \(B^{E}\) (B a positive normalized sum), whose exponents E are themselves canonical forms. It is obtained by applying only identities valid for positive bases and real exponents: (ab)\^{}Y = \(a^{Y}\) \(b^{Y}\), (\(a^{E}\))\^{}Y = \(a^{EY}\), \(a^{E}\) \(a^{F}\) = \(a^{E+F}\), distributivity and \(p^{k+E}\) = \(p^{k}\) \(p^{E}\). So the form is another expression of the same number, and two identical forms denote the same number. Sums used as bases are kept positive: to invert a sum, its sign is decided first with (2) or with Arb.
\item
  \begin{enumerate}
  \def\labelenumi{(\arabic{enumi})}
  \setcounter{enumi}{4}
  \tightlist
  \item
    This is Arb's guarantee. ∎
  \end{enumerate}
\end{itemize}

\textbf{Lemma 14 (queries of B).} B.query(S, {[}lo, hi{]}) returns pairs (value, expression) that include every v ∈ V\(^{+}\)(S) ∩ {[}lo, hi{]}. Each value is exact, and its expression is an E ∈ Expr(S) with val E = v.

\emph{Proof.} If \textbar S\textbar{} ≤ 3 (or ≤ 4 for the large cells), by complete enumeration: literals and all forms of Lemma 3 over all splits, recursively. If \textbar S\textbar{} \textgreater{} 3, as motor A's query but with exact values:

\begin{itemize}
\tightlist
\item
  for each split with \(S_{1}\) the smaller part, all s ∈ V\(^{+}\)(\(S_{1}\)) are enumerated;
\item
  for each family, the windows for l are computed with Arb balls, outward;
\item
  the recursion goes over \(S_{2}\), filtering as in Lemma 13 (5).
\end{itemize}

The range filter keeps values whose ball meets a bound, so it gives a superset. For a power \({x}^{y}\) whose ball is useless (infinite, or containing 0, as 0.\(6^{7^{256}}\)), it decides with the enclosure of y·ln x against ±300·ln 10. Expressions are written with the signs that make the result positive (for example −(a − b) if a \textless{} b), and a base starting with a unary minus is parenthesized. A comparison that cannot be decided is reported unresolved. ∎

\textbf{Theorem C (verification of paths).} For every path π of query(D, \(W_{0}\)), verify(π) returns a set containing \textbar val E\textbar{} for every E ∈ R\(^{+}\)(D) covered by π with \textbar val E\textbar{} ∈ \(W_{0}\).

\emph{Proof.} B follows π step by step:

\begin{itemize}
\tightlist
\item
  In the step (\(S_{1}\), κ, f) it computes \(X_{1}\) = B.query(\(S_{1}\), H(κ)). For each s ∈ \(X_{1}\) it computes the exact windows for l and recurses on the rest of the path.
\item
  In a class leaf it computes B.query(S, H(κ) ∩ window).
\item
  In a ``constant'' leaf it computes f(s, w) with a witness w, the positive literal of \(S_{2}\).
\end{itemize}

If π covers E = (\(E_{1}\) op \(E_{2}\)), then s = \textbar val \(E_{1}\)\textbar{} ∈ \(X_{1}\) (Lemma 14), and l = \textbar val \(E_{2}\)\textbar{} lies in B's exact window and is covered by the rest of π. By induction the recursive call returns l, and f(s, l) = \textbar val E\textbar{} ∈ \(W_{0}\) passes the filter. The leaves are Lemma 14. ∎

\textbf{Lemma 15 (pre-filter, motor C).} In the large cells, before motor B, each path goes through a C program (motorc). It computes a superset of the values covered by the path inside \(W_{0}\), with the same structure as the verification of Theorem C. It uses Arb balls, and exact rationals (fmpq) where the value is rational. A path is discarded only if that superset is empty. The surviving paths go through a second pass at 2048 bits and, finally, motor B.

The first version already included four improvements in speed and memory; none changes the argument:

\begin{itemize}
\tightlist
\item
  \textbf{Tables at lower precision.} The complete tables (up to 4 digits) are built at 128 bits and record the operation and operands of each value. When a value is used, it is recomputed at the working precision (1024 bits) with the same exact operations; a ball at any precision encloses the same exact value. In this mode only equal rationals are merged as duplicates, because two different values can share a ball.
\item
  \textbf{Windows in double precision.} The preimage windows of the query are computed with intervals of doubles rounded outward. Each IEEE 754 operation rounds to nearest with an error of at most half an ulp, so the doubles before and after the result are rigorous bounds. Logarithms and exponentials are computed with Arb and their bounds rounded outward (floor and ceil). Arb is used directly in two cases: near the ends of the double range, and in sums and differences when \textbar s\textbar{} exceeds the window by a factor \(10^{6}\). This affects only efficiency.
\item
  \textbf{Caches.} The values of a step depend only on its class and are computed once, with identical results. In the leaf, the window of the class is queried and then restricted to the path's. The result is a superset, because every value of the intersection lies in both windows. If an undecided value appears, the intersection is queried directly.
\item
  \textbf{Processes.} Two processes share the tables (fork), and each handles half of the paths.
\end{itemize}

Check: with n = 8 in any order, the paths that survive at the end coincide exactly with those of the previous program for π, γ and δ (43, 44 and 60 paths). With n = 7 they coincide with the paths that motor B finds non-empty (1, 2 and 1). The time drops about 3.8 times.

\emph{Proof.} Every exact value that motor B would compute for the path lies in one of motorc's balls, for four reasons:

\begin{itemize}
\tightlist
\item
  \textbf{Same enumeration.} The forms of Lemma 3 are visited over the same splits, and the windows for l are those of Theorem C, computed with Arb outward.
\item
  \textbf{Rigorous operations.} Ball operations are rigorous, and rational ones exact.
\item
  \textbf{Range filter by excess.} A value is discarded only if its ball lies entirely outside {[}\(10^{-300}\), \(10^{300}\){]} ∪ \{0\}.
\item
  \textbf{Values that could be 0.} If a non-rational ball contains 0, the value may be an exact zero. It is still used as an enclosure, and when dividing by it or raising it to a power it is bounded at the corner: \(x^{e}\) ∈ {[}0, max \(r^{e}\){]} if e \textgreater{} 0, \(x^{e}\) ≥ min \(r^{e}\) if e \textless{} 0, and a/x ≥ a/r, with r the upper bound of the ball. If the bound does not decide, the case is marked unknown and the path is kept.
\end{itemize}

So if motor B found a value in \(W_{0}\), motorc would find it too, and the path is not discarded. Each pass is correct on its own, so their composition is too. ∎

As an empirical check, motorc never discarded a path that motor B found non-empty: in n = 7 any order, no wrong discard in the three constants. In γ with n = 8, of 112 880 paths, 10 054 remain after the 256-bit pass, 160 after the 1024-bit one and 44 after the 2048-bit one. The 3 paths that motor B finds non-empty survive all three passes. All of motorc takes 9 min, against 43 min for motor B alone.

\textbf{Optimized version (November 2026).} The runs with n = 9 in any order used a second version of motors A and C. Every change either is a pure speed-up with byte-identical output, or replaces a superset of the covered values by another superset; an undecided event always keeps the path.

\begin{itemize}
\tightlist
\item
  \textbf{Narrow leaves.} The leaf queries query(S, H(κ) ∩ W) instead of querying the class window H(κ) and filtering afterwards. By Theorem B, query returns a superset of the in-range values of S in the window it receives. The ends of the intersection are taken with the ball comparison that is loosest when undecided. This is the route the program already took when an undecided value appeared.
\item
  \textbf{Arb windows.} When a window is narrower than the resolution of doubles, its preimages are computed with Arb, rounded outward, instead of with doubles. Both are rigorous enclosures of the preimage (Lemma 11).
\item
  \textbf{Empty double windows.} A family is skipped only if no ball of the complete table T(\(S_{2}\)), \textbar{}\(S_{2}\)\textbar{} ≤ 4, meets its preimage window in doubles. That window encloses every l with f(s, l) ∈ W, by Lemmas 4 and 5, outward-rounded exp and ln of Arb, and the bounds of the ball of s. The complete table contains a ball for every in-range value of \(S_{2}\). The lookup uses the slack \(10^{-9}\)\textbar x\textbar{} + \(10^{-300}\) of filter\_\allowbreak{}\allowbreak{}full, which exceeds the radius of the balls (at most \(10^{-12}\)\textbar mid\textbar) plus the rounding of the key. Wide balls and balls that may contain 0 always count as meeting the window. Constant families and sets of more than 4 digits are never skipped.
\item
  \textbf{Grouping.} Several values s are verified together over a hull G of their windows. The verification is monotone in the window, so it returns a ball for every covered value in G. Each s then keeps the balls that may lie in its own window, whose ends are exact.
\item
  \textbf{Blocks.} Motor C splits its input into blocks by byte offset, each line belonging to the block that contains its first byte, and processes the blocks in parallel processes. The decision on a line depends only on the line, and a final join concatenates the blocks in order. The output equals that of a single process byte for byte.
\item
  \textbf{Motor A.} An option omits a line identical, as a string, to one already written for the same work item. Theorem B only needs some path to cover each value, so the set of distinct paths does not change. The fused multiply--add of CORE-MATH's pow became the hardware instruction, which IEEE 754 rounds correctly like the software one. Contraction stays disabled.
\end{itemize}

Checks of the second version:

\begin{itemize}
\tightlist
\item
  n = 8 in any order: motor C keeps exactly the same paths as the first version. Motor A reproduces byte for byte the 887 145 candidate lines from which these cells were certified.
\item
  n = 7: motor B finds no discarded path non-empty.
\item
  n = 9: on a sample of 1 \% of motor A's work items and 0.3 \% of motor C's input, the output equals that of the first version.
\item
  52 cells re-certified end to end give the same errors to 25 digits.
\item
  The Windows binaries are built reproducibly, with stated hashes, and give the same outputs as Linux.
\end{itemize}

\subsection*{A.7 Proof of the main theorem}\addcontentsline{toc}{subsection}{A.7 Proof of the main theorem}

Let e\_\allowbreak{}\allowbreak{}rec be the error of the search record. The certifier takes ε = e\_\allowbreak{}\allowbreak{}rec·(1 + \(10^{-6}\)) + 8·ulp(c), which exceeds the rounding error of the search, and \(W_{0}\) = {[}a, b{]} with a ≤ c − ε and b ≥ c + ε (doubles, outward). If it found no value, it would double ε.

\textbf{Step 1: every nearby value appears.} Let E ∈ R(D) with \textbar val E − c\textbar{} ≤ ε. As ε \textless{} c, val E \textgreater{} 0 and val E ∈ \(W_{0}\). By Lemma 2 there is E′ ∈ R\(^{+}\)(D) with \textbar val E′\textbar{} = val E. By Theorem B some path of query(D, \(W_{0}\)) covers E′. In the cells that use the pre-filter of Lemma 15, that path survives it, because it has a value in \(W_{0}\). By Theorem C motor B returns val E. So the set Vals returned by B contains V(D) ∩ {[}c − ε, c + ε{]}.

\textbf{Step 2: the minimum.} Vals is not empty, because it contains the value of the record. Let v* be the value of Vals with minimal \textbar v* − c\textbar; errors are compared with Arb balls of increasing precision. Then \textbar v* − c\textbar{} ≤ e\_\allowbreak{}\allowbreak{}rec \textless{} ε: in every certified cell the optimum is at most 0.9999995 times the radius of \(W_{0}\) (checked in the results files). A value of V(D) outside {[}c − ε, c + ε{]} has error above ε, and one inside is in Vals by step 1, so inf\{\textbar v − c\textbar{} : v ∈ V(D)\} ≥ \textbar v* − c\textbar. The formula of v* is checked independently:

\begin{itemize}
\tightlist
\item
  the string is parsed;
\item
  every subexpression is evaluated with mpmath at 80 digits;
\item
  the digit rule, validity and range are checked, with their margin.
\end{itemize}

So the formula is in R(D) and v* ∈ V(D). Hence m* = \textbar v* − c\textbar{} and the infimum is attained.

\textbf{Step 3: ties.} The values of Vals that cannot be separated from v* are grouped by exact equality (Lemma 13). If there were two different groups, the results file would say so; in every certified cell there is a single group. Each listed value has a formula in R(D), checked as in step 2. By step 1, the list contains all optimal values of V(D), and only those. ∎

\subsection*{A.8 Independent checks}\addcontentsline{toc}{subsection}{A.8 Independent checks}

They are not part of the proof, but they serve to detect programming errors.

\begin{itemize}
\tightlist
\item
  \textbf{t1: motor B against sympy.} An independent enumeration with sympy, for sets of up to 3 digits, agrees with motor B on every set tested: 0 differences.
\item
  \textbf{t2: query against enumeration.} Motor B's window queries agree with its complete enumeration: 606 windows, 0 differences.
\item
  \textbf{t4: random expressions.} Valid expressions are generated at random and a window is centered on the value of each. Motor A's paths, verified by B, must find that value; in the small configurations they must also find everything B finds directly. Results: 5 digits any order, 40 windows without differences; 7 in order, 26 windows without differences; 8 in order, 28 windows, with the central value found in all.
\item
  \textbf{Agreement with the search v5.} In every cell, the certified optimum is the record of the search, which is a different program.
\item
  \textbf{Check of the formulas.} Every winning formula is evaluated with a parser and mpmath independent of the engines (§A.7, step 2).
\end{itemize}

\textbf{A bug found this way (27 October).} When certifying π in order with n = 9, the certifier gave an optimum worse than the record of v5. The cause was a missing absolute value in the rule ``sum or difference with an anchor'' (Lemma 9): \textbar Q + δ\textbar{} was discarded when Q + δ \textless{} 0, losing about 10 \% of the classes. It was fixed, all resolution files were rebuilt, and all cells were repeated. The test t4 was added then, because the previous test used almost empty windows and could not detect it.

\textbf{A second bug (28 October), and the rule it left.} An improvement of the ln enclosure for rationals near 1 was inserted wrongly: an \texttt{else} was attached to the wrong \texttt{if}, and ordinary classes lost their logarithmic enclosure. The control recertification detected it, because π in order with n = 7 came out worse than its already certified value. It was fixed, all resolution files were rebuilt, and the computations made with that binary were invalidated. Since then every change of motor A passes a \textbf{regression test}: recertify the 52 cells already known and require the same errors to 25 digits. With the final binary it gives 52 out of 52, without differences.

\subsection*{A.9 Trusted base}\addcontentsline{toc}{subsection}{A.9 Trusted base}

The theorem assumes the following, and nothing else:

\begin{itemize}
\tightlist
\item
  \textbf{Hardware and compiler:} IEEE 754 binary64 arithmetic with round-to-nearest (SSE2 on x86-64), gcc 13 with \texttt{-O2\ -ffp-contract=off} and without \texttt{-ffast-math} (plus \texttt{-mfma} in the optimized motor A, where only the explicit fused multiply--add of CORE-MATH uses it). On Windows, MinGW with the same options; the identity of results across platforms is checked by the test runs of the n = 9 computation.
\item
  \textbf{CORE-MATH:} \texttt{cr\_\allowbreak{}\allowbreak{}pow} is correctly rounded, as its authors prove {[}8{]}. It was also checked on 28 737 cases against Arb, without any error.
\item
  \textbf{Arb} (python-flint 0.9.0, and FLINT 3.0.1 in motor C): its balls contain the exact value.
\item
  \textbf{Python 3.13:} its integers and fractions are exact. Factorizations are obtained by trial division by small primes, and a large factor is accepted as prime only below \(2^{64}\), where the BPSW test is exact; otherwise the value is treated as numeric. All factorizations are checked by multiplication.
\item
  \textbf{Constants:} \(λ_{1}\) \textless{} lo\(g_{2}\) e \textless{} \(λ_{2}\) (checked with mpmath), and δ with 30 decimals, checked against OEIS A006890 and covered by the rigorous enclosure of {[}19{]} (Section 2.2).
\item
  \textbf{Published mathematics:} the linear independence of real radicals (Besicovitch 1940 {[}9{]}, Mordell 1953 {[}17{]}, Siegel 1972 {[}18{]}).
\end{itemize}

mpmath is used only in the independent checks. Programs, tables and resolution files are identified by their hash, and every computation can be repeated with the instructions of the reproducibility package.

\subsection*{A.10 Correspondence with the code}\addcontentsline{toc}{subsection}{A.10 Correspondence with the code}

{\small
\begin{longtable}[]{@{}>{\raggedright\arraybackslash}p{(\linewidth - 2\tabcolsep) * \real{0.4147}}>{\raggedright\arraybackslash}p{(\linewidth - 2\tabcolsep) * \real{0.5853}}@{}}
\toprule\noalign{}
Result & Code \\
\midrule\noalign{}
\endhead
\bottomrule\noalign{}
\endlastfoot
Lemmas 2 and 3 (forms) & \texttt{derive} and \texttt{build\_\allowbreak{}\allowbreak{}table} (motora.c); \texttt{B.forms} (motorb.py) \\
Lemma 4 (rounding) & \texttt{pred}, \texttt{succ}; Arb \\
Lemma 5 (corners) & \texttt{imul}, \texttt{idiv}, \texttt{pow\_\allowbreak{}\allowbreak{}hull}, general case of \texttt{derive} \\
Lemma 6 (rationals) & \texttt{rat\_\allowbreak{}\allowbreak{}set}, \texttt{rat\_\allowbreak{}\allowbreak{}add} \ldots{} \texttt{rat\_\allowbreak{}\allowbreak{}pow}, \texttt{iroot}, \texttt{ipow\_\allowbreak{}\allowbreak{}lim}, \texttt{finish\_\allowbreak{}\allowbreak{}rat} \\
Lemma 7 (exp and ln) & \texttt{exp\_\allowbreak{}\allowbreak{}lo}, \texttt{exp\_\allowbreak{}\allowbreak{}hi}, \texttt{ln\_\allowbreak{}\allowbreak{}lo}, \texttt{ln\_\allowbreak{}\allowbreak{}hi}, \texttt{L2E\_\allowbreak{}\allowbreak{}LO}, \texttt{L2E\_\allowbreak{}\allowbreak{}HI} \\
Lemma 8 (small arguments) & \texttt{expm1\_\allowbreak{}\allowbreak{}hull}, \texttt{log1p\_\allowbreak{}\allowbreak{}hull} \\
Lemma 9 (derivations) & \texttt{derive}, \texttt{finish\_\allowbreak{}\allowbreak{}real}, \texttt{make\_\allowbreak{}\allowbreak{}anc}, \texttt{anc\_\allowbreak{}\allowbreak{}affine}, \texttt{finish\_\allowbreak{}\allowbreak{}from\_\allowbreak{}\allowbreak{}ln}, \texttt{decomp}, \texttt{code\_\allowbreak{}\allowbreak{}for}, \texttt{hull\_\allowbreak{}\allowbreak{}of} \\
Theorem A & \texttt{build\_\allowbreak{}\allowbreak{}tables}, \texttt{build\_\allowbreak{}\allowbreak{}table}, \texttt{merge}, \texttt{make\_\allowbreak{}\allowbreak{}literals} \\
Lemma 10 (resolution) & resolve.py (\texttt{item\_\allowbreak{}\allowbreak{}of}, \texttt{apply\_\allowbreak{}\allowbreak{}form}), \texttt{parse\_\allowbreak{}\allowbreak{}items}, build\_\allowbreak{}\allowbreak{}res.sh \\
Lemma 11 (preimages) & \texttt{windows}, \texttt{push\_\allowbreak{}\allowbreak{}affine}, \texttt{push\_\allowbreak{}\allowbreak{}anc}, \texttt{push\_\allowbreak{}\allowbreak{}log}, \texttt{mkwin}, \texttt{win\_\allowbreak{}\allowbreak{}decomp} \\
Lemma 12 (leaves) & \texttt{leaf\_\allowbreak{}\allowbreak{}lookup} \\
Theorem B & \texttt{query}, \texttt{query\_\allowbreak{}\allowbreak{}split}, \texttt{query\_\allowbreak{}\allowbreak{}splits}, \texttt{query\_\allowbreak{}\allowbreak{}top} \\
Lemma 13 (arithmetic of B) & exact.py (\texttt{X}, \texttt{sign}, \texttt{add} \ldots{} \texttt{power}, \texttt{in\_\allowbreak{}\allowbreak{}range}, \texttt{tight\_\allowbreak{}\allowbreak{}ball}, \texttt{may\_\allowbreak{}\allowbreak{}lie}), canon.py \\
Lemma 14 (queries of B) & \texttt{B.query}, \texttt{B.full}, \texttt{B.windows}, \texttt{B.family\_\allowbreak{}\allowbreak{}values} \\
Theorem C & \texttt{Verifier} (verify.py) \\
Lemma 15 (pre-filter) & motorc.c (collect\_\allowbreak{}\allowbreak{}verify, query, full, windows, op\_\allowbreak{}\allowbreak{}*); block mode and join; lab/fase2b and lab/fase2c \\
Main theorem & certify.py (\texttt{certify}, \texttt{independent\_\allowbreak{}\allowbreak{}eval}) \\
\end{longtable}
}

\subsection*{A.11 Pseudocode}\addcontentsline{toc}{subsection}{A.11 Pseudocode}

\begin{Verbatim}[fontsize=\footnotesize]
build(S):                              # |S| <= TBL, by increasing |S|
  T <- literals(S);  SING <- {}
  for each split S = S1 u S2, k1 in T(S1), k2 in T(S2), allowed form f:
    d <- derivation(k1, k2, f)          # Lemma 9
    if d is undefined or out of range: continue
    if d is singular: SING <- SING u {d};  else: T <- T u {d}
  T <- merge(T u resolution_B(SING))    # Lemma 10

query(S, W, prefix):                    # Theorem B
  if |S| <= TBL: emit prefix + k for each k in T(S) passing the filters (Lemma 12); return
  emit the literals of S passing the filters
  for each split (S1 smaller, S2), k in T(S1), allowed family f:
    if f(s, .) is constant k and k in W: emit prefix + (S1, k, f, constant)
    for each W' in P_f(k, W): query(S2, W', prefix + (S1, k, f))     # Lemma 11

certify(c, n, mode, e_rec):
  eps <- e_rec (1 + 1e-6) + 8 ulp(c);  W0 <- [c - eps, c + eps] outward
  Pi <- query(D, W0, empty)                          # motor A
  Pi <- prefilter(Pi) at 1024 and 2048 bits          # motor C, Lemma 15 (large cells)
  Vals <- union of verify(pi) for pi in Pi           # motor B, Theorem C
  v* <- the one of minimal error (Arb); group ties by exact equality
  check the formulas of v* independently
  return v*, its error, its formulas and the number of unresolved cases
\end{Verbatim}

\section*{Appendix B. The constant e: error bounds and optimality in F}\addcontentsline{toc}{section}{Appendix B. The constant e: error bounds and optimality in F}

Everything to be proved about e is closed inside the family F, without artificial search limits:

\begin{itemize}
\tightlist
\item
  \textbf{B.1:} the six records for e have a proved error bound and an exact number of correct decimals.
\item
  \textbf{B.3:} inside F, Bamberger's formula is optimal with 1--9, and ours are optimal with 1--7, 1--8, 1--8 in order, 1--9 in order and 0--9.
\item
  \textbf{B.4:} three-story identities exist with 1--9 (for example N = \({2}^{{3}^{49}}\)), but none beats Bamberger's.
\end{itemize}

The proofs use only exact arithmetic and ball arithmetic (Arb). Three independent checks agree with the main program.

\subsection*{B.0 Rules and definitions}\addcontentsline{toc}{subsection}{B.0 Rules and definitions}

\textbf{Grammar.} The grammar is that of Section 2 (any-order mode unless stated). A \emph{block} of k digits is a literal, op(u, v) with u, v expressions over a partition of the block (both orders), or −u.

\textbf{Correct decimals.} An approximation x of e is correct to k decimals if \textbar x − e\textbar{} \textless{} ½·\(10^{-k}\). This definition does not require knowing the digits of e at astronomical positions, which nobody can check.

\textbf{The family F.} F consists of the formulas (U + s·\(A^{X}\))\^{}(t·\(C^{Y}\)), where:

\begin{itemize}
\tightlist
\item
  s, t ∈ \{+1, −1\};
\item
  U, A, X, C and Y are expressions over disjoint, nonempty blocks of digits that together use every digit;
\item
  U has value exactly 1;
\item
  s·\(A^{X}\) · t·\(C^{Y}\) = 1 exactly.
\end{itemize}

The value is then \({(1+1/N)}^{N}\) with N = t·\(C^{Y}\). Every two-story formula (U + A\^{}(\(B^{X}\)))\^{}(C\^{}(\(D^{Y}\))) is in F, taking X′ = \(B^{X}\) and Y′ = \(D^{Y}\). So F contains the formulas of Sabey, Bamberger and Fastfaxr and McKernon's family, as well as our records.

We call (A, X) the \emph{L side} and (C, Y) the \emph{R side}. The identity is equivalent to \(μ_{L}\) + \(μ_{R}\) = 0, with \(μ_{L}\) = X·ln\textbar A\textbar{} and \(μ_{R}\) = Y·ln\textbar C\textbar, and ln\textbar N\textbar{} = \(μ_{R}\).

\subsection*{B.1 Rigorous error bound}\addcontentsline{toc}{subsection}{B.1 Rigorous error bound}

\textbf{Lemma B1.} For every real N ≥ 1,

\[
\frac{e}{2N}\left(1-\frac{11}{12N}\right) \;<\; e-\left(1+\frac{1}{N}\right)^{N} \;<\; \frac{e}{2N}
\]

\emph{Proof.} Let t = 1/N ∈ (0, 1{]}.

\begin{enumerate}
\def\labelenumi{\arabic{enumi}.}
\tightlist
\item
  f = N·ln(1 + t) = 1 − t/2 + t\(^{2}\)/3 − t\(^{3}\)/4 + \ldots{} is an alternating series with strictly decreasing terms. Hence 1 − t/2 \textless{} f \textless{} 1 − t/2 + t\(^{2}\)/3, that is, δ = 1 − f satisfies t/2 − t\(^{2}\)/3 \textless{} δ \textless{} t/2.
\item
  e \(-{(1+1/N)}^{N}\) = e·(1 − \(e^{-δ}\)), and for δ \textgreater{} 0, δ − δ\(^{2}\)/2 \textless{} 1 − \(e^{-δ}\) \textless{} δ.
\item
  Upper bound: e·δ \textless{} e·t/2. Lower bound: g(δ) = δ − δ\(^{2}\)/2 is increasing on (0, 1), so the error exceeds e·(t/2 − t\(^{2}\)/3 − t\(^{2}\)/8) = (e·t/2)·(1 − 11t/12). ∎
\end{enumerate}

\textbf{Lemma B1′ (negative N).} For every real N ≤ −2, with M = \textbar N\textbar,

\[
\frac{e}{2M} \;<\; \left(1+\frac{1}{N}\right)^{N}-e \;<\; \frac{e}{2M}\left(1+\frac{3}{M}\right)
\]

\emph{Proof.} Let t = 1/M ∈ (0, 1/2{]}. Then \({(1+1/N)}^{N}\) = (1 − t)\^{}(−M) = exp(1 + δ′), where δ′ is the sum of \(t^{k-1}\)/k over k ≥ 2, so the difference is e·(\(e^{δ′}\) − 1). The lower bound follows from δ′ \textgreater{} t/2 and \(e^{x}\) − 1 \textgreater{} x. For the upper bound, δ′ ≤ t/2 + t\(^{2}\)/(3(1 − t)) ≤ t/2 + 2t\(^{2}\)/3 ≤ 5/12, and \(e^{x}\) − 1 ≤ x + x\(^{2}\) for 0 ≤ x ≤ 1. Hence \(e^{δ′}\) − 1 ≤ t/2 + (11/12)t\(^{2}\) + (2/3)t\(^{3}\) + (4/9)t\(^{4}\) ≤ t/2 + 1.37t\(^{2}\) \textless{} (t/2)(1 + 3t). ∎

So a formula of F with N \textless{} 0 has an error larger than e/(2\textbar N\textbar), while one with N \textgreater{} 0 has an error smaller than e/(2N). Two formulas whose \textbar N\textbar{} differ by more than a factor 1 + 4/min\textbar N\textbar{} are therefore ordered by \textbar N\textbar, whatever their signs. In every case below the largest \textbar N\textbar{} exceeds the next one by many orders of magnitude and is attained with N \textgreater{} 0, so maximizing \textbar N\textbar{} over F minimizes the error.

\textbf{Consequence.} Let u = lo\(g_{10}\)(N/e). The formula is correct to k decimals for every k ≤ u, and fails for k ≥ u + 0.4/N. So k = ⌊u⌋, unless u lies less than 0.4/N below an integer; the script certifies that this never happens.

For each record, the script:

\begin{itemize}
\tightlist
\item
  checks that each digit is used once (and the order, for the in-order records);
\item
  verifies L·R = 1 with exact vectors of prime exponents, without floating point;
\item
  computes the error and k with Arb, at a precision sufficient for the exponent.
\end{itemize}

As a cross-check, for N = \({2}^{72}\), \({5}^{34}\) and \({5}^{81}\) the error computed directly with mpmath at 200 digits lies inside the bounds.

\textbf{Table B1.} The six records for e.

{\footnotesize
\begin{longtable}[]{@{}>{\raggedright\arraybackslash}p{(\linewidth - 8\tabcolsep) * \real{0.1950}}>{\raggedright\arraybackslash}p{(\linewidth - 8\tabcolsep) * \real{0.1950}}>{\raggedright\arraybackslash}p{(\linewidth - 8\tabcolsep) * \real{0.1715}}>{\raggedright\arraybackslash}p{(\linewidth - 8\tabcolsep) * \real{0.1715}}>{\raggedright\arraybackslash}p{(\linewidth - 8\tabcolsep) * \real{0.2670}}@{}}
\toprule\noalign{}
Digits & Formula & N & −lo\(g_{10}\)(error) & Correct decimals k \\
\midrule\noalign{}
\endhead
\bottomrule\noalign{}
\endlastfoot
1--7, any order & \({(1+{4}^{-36})}^{{.5}^{-72}}\) & \({2}^{72}\) & 21.5408952016 & 21 \\
1--8, any order & \({(1+{2}^{-{8}^{7}})}^{{4}^{{.5}^{-6/.3}}}\) & \({4}^{{2}^{20}}\) = \({2}^{{2}^{21}}\) & 631 305.524202 & 631 305 \\
1--9, any order (Bamberger) & \({\left(1+{.2}^{{9}^{6\cdot 7}}\right)}^{{5}^{{3}^{84}}}\) & \({5}^{{3}^{84}}\) & 8.36842898907·\(10^{39}\) & 8 368 428 989 068 425 943 817 590 916 445 001 887 164 \\
1--8, in order & \({(1+{.2}^{34})}^{{5}^{6\cdot 7-8}}\) & \({5}^{34}\) & 23.6317156612 & 23 \\
1--9, in order & \({(1+{.2}^{{3}^{4}})}^{{5}^{6\cdot (7+8)\cdot .9}}\) & \({5}^{81}\) & 56.4833058650 & 56 \\
0--9, any order & \({(1+{.2}^{{9}^{70}})}^{{5}^{{3}^{84/.6}}}\) & \({5}^{{3}^{140}}\) & 4.37959750359·\(10^{66}\) & 4 379 597 503 586 508 305 506 762 455 287 865 626 188 758 462 213 921 632 491 859 893 686 \\
\end{longtable}
}

The two bounds agree in the digits shown, since their relative difference is below 1/N.

\subsection*{B.2 The size lemma}\addcontentsline{toc}{subsection}{B.2 The size lemma}

\textbf{Lemma B2 (size).} As \textbar L\textbar·\textbar R\textbar{} = 1, ln\textbar N\textbar{} = \textbar{}\(μ_{R}\)\textbar{} = \textbar{}\(μ_{L}\)\textbar. Let SB(k) be the maximum of \textbar X·ln\textbar A\textbar\textbar{} over all sides of k digits. Then ln\textbar N\textbar{} ≤ SB(k) for every formula whose shorter side has k digits.

SB(2) and SB(3) are computed by enumerating all blocks of 1 and 2 digits, with rigorous upper bounds (Arb). With 1--9: SB(2) ≤ \(e^{3}\).032 ≈ 21 and SB(3) ≤ \(e^{21}\).456 ≈ 2.1·\(10^{9}\).

\subsection*{B.3 Optimality in F}\addcontentsline{toc}{subsection}{B.3 Optimality in F}

\textbf{Theorem B3.} With the digits 1--n in any order, the maximum of \textbar N\textbar{} in F is:

{\small
\begin{longtable}[]{@{}>{\raggedright\arraybackslash}p{(\linewidth - 6\tabcolsep) * \real{0.1663}}>{\raggedright\arraybackslash}p{(\linewidth - 6\tabcolsep) * \real{0.3399}}>{\raggedright\arraybackslash}p{(\linewidth - 6\tabcolsep) * \real{0.3275}}>{\raggedright\arraybackslash}p{(\linewidth - 6\tabcolsep) * \real{0.1663}}@{}}
\toprule\noalign{}
Digits & Maximal N & ln(ln N) & Distinct writings by block values \\
\midrule\noalign{}
\endhead
\bottomrule\noalign{}
\endlastfoot
1--7 & \({2}^{72}\) & 3.9101531984343909840 & 8 \\
1--8 & \({2}^{{2}^{21}}\) & 14.189577871177187171 & 16 \\
1--9 & \({5}^{{3}^{84}}\) (Bamberger) & 92.759317243448324698 & 16 \\
\end{longtable}
}

In the three cases the maximum is unique as a value of N. The writings counted are the choices of the blocks (A, X, C, Y) together with the signs s and t that make s·\(A^{X}\)·t·\(C^{Y}\) = 1; with odd exponents half of them need s = −1 or t = −1, and read without those signs, as (U + \(A^{X}\))\^{}(\(C^{Y}\)), they would give \({(1-1/N)}^{N}\) ≈ 1/e. The distinct writings differ only in the signs or in the roles of the blocks, for example \({(1+{.2}^{{3}^{84}})}^{{5}^{{9}^{42}}}\) or \({(1+{5}^{-{3}^{84}})}^{{.2}^{-{9}^{42}}}\).

\emph{Proof.}

\begin{enumerate}
\def\labelenumi{\arabic{enumi}.}
\tightlist
\item
  \textbf{Reduction to 4 + 4 (1--9).} Bamberger has ln N = \({3}^{84}\)·ln 5 ≈ 1.9·\(10^{40}\) \textgreater{} SB(3). So a formula that matches or beats it has both sides of at least 4 digits. With 9 digits this forces U to have a single digit (U = 1) and sides of exactly 4 + 4, so every block has at most 3 digits. The same holds, with its own threshold, for 1--8 (sides 3 + 4 or 3 + 3) and 1--7 (sides 3 + 3).
\item
  \textbf{Exact enumeration.} The complete set of values of every block of 1, 2 and 3 digits is built, with no size or denominator limits. Values of blocks of up to 2 digits are always signed radicals s·∏ \(p^{e_{p}}\) with rational e\_\allowbreak{}\allowbreak{}p, represented exactly. With 3 digits, sums radical + rational and powers with irrational exponent also appear; those are enclosed in Arb balls.
\item
  \textbf{Rigorous sieve.} For every split U \textbar{} L side \textbar{} R side and every pair of blocks, ln\textbar μ\textbar{} is enclosed in an Arb ball (256 bits). An identity requires \(μ_{L}\) = −\(μ_{R}\); if the balls do not overlap, the identity is certainly false. Before that test, candidate pairs are preselected by the floating-point midpoints of the balls, with a relative tolerance of \(10^{-9}\); this loses no pair because the relative radii of the balls are far smaller (at most 3·\(10^{-75}\) over the 451 122 entries of the case 1--7, measured).
\item
  \textbf{Exact verification.} Every overlap is checked with prime exponents: X·\(α_{p}\) + Y·\(γ_{p}\) = 0 for every prime p, where ln\textbar A\textbar{} = Σ \(α_{p}\) ln p and ln\textbar C\textbar{} = Σ \(γ_{p}\) ln p.~That condition suffices for the identity to hold. In the three runs every overlap passed the verification: none was left unresolved, and no transcendence theorem was needed. ∎
\end{enumerate}

\textbf{Independent checks.}

\begin{itemize}
\tightlist
\item
  \textbf{Blocks.} A second enumerator, written separately with sympy, reproduces exactly the value sets of 7 test blocks of 1, 2 and 3 digits, among them \{2,5,9\} with 2296 values and \{3,4,8\} with 2658, with 0 differences.
\item
  \textbf{Old program.} The search of 26 October, with another algorithm and with limits, reaches the same three maxima.
\item
  \textbf{Independent brute force (1--7).} A third program uses the sympy blocks, mpmath at 80 digits and its own matching, reusing no code of the search. It finds the same 8 writings and the same unique value N = \({2}^{72}\).
\end{itemize}

\subsection*{B.4 Towers of three or more stories}\addcontentsline{toc}{subsection}{B.4 Towers of three or more stories}

Three-story identities do exist with 1--9. For example

\[
\left(1+.5^{3^{7^{2}}}\right)^{8^{(9^{6})^{4}}},\qquad .5^{3^{49}} = 2^{-3^{49}},\quad 8^{9^{24}} = 2^{3\cdot 3^{48}} = 2^{3^{49}}
\]

With N = \({2}^{{3}^{49}}\), this formula is correct to about 7.2·\(10^{22}\) decimals: a valid identity, but far worse than Bamberger's (about 8.4·\(10^{39}\)).

\textbf{Proposition B4.} In F with 1--9:

\begin{enumerate}
\def\labelenumi{\arabic{enumi}.}
\tightlist
\item
  No formula with three or more stories matches or beats Bamberger's. This is a direct corollary of Theorem B3.
\item
  The best three-story identity is N = \({2}^{{3}^{49}}\). The three largest values of N in F (5\^{}(3\^{}84) of Bamberger, \({2}^{{3}^{84}}\), and \({3}^{{2}^{85}}\) of Sabey) were checked one by one: none admits a writing in which X or Y is a tower \(d_{1}\)\^{}(\(d_{2}\)\^{}\(d_{3}\)).
\item
  With four or more stories, one side takes at least 5 digits, so the other has at most 3. By the size lemma, ln N ≤ 2.1·\(10^{9}\).
\end{enumerate}

\textbf{Table B2.} The largest identities of F with 1--9 (complete enumeration of those with ln N ≥ \(e^{21}\).5).

{\footnotesize
\begin{longtable}[]{@{}>{\raggedright\arraybackslash}p{(\linewidth - 8\tabcolsep) * \real{0.1215}}>{\raggedright\arraybackslash}p{(\linewidth - 8\tabcolsep) * \real{0.2752}}>{\raggedright\arraybackslash}p{(\linewidth - 8\tabcolsep) * \real{0.1215}}>{\raggedright\arraybackslash}p{(\linewidth - 8\tabcolsep) * \real{0.1215}}>{\raggedright\arraybackslash}p{(\linewidth - 8\tabcolsep) * \real{0.3602}}@{}}
\toprule\noalign{}
Rank & N & ln(ln N) & Stories & Example \\
\midrule\noalign{}
\endhead
\bottomrule\noalign{}
\endlastfoot
1 & \({5}^{{3}^{84}}\) & 92.7593 & 2 & \({(1+{.2}^{{3}^{84}})}^{{5}^{{9}^{6\cdot 7}}}\) (Bamberger) \\
2 & \({2}^{{3}^{84}}\) & 91.9169 & 2 & \({(1+{2}^{-{3}^{84}})}^{{.5}^{-{9}^{6\cdot 7}}}\) \\
3 & \({3}^{{2}^{85}}\) & 59.0116 & 2 & \({(1+{3}^{-{2}^{85}})}^{{9}^{{4}^{6\cdot 7}}}\) (Sabey) \\
4 & \({2}^{{3}^{49}}\) & 53.4655 & 3 & \({\left(1+{.5}^{{3}^{{7}^{2}}}\right)}^{{8}^{{({9}^{6})}^{4}}}\) \\
5 & \({4}^{{.5}^{-69}}\) = \({2}^{{2}^{70}}\) & 48.1538 & 2 & \({\left(1+{2}^{-{8}^{7/.3}}\right)}^{{4}^{{.5}^{-69}}}\) \\
\end{longtable}
}

There are 11 distinct values of N above that threshold in total; the complete list is in the script's output.

\subsection*{B.5 e in order, and the variant 0--9}\addcontentsline{toc}{subsection}{B.5 e in order, and the variant 0--9}

\textbf{Theorem B5 (in order).} With 1--8 and 1--9 in order, the maximum of \textbar N\textbar{} in \(F_{ord}\) is \({5}^{34}\) and \({5}^{81}\) respectively: our formulas \({(1+{.2}^{34})}^{{5}^{6\cdot 7-8}}\) and \({(1+{.2}^{{3}^{4}})}^{{5}^{6\cdot (7+8)\cdot .9}}\). The maximum is unique as a value of N, with 4 and 28 writings. With 1--7 in order there is no identity of \(F_{ord}\) with ln N ≥ 7.4, so there the best value does not come from this family (it is certified by Theorem 1: Table 2).

\(F_{ord}\) is F with the digits in written order: the blocks are consecutive stretches in the order (U, A, X, C, Y) or (A, X, U, C, Y), and the operations do not reorder their operands. The proof is that of Theorem B3, with these differences:

\begin{itemize}
\tightlist
\item
  the size bounds are computed for each concrete configuration;
\item
  the n = 9 record has a side of 5 digits (Y = 6·(7+8)·.9), so all in-order blocks of up to 4 digits are enumerated completely;
\item
  0 configurations needed larger blocks, and 0 overlaps were left unresolved;
\item
  the blocks were checked: on the stretches 123, 234, 345 and 456 they agree exactly with a separate in-order enumerator in sympy, and on the 7 stretches of 3 digits they are a subset of the general case (0 discrepancies).
\end{itemize}

\textbf{Theorem B6 (0--9).} With the digits 0--9 in any order, the maximum of \textbar N\textbar{} in F is \({5}^{{3}^{140}}\), that is, our formula \({(1+{.2}^{{9}^{70}})}^{{5}^{{3}^{84/.6}}}\). It is unique as a value of N, with 16 writings.

\textbf{Zeros.} Friedman fixes no convention: his 2004 problem uses the digits 1--n, and of the variant with 0 he only says that Joseph DeVincentis sent him results, without rules or table. The search therefore uses the most permissive grammar: any permutation is a literal, leading zeros included (``05'', ``0.5'', ``.05''). A maximum under this grammar is also one under any stricter convention in which our formula is still valid, and ours only uses 70 and .6. So the result does not depend on that choice.

\emph{Proof (0--9).}

\begin{enumerate}
\def\labelenumi{\arabic{enumi}.}
\tightlist
\item
  \textbf{Configurations.} SB(3) ≤ \(e^{21}\).46, while ln ln N of the record is 154.28. By the size lemma, either U has 1 digit and the sides are 4 + 5, or U has 2 digits (for example \({2}^{0}\)) and the sides are 4 + 4: 2212 configurations in total.
\item
  \textbf{Side of 5 digits.} It needs blocks of 4 digits, but only those that can reach the threshold. As X with a one-digit A, they need ln\textbar X\textbar{} ≥ θ − ln ln 10 − 1; as A with a one-digit X, ln\textbar ln\textbar A\textbar\textbar{} ≥ θ − ln 9 − 1. They are generated with upper bounds of \textbar ln\textbar{} per operation and margin 1. These are pre-estimates, and every value that passes is computed exactly or with Arb balls.
\item
  \textbf{Check of the filter.} On 3 test blocks of 4 digits it was compared with the complete enumeration, with low thresholds so that thousands of values enter: none is missing.
\item
  \textbf{Astronomical values.} Numbers such as \({2}^{{3}^{{4}^{5}}}\) are not expanded but enclosed in the logarithmic domain, and sums with an astronomical term are bounded with log1p; everything stays rigorous. A side of 4 digits has ln\textbar μ\textbar{} ≤ 2.4·\(10^{9}\), so entries with larger ln\textbar μ\textbar{} cannot be paired.
\item
  \textbf{Result.} A single identity above the threshold, ours, verified exactly; 0 unresolved overlaps, in 6 minutes of computation. ∎
\end{enumerate}

\subsection*{B.6 Scope and open questions}\addcontentsline{toc}{subsection}{B.6 Scope and open questions}

The theorems hold inside F (or \(F_{ord}\)) and with the grammar above. Open:

\begin{itemize}
\tightlist
\item
  \textbf{Outside F:} other forms of the limit, such as \({(1+u/N)}^{N/u}\), or N written as a sum or product of towers. This is a research problem, which is why the text says ``optimal within F''.
\item
  \textbf{e with n ≤ 6, and 1--7 in order:} F produces no useful identity there. Those cells come from the general numerical search and are certified by Theorem 1. Theorem 2 (Appendix C) extends them to expressions with intermediate values out of range.
\end{itemize}

\subsection*{B.7 Reproducibility}\addcontentsline{toc}{subsection}{B.7 Reproducibility}

Everything is in the folder \texttt{demostraciones\_\allowbreak{}\allowbreak{}e/\allowbreak{}} of the package, with the output of each script. Run with Python 3.13.17, python-flint 0.9.0 (Arb), sympy 1.14.0 and mpmath 1.3.0.

{\small
\begin{longtable}[]{@{}>{\raggedright\arraybackslash}p{(\linewidth - 6\tabcolsep) * \real{0.2694}}>{\raggedright\arraybackslash}p{(\linewidth - 6\tabcolsep) * \real{0.3924}}>{\raggedright\arraybackslash}p{(\linewidth - 6\tabcolsep) * \real{0.1885}}>{\raggedright\arraybackslash}p{(\linewidth - 6\tabcolsep) * \real{0.1498}}@{}}
\toprule\noalign{}
Script & What it does & Time & SHA-256 (start) \\
\midrule\noalign{}
\endhead
\bottomrule\noalign{}
\endlastfoot
e\_\allowbreak{}\allowbreak{}records\_\allowbreak{}\allowbreak{}bounds.py & B.1: exact identities, error bounds and k of the 6 records & \textless{} 1 s & 8142915d434e \\
e\_\allowbreak{}\allowbreak{}F.py & B.3 and B.4 in any order: blocks, size lemma, Arb sieve and exact verification; the core used by the others & 3--13 s per case & d707c0255fea \\
e\_\allowbreak{}\allowbreak{}\(F_{ord}\).py & \(F_{ord}\), in order, with blocks of up to 4 digits & 6 s & b8c4abb43239 \\
e\_\allowbreak{}\allowbreak{}\(F_{09}\).py & 0--9, with filtered blocks of 4 digits & 6 min & fd6f0830e2e1 \\
check\_\allowbreak{}\allowbreak{}blocks\_\allowbreak{}\allowbreak{}sympy.py & Second block enumerator (sympy), any order & 2 min & c5b649dabc7e \\
check\_\allowbreak{}\allowbreak{}blocks\_\allowbreak{}\allowbreak{}ord.py & The same in order, plus inclusion in the general case & 1 min & 916fbc3484ed \\
check\_\allowbreak{}\allowbreak{}F7\_\allowbreak{}\allowbreak{}independent.py & Independent brute force for 1--7 & ≈ 5 min & 5f5ad69dc440 \\
check\_\allowbreak{}\allowbreak{}V4big.py & The filter of 4-digit blocks loses nothing & 2 min & 011984ff5025 \\
check\_\allowbreak{}\allowbreak{}3story.py & Three-story writings of the largest identities & 15 s & 264a4fe6ef96 \\
\end{longtable}
}

The numerical argument of each command (\texttt{python3\ e\_\allowbreak{}\allowbreak{}F.py\ 123456789\ 92.0}, \ldots) is the threshold for ln(ln N), which must be below the maximum. Every output must say ``unresolved overlaps: 0'' and, where it appears, ``configurations needing blocks \textgreater{} \ldots: 0''; otherwise the theorem is not proved for that case. All saved outputs satisfy both conditions.

\section*{Appendix C. Expressions with intermediate values out of range}\addcontentsline{toc}{section}{Appendix C. Expressions with intermediate values out of range}

\subsection*{C.1 Statement}\addcontentsline{toc}{subsection}{C.1 Statement}

Theorem 1 (Section 5) certifies the optimum among the set \(R_{n}\) of expressions all of whose intermediate values are 0 or have absolute value in {[}\(10^{-300}\), \(10^{300}\){]}. Here we remove that restriction.

\textbf{Theorem 2.} For each cell of Theorem 1, no valid expression over 1, \ldots, n (Section 2, no range restriction) has a value whose distance to c is smaller than the optimum E* of Theorem 1. Hence the optimum of Theorem 1 is the optimum over all valid expressions. Since the query window W is closed and contains {[}c − E*, c + E*{]} (its radius is at least E* enlarged by \(10^{-9}\) in relative terms; Section 6), the same computation shows that no expression outside \(R_{n}\) attains the error E* either: every optimal expression has the optimal value of Tables 2--3, which is unique.

The proof is a second certified computation (``motor X''). It is independent of the first in its arithmetic (enclosures of logarithms instead of enclosures of values), but it reuses motor A's tables of in-range classes and motor B's exact arithmetic.

\subsection*{C.2 Notation}\addcontentsline{toc}{subsection}{C.2 Notation}

All values are taken as magnitudes, by Lemma 1 (V(S) = −V(S)). Powers have positive bases by Lemma 2, and the forms are those of Lemma 3. For a set of digits S:

\begin{itemize}
\tightlist
\item
  U(S) is the set of magnitudes of all valid expressions over S;
\item
  V(S) ⊆ U(S) is the set of those with every intermediate value 0 or in {[}\(10^{-300}\), \(10^{300}\){]} (the set of Theorem 1);
\item
  NE(S) = U(S) ∖ V(S) is the set of \emph{new or extreme} values;
\item
  X(S) ⊆ U(S) is the set of magnitudes of the expressions that have at least one \emph{node} (the expression itself or a subexpression) whose value is out of range.
\end{itemize}

A value u ∉ {[}\(10^{-300}\), \(10^{300}\){]} ∪ \{0\} is \emph{out of range}. X(S) is defined by the structure of the expressions, not by their values: it also contains values that some other, in-range, expression attains, such as \({9}^{{4}^{8}}/{3}^{{2}^{17}}\) = 1. Four facts are immediate:

\begin{enumerate}
\def\labelenumi{\arabic{enumi}.}
\tightlist
\item
  NE(S) ⊆ X(S): a value outside V(S) has no writing with all nodes in range.
\item
  U(S) = V(S) ∪ X(S): every expression has all its nodes in range or not.
\item
  A literal is never out of range (it lies between 0.1 and \(10^{9}\)), so X(\{d\}) = ∅.
\item
  For the variant without unary minus (Section 8), \(V_{11}\)(S), \(X_{11}\)(S) and \(U_{11}\)(S) are defined in the same way over the expressions without unary minus. Then \(V_{11}\)(S) ⊆ V(S), \(X_{11}\)(S) ⊆ X(S) and \(U_{11}\)(S) = \(V_{11}\)(S) ∪ \(X_{11}\)(S).
\end{enumerate}

Theorem 2 follows from X(\{1,\ldots,n\}) ∩ W = ∅, which is what the computation proves.

\subsection*{C.3 Generation lemma}\addcontentsline{toc}{subsection}{C.3 Generation lemma}

\textbf{Lemma X1 (structural generation).} Let E = (\(E_{1}\) op \(E_{2}\)) be an expression over the split S = \(S_{1}\) ⊔ \(S_{2}\) (in the in-order mode, \(S_{1}\) the left part), reduced as in Lemma 2, which keeps the magnitude of every node. Let s = \textbar val \(E_{1}\)\textbar, l = \textbar val \(E_{2}\)\textbar{} and \textbar val E\textbar{} = f(s, l) with f a form of Lemma 3. If E has a node out of range, then at least one of the following holds:

\begin{enumerate}
\def\labelenumi{(\alph{enumi})}
\tightlist
\item
  \(E_{1}\) has a node out of range, so s ∈ X(\(S_{1}\)) and l ∈ U(\(S_{2}\)) = V(\(S_{2}\)) ∪ X(\(S_{2}\));
\item
  \(E_{2}\) has a node out of range, so l ∈ X(\(S_{2}\)) and s ∈ U(\(S_{1}\));
\item
  all nodes of \(E_{1}\) and \(E_{2}\) are 0 or in range, so s ∈ V(\(S_{1}\)) and l ∈ V(\(S_{2}\)), and f(s, l) is out of range.
\end{enumerate}

\emph{Proof.} The out-of-range node is in \(E_{1}\), in \(E_{2}\) or is the root. A unary minus at the root does not change any magnitude (Lemma 1), and E is not a literal (fact 3). ∎

The lemma is about expressions, not values: no writing is replaced by another one. An earlier version was stated for values (u ∈ NE(S)) and replaced each child by an in-range writing of the same value. That is enough for Theorem 2, but not for the variant without unary minus, where the in-range writing may need the unary minus. The structural form covers both.

\subsection*{C.4 Enclosures in the logarithmic domain}\addcontentsline{toc}{subsection}{C.4 Enclosures in the logarithmic domain}

An \emph{enclosure} of a magnitude u has four parts:

\begin{enumerate}
\def\labelenumi{\arabic{enumi}.}
\tightlist
\item
  \textbf{Interval and flags.} An interval {[}A, B{]} of doubles with A ≤ ln u ≤ B (with \(e^{-∞}\) = 0, so A = −∞ also allows u = 0), and flags Z (u = 0 exactly) and O (u = 1 exactly).
\item
  \textbf{Second level.} S = sign(ln u) ∈ \{−1, 0, +1\}, and, when S ≠ 0, M ≤ ln\textbar ln u\textbar{} ≤ N. This holds values very close to 1, such as 1 + \(10^{-700}\), whose logarithm a double cannot hold. It also holds values whose logarithm exceeds the double range, such as \({3}^{{4}^{7776}}\), with ln ln u = 7776·ln 4 + ln ln 3.
\item
  \textbf{Anchor.} u = (p/q)(1 + δ), with p/q an exact rational, and either δ = 0 (the value is the exact rational) or sign δ known and ln\textbar δ\textbar{} ≤ −20 enclosed. This holds shadows of rationals such as 3 + \(10^{-700}\). Values near 1 with a second level act as anchors at 1.
\item
  \textbf{Third level.} For u = 1 ± t with t below \(e^{-10^{308}}\) and ln t = −\(e^{μ}\) known: ln\textbar ln u\textbar{} = −\(e^{ν}\) with ν enclosed. This holds, for example, \(1+{2}^{-{3}^{{4}^{5}}}\).
\end{enumerate}

Without parts 2--4, expressions such as \({(1+\varepsilon)}^{N}\), the mechanism behind every record for e, could not be evaluated: for them the interval alone gives (0, ∞).

The rest of this section proves Lemmas X2 and X3 rule by rule. Function names refer to \texttt{ivl.py}.

\subsubsection*{C.4.1 Enclosures and what ``outward'' means}\addcontentsline{toc}{subsubsection}{C.4.1 Enclosures and what ``outward'' means}

An enclosure of a magnitude u ≥ 0 is a record of fields, each of which is a statement about u:

{\small
\begin{longtable}[]{@{}>{\raggedright\arraybackslash}p{(\linewidth - 2\tabcolsep) * \real{0.3739}}>{\raggedright\arraybackslash}p{(\linewidth - 2\tabcolsep) * \real{0.6261}}@{}}
\toprule\noalign{}
field & statement \\
\midrule\noalign{}
\endhead
\bottomrule\noalign{}
\endlastfoot
Z & u = 0 \\
O & u = 1 \\
A, B (doubles or ±∞) & A ≤ ln u ≤ B, with ln 0 = −∞ \\
S ∈ \{−1, 0, +1\}; M, N & if S ≠ 0: S·ln u ≥ 0 and M ≤ ln\textbar ln u\textbar{} ≤ N, with ln 0 = −∞ (so ln u ≠ 0, and the sign of ln u is S, only when M \textgreater{} −∞) \\
AP, AQ, DS, DM, DN & if AQ \textgreater{} 0: u = (AP/AQ)(1 + δ) with δ = 0 when DS = 0, and otherwise DS·δ ≥ 0 and DM ≤ ln\textbar δ\textbar{} ≤ DN (so δ ≠ 0 only when DM \textgreater{} −∞) \\
T3, M3, N3 & if T3 = 1: ln\textbar ln u\textbar{} = −\(e^{ν}\) with M3 ≤ ν ≤ N3 \\
\end{longtable}
}

A value v \emph{belongs} to an enclosure when every statement holds for v. Fields that carry no information are neutral (A = −∞, B = +∞, S = 0, AQ = 0, T3 = 0). The rational AP/AQ is always reduced, with AQ \textgreater{} 0 and AP ≥ 0: motor A's rationals are normalized by \texttt{rat\_\allowbreak{}\allowbreak{}set} and made nonnegative by \texttt{finish\_\allowbreak{}\allowbreak{}rat}, the literals by a gcd, and every rational built in \texttt{ivl.py} by \texttt{rat\_\allowbreak{}\allowbreak{}op}/\texttt{ipow} from reduced ones.

The proof of each rule has two parts: an identity or inequality of real numbers, and the remark that the computed double bounds are obtained from bounds of the inputs by expressions that are monotone in each input and evaluated with outward rounding. The primitives are:

\begin{itemize}
\tightlist
\item
  \textbf{(P1) Rounded operations.} For a sum, difference, product or quotient r of two doubles, \texttt{pred(\allowbreak{}fl(\allowbreak{}r))\ ≤\ r\ ≤\ succ(\allowbreak{}fl(\allowbreak{}r))} (Lemma 4). \texttt{addlo}/\texttt{addhi} (and \texttt{divlo}/\texttt{divhi}) return these bounds; a NaN (∞ − ∞, 0·∞) becomes −∞ in a lower bound and +∞ in an upper bound, which is always valid. Multiplying a double by 2 and negating it are exact.
\item
  \textbf{(P2) Exponential.} \texttt{exp\_\allowbreak{}\allowbreak{}lo(\allowbreak{}z)\ ≤\ e\^{}z\ ≤\ exp\_\allowbreak{}\allowbreak{}hi(\allowbreak{}z)} for every double z (Lemma 7, from the correctly rounded CORE-MATH \texttt{pow}), with exp(−∞) = 0. On overflow \texttt{exp\_\allowbreak{}\allowbreak{}hi} returns +∞ and \texttt{exp\_\allowbreak{}\allowbreak{}lo} returns at most DBL\_\allowbreak{}\allowbreak{}MAX; on underflow \texttt{exp\_\allowbreak{}\allowbreak{}lo} returns 0 and \texttt{exp\_\allowbreak{}\allowbreak{}hi} a positive subnormal. All of these are valid bounds.
\item
  \textbf{(P3) Logarithm.} \texttt{ln\_\allowbreak{}\allowbreak{}lo(\allowbreak{}x)\ ≤\ ln\ x\ ≤\ ln\_\allowbreak{}\allowbreak{}hi(\allowbreak{}x)} for x \textgreater{} 0 (Lemma 7: a deterministic guess accepted only after it is checked with the exp bounds), ln 0 = −∞, ln +∞ = +∞. Near 1 (\textbar x − 1\textbar{} \textless{} \(2^{-9}\)), u = x − 1 is exact (Sterbenz) and u − u\(^{2}\) ≤ ln x ≤ u, which tightens both bounds (Lemma 8).
\item
  \textbf{(P4) Small arguments.} For t ≥ 0: t − t\(^{2}\) ≤ ln(1 + t) ≤ t; for 0 ≤ t ≤ 1/2: −t − t\(^{2}\) ≤ ln(1 − t) ≤ −t. \texttt{log1p\_\allowbreak{}\allowbreak{}lo/\allowbreak{}hi} and \texttt{log1m\_\allowbreak{}\allowbreak{}lo/\allowbreak{}hi} use these for t \textless{} \(2^{-20}\) and P3 on 1 ± t (rounded outward) otherwise. Each of them receives a bound of t on the side that makes the result a bound (for example \texttt{log1m\_\allowbreak{}\allowbreak{}lo} receives an upper bound of t), since ln(1 ± t) is monotone in t.
\item
  \textbf{(P5) Log-sum and log-difference.} \texttt{lse\_\allowbreak{}\allowbreak{}lo/\allowbreak{}hi(\allowbreak{}x,\ y)} bound ln(\(e^{x}\) + \(e^{y}\)) = m + ln(1 + \(e^{min - m}\)), m = max(x, y), which is increasing in x and y. \texttt{ldiff\_\allowbreak{}\allowbreak{}lo(\allowbreak{}x,\ y)} and \texttt{ldiff\_\allowbreak{}\allowbreak{}hi(\allowbreak{}x,\ y)} bound ln(\(e^{a}\) − \(e^{b}\)) = a + ln(1 − \(e^{b - a}\)) for a \textgreater{} b: the lower bound is taken from a lower bound of a and an upper bound of b, and the upper bound the other way round. The result is −∞ when a \textgreater{} b is not certain.
\end{itemize}

Because each rule only asserts true statements about every value of the inputs' enclosures, the record returned by \texttt{fwd} (whose fields are set by several rules, each valid on its own, and intersected by \texttt{tighten}) is an enclosure of every value f(\(v_{1}\), \(v_{2}\)).

\textbf{Lemma X2 (forward enclosures).} Let X, Y be enclosures and f one of the forms 0 s+l, 1 \textbar s−l\textbar, 2 s·l, 3 s/l, 4 l/s, 5 \(s^{l}\), 6 s\^{}−l, 7 \(l^{s}\), 8 l\^{}−s. \texttt{fwd(\allowbreak{}f,\ X,\ Y)} returns an enclosure R and a flag ok such that f(\(v_{1}\), \(v_{2}\)) belongs to R for all \(v_{1}\) ∈ X, \(v_{2}\) ∈ Y with f(\(v_{1}\), \(v_{2}\)) defined; ok is false only where f(\(v_{1}\), \(v_{2}\)) is undefined for all such pairs.

\emph{Proof.} Form by form, below; then the anchors (C.4.6) and the final step (C.4.7). ∎

\subsubsection*{C.4.2 Products and quotients (forms 2, 3, 4)}\addcontentsline{toc}{subsubsection}{C.4.2 Products and quotients (forms 2, 3, 4)}

Form 4 is form 3 with the operands exchanged, l/s, so only s·l and s/l are treated. Write \(λ_{i}\) = ln \(v_{i}\) and μᵢ = ln\textbar{}\(λ_{i}\)\textbar.

\begin{enumerate}
\def\labelenumi{\arabic{enumi}.}
\tightlist
\item
  \textbf{Value.} ln(\(v_{1}\)\(v_{2}\)) = \(λ_{1}\) + \(λ_{2}\) and ln(\(v_{1}\)/\(v_{2}\)) = \(λ_{1}\) − \(λ_{2}\). Hence A = addlo(\(A_{1}\), \(A_{2}\)), B = addhi(\(B_{1}\), \(B_{2}\)) for the product, and A = addlo(\(A_{1}\), −\(B_{2}\)), B = addhi(\(B_{1}\), −\(A_{2}\)) for the quotient (P1).
\item
  \textbf{Exact cases.} \(v_{1}\)\(v_{2}\) = 0 iff a factor is 0 (Z = \(Z_{1}\) ∨ \(Z_{2}\)). \(v_{1}\)/\(v_{2}\) is undefined iff \(v_{2}\) = 0 (ok = ¬\(Z_{2}\)), and it is 0 iff \(v_{1}\) = 0 (Z = \(Z_{1}\) ∧ ¬\(Z_{2}\)). 1·1 = 1/1 = 1 (O = \(O_{1}\) ∧ \(O_{2}\)).
\item
  \textbf{Second level from the value.} If A \textgreater{} 0 then ln u \textgreater{} 0 and ln ln u ∈ {[}ln A, ln B{]}; if B \textless{} 0 then ln\textbar ln u\textbar{} ∈ {[}ln(−B), ln(−A){]} (\texttt{mu\_\allowbreak{}\allowbreak{}default}, P3).
\item
  \textbf{A factor exactly 1.} If \(v_{1}\) = 1 the result is \(v_{2}\) or 1/\(v_{2}\): ln u = ±\(λ_{2}\), so S = ±\(S_{2}\) and (M, N) = (\(M_{2}\), \(N_{2}\)). Symmetrically if \(v_{2}\) = 1.
\item
  \textbf{Two values near 1, contributions of the same sign.} If sign \(λ_{1}\) = sign(±\(λ_{2}\)) (the sign that \(λ_{2}\) has in ln u), then \textbar ln u\textbar{} = \textbar{}\(λ_{1}\)\textbar{} + \textbar{}\(λ_{2}\)\textbar{} = e\^{}μ₁ + e\^{}μ₂, so ln\textbar ln u\textbar{} ∈ {[}lse\_\allowbreak{}\allowbreak{}lo(\(M_{1}\), \(M_{2}\)), lse\_\allowbreak{}\allowbreak{}hi(\(N_{1}\), \(N_{2}\)){]} (P5).
\item
  \textbf{Opposite contributions, one dominant.} If the signs are opposite and \(N_{2}\) \textless{} \(M_{1}\), then \textbar{}\(λ_{2}\)\textbar{} ≤ \(e^{N_{2}}\) \textless{} \(e^{M_{1}}\) ≤ \textbar{}\(λ_{1}\)\textbar. So ln u has the sign of \(λ_{1}\) and \textbar ln u\textbar{} = \textbar{}\(λ_{1}\)\textbar(1 − \textbar{}\(λ_{2}\)\textbar/\textbar{}\(λ_{1}\)\textbar), that is, ln\textbar ln u\textbar{} = μ₁ + ln(1 − \(e^{μ_{2} - μ_{1}}\)) with μ₂ − μ₁ ≤ \(N_{2}\) − \(M_{1}\) \textless{} 0. This gives M = addlo(\(M_{1}\), log1m\_\allowbreak{}\allowbreak{}lo(exp\_\allowbreak{}\allowbreak{}hi(addhi(\(N_{2}\), −\(M_{1}\))))) and N = \(N_{1}\) (P2, P4); when \(e^{N_{2} - M_{1}}\) rounds up to 1, log1m\_\allowbreak{}\allowbreak{}lo gives −∞ and only the sign and N remain. Symmetrically when \(λ_{2}\) dominates.
\end{enumerate}

\subsubsection*{C.4.3 Sums (form 0)}\addcontentsline{toc}{subsubsection}{C.4.3 Sums (form 0)}

\begin{enumerate}
\def\labelenumi{\arabic{enumi}.}
\tightlist
\item
  \textbf{Value.} ln(\(v_{1}\) + \(v_{2}\)) = ln(e\^{}\(λ_{1}\) + e\^{}\(λ_{2}\)) is increasing in both λ's. The bounds are A = lse\_\allowbreak{}\allowbreak{}lo(\(A_{1}\), \(A_{2}\)) and B = lse\_\allowbreak{}\allowbreak{}hi(\(B_{1}\), \(B_{2}\)), written out in the code as m + log1p(exp(min − m)) (P5). If \(v_{1}\) = 0 the sum is \(v_{2}\) (A = \(A_{2}\), B = \(B_{2}\)), and symmetrically. Z = \(Z_{1}\) ∧ \(Z_{2}\), and O holds when one operand is 0 and the other is exactly 1.
\item
  \textbf{Dominant operand (second level).} Let a be the operand with \(S_{a}\) = +1 and \(M_{a}\) \textgreater{} 2, so that \(λ_{a}\) = ln a \textgreater{} \(e^{2}\) and a \textgreater{} 1. Suppose the other operand b is smaller: either \(B_{b}\) \textless{} \(A_{a}\), or \(S_{b}\) = +1 and \(N_{b}\) \textless{} \(M_{a}\) (then \(λ_{b}\) ≤ \(e^{N_{b}}\) \textless{} \(e^{M_{a}}\) ≤ \(λ_{a}\)). Then a ≤ u = a + b ≤ 2a, so \(λ_{a}\) ≤ ln u ≤ \(λ_{a}\) + ln 2 ≤ \(B_{a}\) + 1. Taking logarithms, ln \(λ_{a}\) ≤ ln ln u ≤ ln(\(λ_{a}\) + 1) ≤ \(μ_{a}\) + 1/\(λ_{a}\) ≤ \(N_{a}\) + \(e^{-M_{a}}\), since ln(1 + 1/λ) ≤ 1/λ and \(λ_{a}\) ≥ \(e^{M_{a}}\). This rule keeps the second level of values such as \({3}^{{4}^{7776}}+5\), whose logarithm is not a double.
\item
  \textbf{1 + t with t small (second level).} If \(v_{1}\) = 1 and t = \(v_{2}\) ≤ \(e^{B_{2}}\) \textless{} \(e^{-1}\), then t − t\(^{2}\)/2 ≤ ln(1 + t) ≤ t gives ln t − t ≤ ln ln(1 + t) ≤ ln t, using ln(1 − t/2) ≥ −t for t ≤ 1. So S = 1, N = \(B_{2}\) and M = addlo(\(A_{2}\), −exp\_\allowbreak{}\allowbreak{}hi(\(B_{2}\))). If t may be 0 without the flag Z, then \(A_{2}\) = −∞ and M = −∞, so the record only states ln u ≥ 0 (C.4.1); if t = 0 is known (flag Z), the flag O takes over (C.4.7).
\item
  \textbf{1 + t with t below the double range (third level).} If moreover \(B_{2}\) ≤ −\(10^{308}\), \(S_{2}\) = −1 and \(μ_{t}\) = ln\textbar ln t\textbar{} ∈ {[}\(M_{2}\), \(N_{2}\){]} (finite), then −ln ln(1 + t) ∈ {[}e\^{}\(μ_{t}\), e\^{}\(μ_{t}\) + t{]}. Hence ν = ln(−ln ln(1 + t)) ∈ {[}\(μ_{t}\), \(μ_{t}\) + t·\(e^{-μ_{t}}\){]} ⊂ {[}pred(\(M_{2}\)), succ(\(N_{2}\)){]}, because t·\(e^{-μ_{t}}\) \textless{} \(e^{-10^{308}}\) is far below one ulp of \(μ_{t}\) ≥ 709.
\end{enumerate}

\subsubsection*{C.4.4 Differences (form 1, \textbar s − l\textbar)}\addcontentsline{toc}{subsubsection}{C.4.4 Differences (form 1, \textbar s − l\textbar)}

\begin{enumerate}
\def\labelenumi{\arabic{enumi}.}
\item
  \textbf{Upper bound.} \textbar{}\(v_{1}\) − \(v_{2}\)\textbar{} ≤ max(\(e^{B_{1}}\) − \(e^{A_{2}}\), \(e^{B_{2}}\) − \(e^{A_{1}}\)). Each term is ln(\(e^{x}\) − \(e^{y}\)) = x + ln(1 − \(e^{y - x}\)), bounded above with x an upper and y a lower bound. When y ≥ x the term is clipped to x + ln(1 − 1), and ln(0) is replaced by an upper bound; this overestimates, so it is safe (P2, P4).
\item
  \textbf{Lower bound.} In general A = −∞, since the difference may be 0. If \(B_{2}\) \textless{} \(A_{1}\) (certainly \(v_{2}\) \textless{} \(v_{1}\)), then \(v_{1}\) − \(v_{2}\) ≥ \(e^{A_{1}}\) − \(e^{B_{2}}\), and A = addlo(\(A_{1}\), log1m\_\allowbreak{}\allowbreak{}lo(exp\_\allowbreak{}\allowbreak{}hi(addhi(\(B_{2}\), −\(A_{1}\))))). Symmetrically if \(B_{1}\) \textless{} \(A_{2}\).
\item
  \textbf{Exact rationals.} If both operands are exact reduced rationals \(p_{1}\)/\(q_{1}\) and \(p_{2}\)/\(q_{2}\) (passed in R), then they are equal iff (\(p_{1}\), \(q_{1}\)) = (\(p_{2}\), \(q_{2}\)), and the difference is then 0 (Z). Otherwise \textbar{}\(p_{1}\)/\(q_{1}\) − \(p_{2}\)/\(q_{2}\)\textbar{} = \textbar{}\(p_{1}\)\(q_{2}\) − \(p_{2}\)\(q_{1}\)\textbar/(\(q_{1}\)\(q_{2}\)) ≥ 1/(\(q_{1}\)\(q_{2}\)), and A = −ln\_\allowbreak{}\allowbreak{}hi(succ(fl(\(q_{1}\))·fl(\(q_{2}\)))). When \(q_{1}\), \(q_{2}\) \textless{} \(2^{53}\) the argument of ln\_\allowbreak{}\allowbreak{}hi is at least \(q_{1}\)\(q_{2}\). Otherwise the conversion of \(q_{1}\) or \(q_{2}\) to a double may round down, and the argument may fall below \(q_{1}\)\(q_{2}\) by a relative amount below 2.3·\(10^{-16}\). But then \(q_{1}\)\(q_{2}\) \textgreater{} \(2^{53}\) and ln(\(q_{1}\)\(q_{2}\)) \textgreater{} 36, and at such arguments ln\_\allowbreak{}\allowbreak{}hi exceeds the logarithm by more than \(10^{-15}\): its exponent is rounded up twice (Lemma 7), by at least a quarter of an ulp of a number above 32. So A is still a lower bound of ln\textbar{}\(p_{1}\)/\(q_{1}\) − \(p_{2}\)/\(q_{2}\)\textbar.
\item
  \textbf{Near 1.} Write \textbar{}\(e^{a}\) − \(e^{b}\)\textbar{} = \(e^{c}\)·\textbar{}\(e^{d}\) − 1\textbar, where d = a − b and c is either a or b. For every real d, (\(e^{d}\) − 1)/d lies between \(e^{-|d|}\) and \(e^{|d|}\), so ln\textbar{}\(e^{d}\) − 1\textbar{} ∈ {[}ln\textbar d\textbar{} − \textbar d\textbar, ln\textbar d\textbar{} + \textbar d\textbar{]}. The code uses the margin 2\textbar d\textbar{} ≥ \textbar d\textbar{} and applies the rule only when \textbar d\textbar{} \textless{} 1/2, so the guard ln(1/2) needs no rounding care. Two cases use it:

  \begin{itemize}
  \tightlist
  \item
    \(v_{1}\) = 1 and \(v_{2}\) near 1: d = \(λ_{2}\), ln\textbar d\textbar{} ∈ {[}\(M_{2}\), \(N_{2}\){]}, \textbar d\textbar{} ≤ \(e^{N_{2}}\). The base c ∈ {[}min(\(A_{2}\), 0), max(\(B_{2}\), 0){]} covers both c = 0 and c = \(λ_{2}\).
  \item
    Both near 1, with logarithms of opposite signs: \textbar d\textbar{} = \textbar{}\(λ_{1}\)\textbar{} + \textbar{}\(λ_{2}\)\textbar, so ln\textbar d\textbar{} ∈ {[}lse\_\allowbreak{}\allowbreak{}lo(\(M_{1}\), \(M_{2}\)), lse\_\allowbreak{}\allowbreak{}hi(\(N_{1}\), \(N_{2}\)){]} and c ∈ {[}min(\(A_{1}\), \(A_{2}\)), max(\(B_{1}\), \(B_{2}\)){]}.
  \end{itemize}

  The bounds are A ≥ c\_\allowbreak{}\allowbreak{}lo + ln\textbar d\textbar\_\allowbreak{}\allowbreak{}lo − 2e\^{}(\(N_{d}\)) and B ≤ c\_\allowbreak{}\allowbreak{}hi + ln\textbar d\textbar\_\allowbreak{}\allowbreak{}hi + 2e\^{}(\(N_{d}\)), computed with addlo/addhi and exp\_\allowbreak{}\allowbreak{}hi.
\item
  \textbf{Exact cases.} \textbar0 − v\textbar{} = v. The difference is 0 if both operands are 0, both are exactly 1, or they are equal exact rationals; it is exactly 1 if one operand is 0 and the other 1.
\item
  \textbf{1 − t with t below the double range (third level).} If \(v_{1}\) = 1 and t = \(v_{2}\) ≠ 0 has \(B_{2}\) ≤ −\(10^{308}\) and \(μ_{t}\) ∈ {[}\(M_{2}\), \(N_{2}\){]}, then −ln(1 − t) ∈ {[}t, t + t\(^{2}\){]}. As in C.4.3(4), ν ∈ {[}\(μ_{t}\) − t·\(e^{-μ_{t}}\), \(μ_{t}\){]} ⊂ {[}pred(\(M_{2}\)), succ(\(N_{2}\)){]}.
\item
  \textbf{Dominant operand (second level).} Let a satisfy \(S_{a}\) = 1 and \(M_{a}\) \textgreater{} 2, and let the other operand satisfy b ≤ a/e, i.e.~\(λ_{a}\) − \(λ_{b}\) ≥ 1: either \(B_{b}\) \textless{} \(A_{a}\) − 1, or \(S_{b}\) = 1, \(N_{b}\) \textless{} \(M_{a}\) and L = addlo(\(M_{a}\), log1m\_\allowbreak{}\allowbreak{}lo(exp\_\allowbreak{}\allowbreak{}hi(addhi(\(N_{b}\), −\(M_{a}\))))) ≥ 0. L is a lower bound of ln(\(e^{M_{a}}\) − \(e^{N_{b}}\)) ≤ ln(\(λ_{a}\) − \(λ_{b}\)), so L ≥ 0 gives \(λ_{a}\) − \(λ_{b}\) ≥ 1. Then a(1 − 1/e) ≤ u ≤ a, so \(λ_{a}\) − 1 ≤ ln u ≤ \(λ_{a}\), using ln(1 − 1/e) \textgreater{} −1. Hence ln ln u ∈ {[}ln(\(λ_{a}\) − 1), ln \(λ_{a}\){]} ⊂ {[}\(μ_{a}\) − 2/\(λ_{a}\), \(μ_{a}\){]}, by ln(1 − x) ≥ −2x for x ≤ 1/2. The code sets M = addlo(\(M_{a}\), −2·exp\_\allowbreak{}\allowbreak{}hi(−\(M_{a}\))), N = \(N_{a}\), A = \(A_{a}\) − 1 and B = \(B_{a}\). This is the rule for \textbar s − l\textbar{} with s huge and l moderate.
\item
  \textbf{1 − t with t small (second level).} If \(v_{1}\) = 1 and t = \(v_{2}\) ≤ \(e^{B_{2}}\) \textless{} \(e^{-1}\), then u = 1 − t \textless{} 1 and t ≤ −ln u ≤ t + t\(^{2}\). So S = −1, M = \(A_{2}\) and N = addhi(\(B_{2}\), exp\_\allowbreak{}\allowbreak{}hi(\(B_{2}\))), from ln(1 + t) ≤ t. Then ln u ∈ {[}−\(e^{N}\), −\(e^{M}\){]}.
\end{enumerate}

\subsubsection*{C.4.5 Powers (forms 5--8)}\addcontentsline{toc}{subsubsection}{C.4.5 Powers (forms 5--8)}

Forms 7 and 8 are forms 5 and 6 with the roles exchanged; below, x is the base and y the exponent. ln(\(x^{±y}\)) = ±y·ln x, with y = \(e^{λ_{y}}\) ∈ {[}exp\_\allowbreak{}\allowbreak{}lo(\(A_{y}\)), exp\_\allowbreak{}\allowbreak{}hi(\(B_{y}\)){]} and ln x ∈ {[}\(A_{x}\), \(B_{x}\){]}.

\begin{enumerate}
\def\labelenumi{\arabic{enumi}.}
\item
  \textbf{Value.} ±y·ln x is bilinear in (y, ln x), so on the box its extremes are at the four corners (Lemma 5). The bounds are the predecessor of the minimum and the successor of the maximum of the four rounded products. This is valid because, if p\_\allowbreak{}\allowbreak{}j has the smallest rounded value, then pred(fl(p\_\allowbreak{}\allowbreak{}j)) ≤ pred(fl(p\_\allowbreak{}\allowbreak{}k)) ≤ p\_\allowbreak{}\allowbreak{}k for every k. A NaN corner (0·∞) gives the whole line.
\item
  \textbf{Second level.} ln\textbar ln u\textbar{} = ln y + ln\textbar ln x\textbar{} = \(λ_{y}\) + \(μ_{x}\), with sign ±sign(ln x). So M = addlo(\(A_{y}\), \(M_{x}\)) and N = addhi(\(B_{y}\), \(N_{x}\)). This is what evaluates \({(1+\varepsilon)}^{N}\) for ε near \(10^{-700}\).
\item
  \textbf{Base at the third level, deep exponent.} Suppose ln\textbar ln x\textbar{} = −\(e^{ν}\) with ν ∈ {[}M3, N3{]}, and the exponent has \(λ_{y}\) \textgreater{} 700 and \(μ_{y}\) = ln \(λ_{y}\) ∈ {[}\(M_{y}\), \(N_{y}\){]}. Then ln\textbar ln u\textbar{} = \(e^{μ_{y}}\) − \(e^{ν}\).

  \begin{itemize}
  \tightlist
  \item
    If \(M_{y}\) \textgreater{} N3 + 1, this quantity lies in {[}\(e^{M_{y}}\)(1 − \(e^{N3 - M_{y}}\)), \(e^{N_{y}}\){]} and is positive.
  \item
    If \(N_{y}\) \textless{} M3 − 1, it is at most −\(e^{M3}\)(1 − \(e^{N_{y} - M3}\)) \textless{} 0, so \textbar ln u\textbar{} ≤ \(e^{-Λ}\) with Λ that bound.
  \item
    Otherwise the sign is left undecided (S = 0) and only the corner bounds of (1) remain.
  \end{itemize}

  The products are rounded down (pred of a product of lower bounds of positive numbers), and ln u = S·\(e^{ln|ln u|}\) is then bounded with P2.
\item
  \textbf{Exponent below \(e^{-10^{308}}\), base with known \(μ_{x}\) (third level).} Let ln y = −\(e^{μ_{y}}\) with \(μ_{y}\) ∈ {[}\(M_{y}\), \(N_{y}\){]} (so \(μ_{y}\) ≥ ln \(10^{308}\) ≈ 709.2), and let \(μ_{x}\) ∈ {[}\(M_{x}\), \(N_{x}\){]}. Then ln\textbar ln u\textbar{} = −(\(e^{μ_{y}}\) − \(μ_{x}\)) = −\(e^{ν}\) with ν = ln(\(e^{μ_{y}}\) − \(μ_{x}\)). This is defined when \(e^{μ_{y}}\) \textgreater{} \(μ_{x}\), which the condition \(μ_{y}\) \textgreater{} ln max(\textbar{}\(μ_{x}\)\textbar, 1) + 1 ensures (and then \textbar{}\(μ_{x}\)\textbar·\(e^{-μ_{y}}\) \textless{} \(e^{-1}\)). With \(t_{1}\) ≥ \(N_{x}\)\(^{+}\)·\(e^{-M_{y}}\) and \(t_{2}\) ≥ \(M_{x}\)\(^{-}\)·\(e^{-N_{y}}\) (products rounded up):

  \begin{itemize}
  \tightlist
  \item
    ν ≥ \(M_{y}\) + ln(1 − \(t_{1}\)), giving M3 = addlo(\(M_{y}\), log1m\_\allowbreak{}\allowbreak{}lo(\(t_{1}\)));
  \item
    ν ≤ \(N_{y}\) + ln(1 + \(t_{2}\)), giving N3 = addhi(\(N_{y}\), log1p\_\allowbreak{}\allowbreak{}hi(\(t_{2}\)));
  \item
    ln\textbar ln u\textbar{} ∈ {[}−exp\_\allowbreak{}\allowbreak{}hi(N3), −exp\_\allowbreak{}\allowbreak{}lo(M3){]}, which is used for M and N.
  \end{itemize}
\item
  \textbf{Exact cases.} \({1}^{y}\) = 1 (\(O_{1}\)). x\(^{0}\) = 1 for x ≠ 0. \({0}^{y}\) = 0 for y \textgreater{} 0. \({0}^{-y}\) and \(0^{0}\) are undefined (ok false).
\end{enumerate}

\subsubsection*{C.4.6 Rational anchors}\addcontentsline{toc}{subsubsection}{C.4.6 Rational anchors}

The anchor of the result is computed from the anchors of the operands. Values near 1 that carry a second level are first seen as anchors at 1 (\texttt{\_\allowbreak{}\allowbreak{}virtual\_\allowbreak{}\allowbreak{}anchor}): if \textbar ln u\textbar{} ≤ \(e^{N}\) \textless{} \(e^{-20}\), then δ = \(e^{ln u}\) − 1 has the sign of ln u, and ln\textbar δ\textbar{} = ln\textbar ln u\textbar{} + ln((\(e^{λ}\) − 1)/λ) ∈ {[}M − 2e\^{}N, N + 2e\^{}N{]} (C.4.4(4)). Anchors whose \textbar δ\textbar{} may exceed \(e^{-20}\) are dropped. Every rule below therefore has \textbar δ\textbar{} \textless{} 2·\(10^{-9}\), which justifies \textbar ln(1 + δ)\textbar{} ≤ 2\textbar δ\textbar{} and \textbar ln(ln(1 + δ)/δ)\textbar{} ≤ 2\textbar δ\textbar.

\begin{enumerate}
\def\labelenumi{\arabic{enumi}.}
\item
  \textbf{Exact q plus or minus a tiny t} (t not exact, not anchored, \(e^{B_{t}}\) \textless{} q·\(e^{-40}\)): q ± t = q(1 ± t/q), so δ = ±t/q and ln\textbar δ\textbar{} = ln t − ln q ∈ {[}\(A_{t}\) − ln\_\allowbreak{}\allowbreak{}hi q, \(B_{t}\) − ln\_\allowbreak{}\allowbreak{}lo q{]}. The difference is \textbar q − t\textbar{} = q − t because t \textless{} q.
\item
  \textbf{Anchored q(1 + d) plus a tiny t of the same sign:} δ = d ± t/q with both terms of the same sign, so ln\textbar δ\textbar{} = ln(\textbar d\textbar{} + t/q) ∈ {[}lse\_\allowbreak{}\allowbreak{}lo, lse\_\allowbreak{}\allowbreak{}hi{]} (P5).
\item
  \textbf{Exact + exact, \textbar exact − exact\textbar, exact·exact, exact/exact:} exact rationals (\texttt{rat\_\allowbreak{}\allowbreak{}op}, reduced, with an overflow check). A zero difference is the exact zero of C.4.4(3).
\item
  \textbf{Anchored + exact:} \(q_{1}\)(1 + \(d_{1}\)) + \(q_{2}\) = (\(q_{1}\) + \(q_{2}\))(1 + δ) with δ = \(q_{1}\)\(d_{1}\)/(\(q_{1}\) + \(q_{2}\)): the same sign, and ln\textbar δ\textbar{} = ln\textbar{}\(d_{1}\)\textbar{} + ln \(q_{1}\) − ln(\(q_{1}\) + \(q_{2}\)).
\item
  \textbf{\textbar exact − anchored\textbar:} if \(q_{1}\) = \(q_{2}\) = q, then \textbar q − q(1 + d)\textbar{} = q\textbar d\textbar, exactly enclosed by ln q + ln\textbar d\textbar, and it is not 0. Otherwise, with D = \(q_{1}\) − \(q_{2}\) ≠ 0, \textbar{}\(q_{1}\) − \(q_{2}\)(1 + \(d_{2}\))\textbar{} = \textbar D\textbar·(1 − (\(q_{2}\)/D)\(d_{2}\)): an anchor at \textbar D\textbar{} with δ = −(\(q_{2}\)/D)\(d_{2}\), whose sign is −sign(D)·sign(\(d_{2}\)), and ln\textbar δ\textbar{} = ln \(q_{2}\) − ln\textbar D\textbar{} + ln\textbar{}\(d_{2}\)\textbar.
\item
  \textbf{Exact × anchored, anchored / exact:} \(q_{1}\)·\(q_{2}\)(1 + d) and (\(q_{1}\)/\(q_{2}\))(1 + d) keep δ. \textbf{Exact / anchored:} \(q_{1}\)/(\(q_{2}\)(1 + d)) = (\(q_{1}\)/\(q_{2}\))(1 + δ) with δ = −d/(1 + d): the opposite sign, and ln\textbar δ\textbar{} = ln\textbar d\textbar{} − ln(1 + d) ∈ ln\textbar d\textbar{} ± 2\textbar d\textbar.
\item
  \textbf{Exact base, exponent anchored at an integer k ≤ 6:} q\^{}(±k(1+d)) = \(q^{±k}\)·\(e^{z}\) with z = ±k·d·ln q. So δ = \(e^{z}\) − 1 has sign ±sign(d)·sign(ln q), and ln\textbar δ\textbar{} = ln k + ln\textbar d\textbar{} + ln\textbar ln q\textbar{} + θ with \textbar θ\textbar{} ≤ \textbar z\textbar{} ≤ \(e^{dn}\). The three-term sums are outward (\texttt{addlo}/\texttt{addhi}).
\item
  \textbf{Exact or anchored base, exact integer exponent k ≥ 1.} (q(1 + d))\^{}(±k) = \(q^{±k}\)·(1 + d)\^{}(±k). For d = 0 this is the exact rational \(q^{±k}\) = \(P^{k}\)/\(Q^{k}\) (or \(Q^{k}\)/\(P^{k}\)), kept when it is below \(2^{60}\) (\texttt{ipow}, exact int64 powers). It is reduced because gcd(P, Q) = 1. For d ≠ 0 it is an anchor at \(q^{±k}\) with δ = \(e^{z}\) − 1, z = ±k·ln(1 + d): sign δ = ±sign d, and ln\textbar δ\textbar{} = ln k + ln\textbar d\textbar{} + ln(ln(1 + d)/d) + ln((\(e^{z}\) − 1)/z), where the last two terms have absolute value at most 2\textbar d\textbar{} and \textbar z\textbar{} ≤ 2k\textbar d\textbar. Hence ln\textbar δ\textbar{} ∈ {[}ln k + DM − 2(k + 1)\(e^{DN}\), ln k + DN + 2(k + 1)\(e^{DN}\){]}, with the margin rounded up. \emph{Why it is needed:} in n = 8 (any order), a chain such as \textbar{}\(8-{(.5+{.4}^{{7}^{6}})}^{-3}\)\textbar{} lost its anchor at \({(.5+\varepsilon)}^{-3}\). The difference then had the enclosure {[}0, 8·\(2^{-52}\){]} instead of \(48\cdot {.4}^{{7}^{6}}\), and the next power, (·)\^{}(−.01), produced the hits ln u ∈ {[}0.305, ∞) of split 8.
\item
  \textbf{Exact powers of rationals: the factored channel (FP, FQ, FE, FD).} It states u = (FP/FQ)\^{}(FE/FD) exactly, where FP/FQ \textgreater{} 1 is a reduced rational that is not a perfect power (\texttt{primitive}: P/Q = \(r^{m}\), with every root checked by exponentiation) and FE/FD is a nonzero reduced fraction, FD ≥ 1. The channel is set or used in these cases:

  \begin{itemize}
  \tightlist
  \item
    a power of an exact rational q = \(r^{m}\) to an exact positive rational ±p/q′ gives \(r^{±mp/q′}\), unless the result is itself an exact rational (rule 8);
  \item
    a power of \(r^{e}\) to an exact positive rational ±p/q′ gives \(r^{±ep/q′}\). This holds because (\(r^{e}\))\^{}t = \(r^{et}\) for r \textgreater{} 0 and real t;
  \item
    for products and quotients of two values of the same r (a triple, or an exact rational, taken through its own primitive form \(r^{m}\)), \(r^{e_{1}}\)·\(r^{e_{2}}\) = \(r^{e_{1}+e_{2}}\) and \(r^{e_{1}}\)/\(r^{e_{2}}\) = \(r^{e_{1}-e_{2}}\); an exponent 0 means the value is exactly 1;
  \item
    for \textbar s − l\textbar, the same triple on both sides means the same value, hence exactly 0; an exact rational on one side is compared through its primitive triple;
  \item
    for \textbar s − l\textbar{} with two \emph{different} triples the values are different, so x ≠ y. Then L = ln y − ln x = \(e_{2}\) ln \(r_{2}\) − \(e_{1}\) ln \(r_{1}\) ≠ 0 is computed with Arb balls, at 256 bits and more if the ball contains 0; it cannot stay at 0, since L ≠ 0. This gives ln\textbar x − y\textbar{} = max(ln x, ln y) + ln(1 − \(e^{-|L|}\)), rounded outward (\texttt{\_\allowbreak{}\allowbreak{}factor\_\allowbreak{}\allowbreak{}gap}). For instance \textbar{}\({5}^{4096}\) − v\textbar{} with v = \(e^{6592.26…}\) from other digits is then known to be of size \(e^{6592}\) and not ``anything in {[}0, e\^{}6592{]}''.
  \end{itemize}

  The orientation is canonical, r = FP/FQ \textgreater{} 1: a base r \textless{} 1 is written (1/r)\^{}(−e). Without it, \(2^{4^{6}}\) and \({(.5)}^{-4096}\) would get different triples for the same value.

  \emph{Why the triple is canonical:} suppose \(r_{1}\)\^{}(\(e_{1}\)) = \(r_{2}\)\^{}(\(e_{2}\)) with \(r_{1}\), \(r_{2}\) \textgreater{} 1 primitive and \(e_{1}\), \(e_{2}\) rational. Multiply both exponents by a common denominator D; then \(r_{1}\)\^{}(D\(e_{1}\)) = \(r_{2}\)\^{}(D\(e_{2}\)) with integer exponents. By unique factorization \(r_{1}\) = \(ρ^{a}\) and \(r_{2}\) = \(ρ^{b}\) with D\(e_{1}\)a = D\(e_{2}\)b, and primitivity forces a = b = 1, so \(r_{1}\) = \(r_{2}\) and \(e_{1}\) = \(e_{2}\). All operations are exact integer arithmetic, and a triple is kept only while \textbar FE\textbar, FD \textless{} \(2^{62}\). \emph{Why it is needed:} \(2^{8^{5}}\) and 4\^{}((1+3)\(^{7}\)) are both \({2}^{32768}\) ≈ e\^{}22713, and their double enclosures can only show that the difference lies in {[}0, e\^{}22713{]}; this produced the hits of n = 8 (any order) of type ln u ∈ {[}ln .6, 22 691{]}. Test \texttt{t\_\allowbreak{}\allowbreak{}factor.py}: 6000 powers. All 159 pairs of equal values written differently are detected, including inverted bases. For the 2841 different pairs, the enclosure of \textbar x − y\textbar{} contains the exact value computed with mpmath at 80 digits. There are 0 false equalities and 0 failures. \emph{Rational exponents:} hitcheck found in split 129 of n = 9 hits that \texttt{lvl.py} left undecided, all of them \textbar x − y\textbar{} of two equal powers. Examples are (\(2^{1000}\))\^{}(2/5) = \(2^{400}\) = (\(2^{25}\))\(^{16}\) and (\(3^{10000}\))\^{}(1/5) = \(3^{2000}\) = (\(3^{16}\))\(^{125}\): the exponent 2/5 or 1/5 (= .4, .2) was not an integer, so no triple was set. With rational exponents both sides get the triple (2, 1, 400) or (3, 1, 2000), and the difference is the exact zero. Test \texttt{t\_\allowbreak{}\allowbreak{}factor2.py}: 6000 powers built as (b\^{}\(t_{1}\))\^{}\(t_{2}\) with rational t's. The 300 pairs of equal values written differently are all detected (zero difference, quotient exactly 1). All 2700 different pairs have a lower bound that contains the mpmath value at 120 digits. In the 1348 products with an exact rational of the same base, the triple is the exact product. 0 failures.
\item
  \textbf{Use of the anchor.} ln u = ln q + ln(1 + δ) ∈ {[}ln\_\allowbreak{}\allowbreak{}lo q − 2e\^{}(DN), ln\_\allowbreak{}\allowbreak{}hi q + 2e\^{}(DN){]}, which tightens A and B. An anchor at 1 gives the second level: S = DS, ln\textbar ln u\textbar{} = ln\textbar δ\textbar{} + ln(ln(1 + δ)/δ) ∈ {[}DM − 2e\^{}(DN), DN + 2e\^{}(DN){]}. An exact anchor fixes Z = false and O = (AP = AQ).
\end{enumerate}

\subsubsection*{C.4.7 The final step}\addcontentsline{toc}{subsubsection}{C.4.7 The final step}

The flags override the intervals: Z gives A = B = −∞ and O gives A = B = 0, and S = 0 for both. Then \texttt{tighten} intersects {[}A, B{]} with S·{[}\(e^{M}\), \(e^{N}\){]} (P2). Each field holds for every value, so the intersection does too. ∎ (Lemma X2)

\subsubsection*{C.4.8 Preimage windows: proof of Lemma X3}\addcontentsline{toc}{subsubsection}{C.4.8 Preimage windows: proof of Lemma X3}

\textbf{Lemma X3 (preimage windows).} \texttt{pre(f,\ X,\ {[}α,\ β{]})} returns windows (in ln l) containing every l with f(s, l) ∈ e\footnote{α, β} for some s ∈ X. The window {[}−∞, −∞{]} (only l = 0) is returned only where f(s, 0) is defined and may lie in the window.

\emph{Proof.} Write a = \(A_{s}\) and b = \(B_{s}\).

\begin{itemize}
\tightlist
\item
  \textbf{s·l:} \(λ_{l}\) = ln u − \(λ_{s}\) ∈ {[}α − b, β − a{]}. If s = 0, then s·l = 0 for every l, so the window is everything if 0 ∈ window (α = −∞) and empty otherwise.
\item
  \textbf{s/l:} \(λ_{l}\) = \(λ_{s}\) − ln u ∈ {[}a − β, b − α{]}. \textbf{l/s:} \(λ_{l}\) ∈ {[}α + a, β + b{]}. l/0 is undefined.
\item
  \textbf{s + l:} l ≤ \(e^{β}\) − s ≤ \(e^{β}\) − \(e^{a}\) and l ≥ \(e^{α}\) − \(e^{b}\), so the window is {[}ldiff\_\allowbreak{}\allowbreak{}lo(α, b), ldiff\_\allowbreak{}\allowbreak{}hi(β, a){]}. If s itself may lie in the window, l may be 0 or arbitrarily small, and the lower end is −∞. If s = 0, the window is {[}α, β{]}.
\item
  \textbf{\textbar s − l\textbar:} l = s + d or l = s − d with d ∈ {[}\(e^{α}\), \(e^{β}\){]}. The first branch is {[}lse\_\allowbreak{}\allowbreak{}lo(a, α), lse\_\allowbreak{}\allowbreak{}hi(b, β){]}. The second has upper end ln(\(e^{b}\) − \(e^{α}\)), and the code takes the maximum with ln(\(e^{β}\) − \(e^{a}\)), which only enlarges it. Its lower end is ln(\(e^{a}\) − \(e^{β}\)) when a \textgreater{} β, and −∞ otherwise. When b \textless{} α the code still gives this branch the lower end ln(\(e^{α}\) − \(e^{b}\)): the branch is then empty, and the extra window is harmless. If s may lie in the window, l = 0 is admitted.
\item
  \textbf{\(s^{l}\), s\^{}−l:} ±l·\(λ_{s}\) ∈ {[}α, β{]}, that is, l·\(λ_{s}\) ∈ {[}α', β'{]} with {[}α', β'{]} = {[}α, β{]} for \(s^{l}\) and {[}−β, −α{]} for s\^{}−l. So ln l lies in the window of \{y \textgreater{} 0 : y·μ ∈ {[}α', β'{]} for some μ ∈ {[}a, b{]}\} (\texttt{\_\allowbreak{}\allowbreak{}ratio\_\allowbreak{}\allowbreak{}window}). For a \textgreater{} 0 (μ \textgreater{} 0) it is y ∈ {[}α'/b, β'/a{]} when α' \textgreater{} 0 and (0, β'/a{]} otherwise, empty if β' ≤ 0. For b \textless{} 0 it is the mirror image. If 0 ∈ {[}a, b{]}, every y. Divisions are rounded outward and P3 gives the logarithms. Further cases:

  \begin{itemize}
  \tightlist
  \item
    l = 0 gives s\(^{0}\) = 1, admitted when 0 ∈ {[}α, β{]};
  \item
    s = 1 gives 1 for every l;
  \item
    s = 0 gives 0 for l \textgreater{} 0 (form 5).
  \end{itemize}
\item
  \textbf{\(l^{s}\), l\^{}−s:} ±s·\(λ_{l}\) ∈ {[}α, β{]} with s ∈ {[}\(e^{a}\), \(e^{b}\){]}, so \(λ_{l}\) ∈ {[}α'/s, β'/s{]} over that range, with {[}α', β'{]} as above. The extremes are α'/\(e^{b}\) or α'/\(e^{a}\) (according to the sign of α'), and likewise for β'. With s = 0, l\(^{0}\) = 1 for every l; with l = 0, 0\^{}s = 0 (form 7).
\item
  \textbf{The window ``only l = 0''} ({[}−∞, −∞{]}) is kept only where f(s, 0) is defined and may lie in the window (\texttt{pre}): s + 0 = s, \textbar s − 0\textbar{} = s, s·0 = 0, 0/s = 0 (s ≠ 0), s\(^{0}\) = 1, 0\^{}s = 0. It is removed for s/0 and \(0^{-s}\).
\end{itemize}

All the bounds are monotone in a, b, α and β, evaluated with P1--P5, so the windows are outward. ∎

\subsubsection*{C.4.9 Tests}\addcontentsline{toc}{subsubsection}{C.4.9 Tests}

\begin{itemize}
\tightlist
\item
  \texttt{t\_\allowbreak{}\allowbreak{}fwd.py} (fwd against mpmath at 60 digits on moderate, near-1, extreme and beyond-double values), \texttt{t\_\allowbreak{}\allowbreak{}anchor.py} (anchors, third level; mpmath at 520 digits, 3000 cases) and \texttt{t\_\allowbreak{}\allowbreak{}pre.py} (1.8·\(10^{6}\) random cases): 0 failures.
\item
  \texttt{t\_\allowbreak{}\allowbreak{}index.py}: the new index (C.9) returns exactly the same pairs as the original one on random tables with every kind of entry, and exactly the pairs of the plain test on tables of one-signed entries of every magnitude, with signed zeros, subnormals and NaN (version ixT2).
\item
  Further tests, all with 0 failures (Table C1).
\end{itemize}

\textbf{Table C1.} Further tests of motor X.

{\small
\begin{longtable}[]{@{}>{\raggedright\arraybackslash}p{(\linewidth - 4\tabcolsep) * \real{0.2156}}>{\raggedright\arraybackslash}p{(\linewidth - 4\tabcolsep) * \real{0.5298}}>{\raggedright\arraybackslash}p{(\linewidth - 4\tabcolsep) * \real{0.2546}}@{}}
\toprule\noalign{}
Test & What it checks & Cases \\
\midrule\noalign{}
\endhead
\bottomrule\noalign{}
\endlastfoot
\texttt{t\_\allowbreak{}\allowbreak{}fwd2.py} & sums, differences and products beyond the double range with close second levels, against mpmath at 80 digits & 20 000 \\
\texttt{t\_\allowbreak{}\allowbreak{}pre2.py} & second-level preimages of \(s^{l}\), s\^{}−l & 14 007 \\
\texttt{t\_\allowbreak{}\allowbreak{}level2.py} & second-level constraint on the bases of \(l^{s}\) & 20 000 \\
\texttt{t\_\allowbreak{}\allowbreak{}unique.py} & one logarithm per distinct window bound, bit-identical & 200 batches \\
\texttt{t\_\allowbreak{}\allowbreak{}near1.py} & near-1 rows against motor X's own recursion; the index; the rows against their writings & 441 784 values \\
\texttt{t\_\allowbreak{}\allowbreak{}factor2.py} & triples with rational exponents, equal pairs, products with rationals & 6000 powers \\
\texttt{t\_\allowbreak{}\allowbreak{}offset.py} & exact offsets X + c & 19 950 \\
\texttt{t\_\allowbreak{}\allowbreak{}lvlpf.py} & lvl: exact powers with rational exponents, rational and irrational offsets & 2000 + 704 + 176 \\
\end{longtable}
}

\subsection*{C.5 Tables and exact decisions}\addcontentsline{toc}{subsection}{C.5 Tables and exact decisions}

\textbf{V-tables.} For \textbar S\textbar{} ≤ t\_\allowbreak{}\allowbreak{}V, the V-tables are the classes of motor A's table T(S), each converted to an enclosure; by Theorem A they cover V(S). Exact rationals, the exact 0 and the exact 1 are flagged; the rationals of tables of size 5 are kept too.

\textbf{NE-tables.} For \textbar S\textbar{} ≤ t\_\allowbreak{}\allowbreak{}N, N̂(S) is a finite list of enclosures built by Lemma X1 from the tables of the parts, with all forms over all splits, and without any test of whether a value lies in V(S). It consists of:

\begin{itemize}
\item
  \begin{enumerate}
  \def\labelenumi{(\alph{enumi})}
  \tightlist
  \item
    N̂(\(S_{1}\)) × (V̂ ∪ N̂)(\(S_{2}\));
  \end{enumerate}
\item
  \begin{enumerate}
  \def\labelenumi{(\alph{enumi})}
  \setcounter{enumi}{1}
  \tightlist
  \item
    V̂(\(S_{1}\)) × N̂(\(S_{2}\));
  \end{enumerate}
\item
  \begin{enumerate}
  \def\labelenumi{(\alph{enumi})}
  \setcounter{enumi}{2}
  \tightlist
  \item
    V̂(\(S_{1}\)) × V̂(\(S_{2}\)), keeping only the results that may be out of range.
  \end{enumerate}
\end{itemize}

\textbf{Exact decisions (motor B).} Three situations are decided with the exact members of the classes involved (Lemmas 13 and 14) instead of their enclosures:

\begin{enumerate}
\def\labelenumi{\arabic{enumi}.}
\tightlist
\item
  \textbf{Overlapping differences.} \textbar s − l\textbar{} of two overlapping in-range classes is needed only if it may be out of range. It is decided exactly: an exact zero, or an in-range difference, is not a new value. An out-of-range nonzero difference is replaced by its exact ball. This also decides coincidences between irrational values, such as \({4}^{-.3}\) = \({.5}^{.6}\) = \({2}^{-3/5}\).
\item
  \textbf{Classes that may be 1.} A V-class whose enclosure of ln contains 0 without being the exact 1 is replaced by its exact members: the exact 1, or balls with the sign of ln known.
\item
  \textbf{Bases that may be 1.} In a power with a huge exponent, an in-range base whose enclosure contains 1 gives no information. Its exact members are computed (B.query on its digit set and enclosure); each is either the exact 1, giving the value 1, or a ball with \textbar ln\textbar{} \textgreater{} 0, which is carried through \texttt{fwd}.
\end{enumerate}

\textbf{Lemma X4.} For every S with \textbar S\textbar{} ≤ t\_\allowbreak{}\allowbreak{}N, every value of X(S) lies in an enclosure of N̂(S).

\emph{Proof:} induction on \textbar S\textbar, with Lemmas X1 and X2, the covering property of V̂ (Theorem A) and the correctness of the exact decisions (Lemmas 13 and 14). In case (c) of Lemma X1 the value f(s, l) is out of range, so its enclosure may be out of range and is kept. The exact decisions on overlapping differences drop only an exact zero or an in-range difference, which are values of a node that is 0 or in range, never case (c). Merging equal enclosures loses no value. ∎

\subsection*{C.6 Window queries}\addcontentsline{toc}{subsection}{C.6 Window queries}

Q(S, kind, w) returns enclosures covering the values of the given kind (V, NE or both) that lie in the window w:

\begin{itemize}
\tightlist
\item
  for small S, by lookup in the sorted tables;
\item
  otherwise, by recursion over the splits S = \(S_{1}\) ⊔ \(S_{2}\) with \(S_{1}\) the table side. For every s in the \(S_{1}\)-table and every form f, \texttt{pre} gives the windows for l, Q(\(S_{2}\), ·) answers them, and \texttt{fwd} keeps the pairs whose value meets w.
\end{itemize}

For kind NE the pairs are those of Lemma X1: (NE, any), (V, NE), and (V, V) with an out-of-range result; the NE side is answered by the NE-tables or by recursion. For kind V they are (V, V) with an in-range result, plus the literals of S. Overlapping (V, V) differences that may be out of range, and in-range bases that may be 1, are decided exactly, as in C.5.

\textbf{Theorem X.} For every S and window w, every value of X(S) in w lies in an enclosure returned by Q(S, NE, w), and every value of V(S) in w lies in an enclosure returned by Q(S, V, w). \emph{Proof:} induction on \textbar S\textbar. The base cases are the tables (Theorem A and Lemma X4). The induction step uses Lemma X1 (every value arises from a split, with one side from a table: in any order the smaller side, of at most 4 digits, and in order the shorter part), Lemma X3 (the partner lies in a preimage window), the induction hypothesis for \(S_{2}\), and Lemma X2 (the forward enclosure meets w). Every filter of the query removes a candidate only when the value \emph{of that node} is certainly not what the case needs: outside the window or its second-level constraint, an in-range result of a (V, V) pair in a query of kind NE (case (c)), an exact zero or an in-range difference of a (V, V) pair, or a base whose exact members give no value in the window. No step tests whether a value belongs to V(S). ∎

\textbf{Proof of Theorem 2.} Run Q(\{1,\ldots,n\}, NE, W) with W = {[}c − r, c + r{]}, where r is the certified optimum E* enlarged by \(10^{-9}\) relative or, for n = 9 in any order, an upper bound of E* (Section 6). A larger window only strengthens the conclusion. The query returns a finite list of enclosures (the \emph{hits}). The independent evaluator \texttt{lvl} evaluates every derivation of every hit over the exact members of the table classes, so the value of every expression covered by Theorem X is among the values it evaluates. For every cell listed in Table 6 every evaluated value is certainly outside W. By Theorem X, X(\{1,\ldots,n\}) ∩ W = ∅, and by fact 1, NE(\{1,\ldots,n\}) ∩ W = ∅. So every value within E* of c is in V(\{1,\ldots,n\}), where Theorem 1 applies. In fact more holds: no expression with an out-of-range intermediate comes within E* of c.~∎

\textbf{The variant without unary minus.} Let \(E_{11}\) ≥ E* be the optimum of Section 8 (certified in range), enlarged in the same way, or an independently checked upper bound of it (Section 6). The same query is run with the window of radius max(E*, \(E_{11}\)), and every evaluated value is checked against both windows. If none is undecided or certainly inside W(\(E_{11}\)), then X(\{1,\ldots,n\}) ∩ W(\(E_{11}\)) = ∅. By fact 4, every expression without unary minus within \(E_{11}\) of c then has all its nodes in range, so its value lies in \(V_{11}\)(\{1,\ldots,n\}), where the certification of Section 8 applies. In particular, a value that is within \(E_{11}\) of c and can be written in range only with the unary minus has no writing without unary minus at all. The structural form of Lemma X1 is what makes this work. A test (\texttt{t\_\allowbreak{}\allowbreak{}coinc.py}) checks it in practice. Exact coincidences written through the out-of-range power \({9}^{{8}^{7}}\) are found by the query as exact values: 1 = \({1}^{{9}^{{8}^{7}}}\) over \{1, 7, 8, 9\}, and 2 and 3 over \{1, 2, 7, 8, 9\}, from the 4- and 5-digit tables, and 6 = \(2+\left(3+{1}^{{.9}^{{8}^{7}}}\right)\) over six digits, through the recursion. ∎

\subsection*{C.7 What is trusted, and tests}\addcontentsline{toc}{subsection}{C.7 What is trusted, and tests}

Trusted, in addition to the base of Appendix A.9:

\begin{itemize}
\tightlist
\item
  numpy's IEEE arithmetic (+ − × / and nextafter, which are correctly rounded);
\item
  numpy's \texttt{log} and \texttt{log2}, used only where a relative error below \(2^{-40}\) suffices: to select candidates in the indices of the tables (the selection then keeps a margin of that size) and to estimate sizes before an exact integer power;
\item
  the correctness of the code of \texttt{ivl.py}, \texttt{motorx.py}, \texttt{qx.py}, \texttt{hitcheck.py} and \texttt{lvl.py}; for n = 9 in any order, in each of the five versions of \texttt{hitcheck.py} and \texttt{lvl.py} that checked its records (Appendix C.9).
\end{itemize}

\texttt{lvl} discards a derivation as undefined only when it proves that an operation is undefined (for instance 0 to a negative power); its exact zeros use only the direction ``same form ⇒ same value''.

Independent tests:

\begin{enumerate}
\def\labelenumi{\arabic{enumi}.}
\tightlist
\item
  The tests of \texttt{fwd}, \texttt{pre} and the anchors listed in C.4.9.
\item
  Coverage of an approximate exhaustive enumeration in floating point (\texttt{measure.py}, \texttt{t\_\allowbreak{}\allowbreak{}cover.py}). Every approximate new value for sets of up to 5 digits is covered. The only exceptions are floating-point ghosts such as (4.4·\(10^{-17}\))\^{}(−0.5), where 4.4·\(10^{-17}\) is a rounding residue of an exact zero.
\item
  Agreement of the recursive queries with the materialized tables (\texttt{t\_\allowbreak{}\allowbreak{}q.py}, n = 5: 0 values not covered).
\item
  The forward method (full materialization, \texttt{cells\_\allowbreak{}\allowbreak{}x.py}) and the recursive method (\texttt{cells\_\allowbreak{}\allowbreak{}q.py}) agree on the cells n ≤ 6.
\item
  Structural coverage (\texttt{t\_\allowbreak{}\allowbreak{}coinc.py}): exact coincidences written through an out-of-range node, as in C.6, are found by the query.
\item
  Regression: the cells closed with the first version give the same verdicts with the optimized one.
\item
  Pieces: a split computed in pieces gives the same hits as the whole split (\texttt{t\_\allowbreak{}\allowbreak{}pieces.py}; at n = 6 the independent review of the manuscript also compared 3 whole splits, 4.1 million hits, with their pieces: identical multisets).
\end{enumerate}

\textbf{A latent issue found by the independent review.} The cache of exact members of a table class is keyed by the position of the class. If a class whose logarithm encloses 0 without being the exact 1 were refined (C.5, item 2), positions would shift and a key could designate another class. This needs such an ambiguous class; a count over every V-table used (n = 9 in order up to 5 digits, 1 233 008 classes; n = 7, 8 and 9 in any order up to 4 digits, 2.3, 4.9 and 9.3 million classes) found none, so no computation of Table 6 is affected. The published code is unchanged and the README describes the issue.

The cells, hits and times are those of Table 6 (Section 6).

\subsection*{C.8 Measurements}\addcontentsline{toc}{subsection}{C.8 Measurements}

From an approximate floating-point enumeration of all expressions, with no range limit:

\begin{itemize}
\tightlist
\item
  \textbf{Up to 4 digits:} no in-range value reached only through an extreme intermediate (extreme values such as \({9}^{{8}^{7}}\) do occur).
\item
  \textbf{Five digits:} 20--300 new values per set, all with \textbar lo\(g_{10}\) u\textbar{} ≥ 30.
\item
  \textbf{Digits 1--6:} 41 247 new values. The closest to the constants are at distance 3.8·\(10^{-3}\) (π), 5.9·\(10^{-3}\) (γ), 3.3·\(10^{-2}\) (δ) and 1.1·\(10^{-2}\) (e), against certified optima of order \(10^{-8}\).
\end{itemize}

Most new values are near the range boundaries. The mechanisms that bring a value back near c are:

\begin{enumerate}
\def\labelenumi{\arabic{enumi}.}
\tightlist
\item
  Damping: \({x}^{y}\) with x extreme and y small.
\item
  Pairs of extremes: quotients, or near-cancelling differences.
\item
  Thresholds: values just beyond \(10^{300}\) times moderate factors.
\item
  Shadows: y + ε, and 1 + ε raised to huge powers.
\end{enumerate}

\subsection*{C.9 The optimized motor X}\addcontentsline{toc}{subsection}{C.9 The optimized motor X}

The runs of Table 6 used a version of motor X optimized between 2 and 18 November 2026. With the first version, a single top-level split of n = 8 in any order took 40--70 minutes (127 splits per constant). The optimized version closes the whole cell in about 17 minutes of CPU for the three constants together. The supplementary document \texttt{motorx\_\allowbreak{}\allowbreak{}changes.md} of the repository describes the 50 changes one by one, each with its proof and its test. They are of four kinds:

\begin{enumerate}
\def\labelenumi{\arabic{enumi}.}
\tightlist
\item
  \textbf{Same output.} Changes that return exactly the same set of pairs, each checked by a test of identity:

  \begin{itemize}
  \tightlist
  \item
    tiered indices of the tables, also keyed by ln\textbar·\textbar;
  \item
    the exchange of the roles of the two sides of a split;
  \item
    a persistent cache of the NE-tables, keyed by a SHA-256 fingerprint of their inputs;
  \item
    parallel exp/ln bounds;
  \item
    one logarithm per distinct window bound;
  \item
    grouped exact queries;
  \item
    a single pass with the check of the hits;
  \item
    queries cut in pieces, with the three constants together;
  \item
    bounded lookups;
  \item
    recursive sub-queries that receive the windows of every form at once;
  \item
    repeated windows answered once.
  \end{itemize}
\item
  \textbf{Smaller supersets.} A superset replaced by a smaller one, each by a one-line covering argument:

  \begin{itemize}
  \tightlist
  \item
    the (V, V) part of a query of kind NE restricted to the out-of-range parts of its windows (case (c) of Lemma X1);
  \item
    preimages and constraints at the second level, passed down the recursion;
  \item
    the exact narrowing of bases near 1 to their preimage window (Lemma X3);
  \item
    one table entry per distinct statement.
  \end{itemize}
\item
  \textbf{New exact rules of \texttt{ivl}} (C.4.2--C.4.6):

  \begin{itemize}
  \tightlist
  \item
    the margins of the dominance rules;
  \item
    exact or anchored bases with exact integer or rational exponents;
  \item
    the factored channel (FP/FQ)\^{}(FE/FD) and the exact offsets of exact powers;
  \item
    the exact 1 as the rational 1/1;
  \item
    triples whose value is an exact rational.
  \end{itemize}
\item
  \textbf{Exact decisions with motor B's arithmetic} (motor B itself is unchanged):

  \begin{itemize}
  \tightlist
  \item
    the exact zero in recursive V queries;
  \item
    exact differences of overlapping classes;
  \item
    the exact members of V(S) near 1;
  \item
    the exact members of every table class.
  \end{itemize}
\end{enumerate}

The independent evaluator \texttt{lvl} (used for the hits) received its own rules, written without copying those of \texttt{ivl}:

\begin{itemize}
\tightlist
\item
  exact powers through their vectors of prime exponents;
\item
  first-order anchors at exact powers;
\item
  exact offsets;
\item
  the equality of two writings of the same exact power.
\end{itemize}

It only uses the direction ``same form ⇒ same value''.

For n = 9 in any order three more rules were added to \texttt{lvl}, after some hits of the laboratory run stayed undecided. All of them were differences \textbar a − b\textbar{} that are exact zeros through identities between powers, such as \({4}^{0.9}\) = \({8}^{0.6}\) = \({2}^{9/5}\) and \({8}^{0.3}\) = \({0.5}^{-0.9}\); raising the precision can never decide an exact zero. The rules are:

\begin{itemize}
\tightlist
\item
  \textbf{L1, powers of canonical towers.} A value \(P^{±E}\), with P an exact rational and E an exact product of prime powers with rational exponents, raised to another such product, is computed exactly, with the exponent vectors added; (\(u^{a}\))\^{}b = \(u^{ab}\) for u \textgreater{} 0.
\item
  \textbf{L2, exact linear forms.} A value may carry the exact form \(c_{0}\) + Σ a\_\allowbreak{}\allowbreak{}k \(P_{k}\), with rational coefficients and products of prime powers with exponents in (0, 1). Sums, differences and products combine the forms exactly. When every term cancels, the value is the rational \(c_{0}\); a form that does not cancel claims nothing.
\item
  \textbf{H, exact forms of table members.} \texttt{lvl} evaluates the writing given by motor B for each member of a table, and attaches the form obtained (L1 or L2) when its enclosure meets motor B's ball of the member. A difference in the writing is accepted only when \texttt{lvl} proves it positive.
\end{itemize}

Each rule is used only in the direction ``identity ⇒ exact value''. None removes a value or drops a derivation, and each has its own randomized test against mpmath. With them, the 200 derivations left undecided by the laboratory are decided: 160 have no value (they reach 0 to a negative power) and 40 are exactly 0.

\textbf{Later changes (14--18 November).} The largest splits of n = 9 needed more changes; L1, L2 and H above are items 38--40 of \texttt{motorx\_\allowbreak{}\allowbreak{}changes.md}, and these are items 41--50 and the corrections of error 8. The changes to the query return the same hits or a smaller superset:

\begin{itemize}
\tightlist
\item
  a base exactly 1, against a window whose second-level constraint fixes the sign of ln w, has an empty preimage. This is the rule that already removed such pairs at the end, applied before the enumeration;
\item
  the sub-hits of a recursive query, and its distinct windows, are handled in bounded batches. Each window receives the same sub-hits; only their order changes;
\item
  for \textbar s − l\textbar{} = u with s exactly 1 and u tiny, each branch l = 1 ± u receives a second-level constraint on ln\textbar ln l\textbar, from u(1 − u/2) ≤ log1p(u) ≤ u and u ≤ −log1p(−u) ≤ u/(1 − u);
\item
  for each exact member x of s, the exact members y of l with \textbar x − y\textbar{} = u are looked up in the window \textbar ln y − ln x\textbar{} ≤ 2u/x instead of in the whole enclosure of l.
\end{itemize}

They were checked by equality of hits with the earlier code: n = 8 in any order (381 of 381 pairs) and splits 9, 13, 70 and 129 of n = 9.

The other changes are rules of \texttt{lvl}. Each is an identity between exact forms, used only in the direction ``identity ⇒ exact value'', with its own test against mpmath or exact rationals:

\begin{itemize}
\tightlist
\item
  X/(X ± c) and (X ± c)/X, with X a huge exact power, are anchors at 1 (the value 1 ∓ c/X, up to a bounded factor).
\item
  A huge exact rational q close to a small one is written exactly as M(1 + d), so that its integer powers are anchored.
\item
  A linear form that reduces to one term a·P is the exact power a·P, and two proportional forms have a rational quotient.
\item
  Linear forms may have exact constants of up to 600 000 bits. Huge exact powers are kept inside them as vectors of prime exponents. Integer powers k ≤ 8 of a form are expanded term by term, so that, for example, the leading term of \(3^{10^{7}}\) − (\(9^{10^{6}}\) + c)\(^{5}\) cancels exactly.
\item
  (Patches 7 and 7c.) The rules are:

  \begin{itemize}
  \tightlist
  \item
    the correction of error 8;
  \item
    a check of every exact form against the enclosure of its value, which leaves the value undecided when they disagree;
  \item
    bounds \(e^{min ln}\)·\textbar d\textbar\_\allowbreak{}\allowbreak{}min·\(e^{-|d|_max}\) ≤ \textbar A − B\textbar{} ≤ \(e^{max ln}\)·\textbar d\textbar\_\allowbreak{}\allowbreak{}max, with d = ln A − ln B, used only to compare with a window;
  \item
    differences, and roots of ``power + constant'', anchored at the exact power;
  \item
    exact offsets over towers with irrational exponents;
  \item
    exact forms kept through x·1, x/1 and 1/x.
  \end{itemize}
\end{itemize}

Five versions of \texttt{lvl} and \texttt{hitcheck} checked the records used for n = 9 in any order: the first version (1 303 of the 5 949 pieces, splits 0--16 of the laboratory; these records predate the fingerprint field), and patches 4 (95 pieces), 5 (1 086), 7 (3 459) and 7c (6). Each patch adds rules of the kind described above or corrects error 8, and each was tested on its own. Error 8 was introduced in patch 6, and no record made with patch 6 is used. The repository gives every version with its fingerprint, and the version that checked every record.

\textbf{Errors found and corrected.} As for the certifier (Appendix A.8), the work found errors, here eight. Errors 1--7 were corrected before the runs of Table 6, and every closed cell was run again with the corrected code; none changed a result. Error 8 was found during the run of n = 9 in any order, and everything evaluated with the faulty version was recomputed.

\begin{enumerate}
\def\labelenumi{\arabic{enumi}.}
\tightlist
\item
  \textbf{Third level of a power whose exponent is below \(e^{-10^{308}}\)} (C.4.5(4)). The bounds ν ≥ pred(\(μ_{y}\)) and ln\textbar ln u\textbar{} ≤ −DBL\_\allowbreak{}\allowbreak{}MAX do not follow from ν = ln(\(e^{μ_{y}}\) − \(μ_{x}\)). The first fails when \textbar{}\(μ_{x}\)\textbar{} \textgreater{} \(10^{295}\), the second when \(μ_{y}\) \textless{} ln DBL\_\allowbreak{}\allowbreak{}MAX ≈ 709.78. In the closed cells the rule had been applied 20 242 times, never with an invalid old bound.
\item
  \textbf{Three-term sums in the anchors of C.4.6(7).} They were plain double sums, which may round inward by one ulp, and are now outward. In 249 284 applications, compared in exact rational arithmetic, no old bound was invalid.
\item
  \textbf{The key used to merge equal enclosures.} It omitted the third level and compared the anchor rationals as doubles. Two different enclosures could therefore be merged, and the surviving one would attribute its statements to the values of the other. The key is now the bytes of every field.
\item
  \textbf{The exact zero in recursive V queries.} Windows {[}−∞, h{]} were clipped to the range and lost 0 ∈ V(S), which Theorem X needs (u + 0 = u).
\item
  \textbf{Undecided operations of motor B.} Motor B skips an operation it cannot decide and records it, and its answer is then no longer a superset; motor X did not look at that record. Answers with a nonempty record are now rejected. The record had been empty in every run, so no result changes.
\item
  \textbf{The sign of a near-1 member} was taken from a rounded bound, which is wrong for 0 \textless{} ln v \textless{} \(2^{-1075}\). It is now taken from the Arb ball. This was corrected before those members were first used.
\item
  \textbf{The second-level constraint passed to a sub-query} was stored in the same variable as the constraint of the query itself, so true hits could be dropped. It was found through an index error and corrected; every cell run with that code was run again.
\item
  \textbf{The sign of a linear form attached to a huge exact rational} (\texttt{lvl}, patch 6, 17 November).

  \begin{itemize}
  \tightlist
  \item
    \emph{The error.} This version began to attach linear forms to exact rationals of more than 600 000 bits, but in \textbar a − b\textbar{} it corrected the sign of the form only when the result was not exact. So \textbar c − P\textbar, with P a huge exact power and c \textless{} P, received the form c − P.
  \item
    \emph{The risk.} Where a value is taken from its form, a value inside the window could be declared outside without leaving a trace.
  \item
    \emph{The negative control.} It finds the wrong sign in 50 of 146 random forms with patch 6, and in none with the corrected version.
  \item
    \emph{Why the earlier versions are not affected.} They attach forms only to results that are not exact, and there the sign was corrected.
  \item
    \emph{What was done.} Every record made from 13 November on carries the fingerprint of the evaluator that checked it, and every split evaluated with patch 6 was recomputed with the corrected version.
  \end{itemize}
\end{enumerate}

\textbf{Changes in the hits of n = 9 in order.} The first version gave 22 hits per constant and the corrected one gives 40. The difference has three parts:

\begin{itemize}
\tightlist
\item
  \textbf{+24, split \(S_{1}\) = \{7, 8, 9\}:} the old merge (error 3) had fused four entries of N̂(23456) into one. The hits are of the form \(l^{-s}\) with l = \textbar1 − w\textbar, and the entries w differ only in the third level.
\item
  \textbf{−2, split \(S_{1}\) = \{6, 7, 8, 9\}:} the new exact rules for powers and offsets (C.4.6) now decide these differences of two huge entries.
\item
  \textbf{−4, split \(S_{1}\) = \{8, 9\}:} the exact treatment of bases near 1 removes these.
\end{itemize}

All 40 hits of each constant are proved outside the window by \texttt{lvl}.

\end{document}